%% file: main.tex
\documentclass[a4paper]{article}
\pdfoutput=1

\usepackage{lg_style}

\usepgfplotslibrary{groupplots}

\pgfplotsset{compat=newest}
\newcommand\figfontsize\footnotesize
\pgfplotsset{every axis legend/.append style={font=\figfontsize}}
\pgfplotsset{every axis title/.append style={font=\figfontsize, yshift=-1.5ex}}
\pgfplotsset{every x tick label/.append style={font=\figfontsize, yshift=0.5ex}}
\pgfplotsset{every y tick label/.append style={font=\figfontsize, xshift=0.5ex}}
\pgfplotsset{every axis x label/.append style={font=\figfontsize, yshift=1ex}}
\pgfplotsset{every axis y label/.append style={font=\figfontsize, yshift=-1.5ex}}

\newcommand\sqrtm[1]{#1^{\frac{1}{2}}}

\newcommand{\kld}[2]{\rmD_{\mathrm{KL}}\left(#1\parallel#2\right)}

\newcommand{\ldd}[2]{\rmD_{\ell d}\left(#1, #2\right)}
\newcommand{\lddt}[2]{\rmD_{\ell d} (#1, #2)}

\newcommand{\md}[3]{\rmD_{#1}\left(#2, #3\right)}
\newcommand{\mdt}[3]{\rmD_{#1} (#2, #3)}

\newcommand{\Grn}{\Gr(r,n)}
\newcommand{\Stn}{\St(r,n)}

\newcommand{\glr}{\mathrm{GL}_r}

\newcommand\x{X}
\newcommand\e{E}
\newcommand\ma{m_A}
\newcommand\mx{m_\x}
\newcommand\me{m_\e}
\newcommand\my{m_Y}
\newcommand\mw{m_W}
\newcommand\ca{C_A}
\newcommand\cx{C_\x}
\newcommand\ce{C_\e}
\newcommand\cy{C_Y}
\newcommand\cw{C_W}
\newcommand\cax{C_{A\x}}

\newcommand\ms{m_\star}
\newcommand\cs{C_\star}
\newcommand\gs{G_\star}
\newcommand\hs{h_\star}

\newcommand\mv{m_V}
\newcommand\cv{C_V}
\newcommand\gv{G_V}
\newcommand\hv{h_V}

\newcommand\revis[1]{{\color{black}#1\color{black}}}

\title{Optimal projection of observations\\in a Bayesian setting}
\date{}
\author{L. Giraldi, O. P. Le Ma\^{\i}tre, I. Hoteit, O. M. Knio}

\begin{document}
\maketitle

\begin{abstract}
    \input{abstract.tex}
\end{abstract}

\input{intro.tex}

\input{linear.tex}

\input{functionals.tex}

\input{numerical.tex}

\input{nonlinear.tex}

\input{conclusions.tex}

\input{ack.tex}

\input{proofs.tex}

\bibliography{biblio}{}
\bibliographystyle{plain}

\end{document}

%% file: abstract.tex
%!TEX root = main.tex
Optimal dimensionality reduction methods are proposed for
% This work proposes optimal dimensionality reduction methods adapted to
the Bayesian inference of a
Gaussian linear model with additive noise in presence of overabundant data. Three different optimal
projections of the observations are proposed based on information theory: the projection that
minimizes the Kullback-Leibler divergence between the posterior distributions of the original and
the projected models, the one that minimizes the expected Kullback-Leibler divergence between the
same distributions, and the one that maximizes the mutual information between the parameter of
interest and the projected observations. The first two optimization problems are formulated as the
determination of an optimal subspace and therefore the solution is computed using Riemannian
optimization algorithms on the Grassmann manifold. Regarding the maximization of the mutual
information, it is shown that there exists an optimal subspace that minimizes the entropy of the posterior
distribution of the reduced model; a basis of the subspace can be computed as the solution to a
generalized eigenvalue problem; an a priori error estimate on the mutual information is available
for this particular solution; and that the dimensionality of the subspace to exactly conserve the
mutual information between the input and the output of the models is less than the number of
parameters to be inferred. Numerical applications to linear
and nonlinear models are used to assess the efficiency of the proposed approaches, and to highlight their
advantages compared to standard approaches based on the principal component analysis of the observations.

% {\bf Keywords:} Optimal dimensionality reduction, Optimal data reduction, Gaussian linear model, Information theory

%% file: intro.tex
%!TEX root = main.tex

\section{Introduction}
\label{sec:introduction}

We consider the problem of Bayesian inference in the case of overabundant data. The goal is to
compute an optimal approximation of the posterior distribution by projection of the observations.
These projections are computed solving the following optimization problem
\begin{equation*}
    \min_V \cJ\left(P(X\mid Y=y) , P(X\mid W = V^T y)\right),
\end{equation*}
where $X$ is the inferred parameter, $Y$ the observations of the Bayesian model with values in
$\bR^n$, $y\in\bR^n$ the data, $W = V^T Y$ the reduced observations with values in $\bR^r$,
$V\in\bR^{n\times r}$ the deterministic matrix defining the projection from the full to the reduced
observations, and $\cJ$ a functional defining the optimality criterion.

In the literature, the most popular dimensionality reduction techniques result from the optimal
approximation of the observations $Y$ with respect to the $L^2$ norm defined by
\begin{equation*}
    \norm{Y}_{L^2}^2 = \bE(Y^T Y),
\end{equation*}
with $\bE$ the expectation operator.
A best low-rank approximation of $Y$ with respect to the $L^2$ norm is computed by a singular
value decomposition. Then, with $V$ being the matrix composed of the dominant left eigenvectors of
$Y$, the approximation is given by $Y \approx V W$ where $W = V^T Y$. The reader could refer to
\cite[Section 4.4.3]{Hackbusch2012} for the presentation of the singular value decomposition in a
general case.  When applied to the centered random vector $Y-\bE(Y)$, this decomposition
is also called truncated Karhunen-Lo\`eve expansion~\cite{karhunen1947,LeMaitre2010,loeve1948}, or
principal component analysis~\cite{Hotelling1933,Jolliffe2002,Pearson1901}.

Let $\Omega$ denote a sample space. When considering the random vector as a map
$\Omega\to \bR^n$, an approximation of the observations $Y$ can be used to define the matrix $V$
based on the $L^\infty$ norm defined by
\begin{equation*}
    \norm{Y}_{L^\infty} = \sup_{\omega \in \Omega} \max_{1 \le i \le n} |Y_i(\omega)|.
\end{equation*}
The empirical interpolation method~\cite{BARRAULT2004667} provides an approximation of the form
$\omega\mapsto VW(\omega)$ of a parametric vector $\omega \mapsto Y(\omega)$ based on this supremum
norm. Given that computing the supremum is not tractable, the sample space is restricted to a
finite sample $\Omega_N = \{\omega_i\}_{i=1}^N \subset \Omega$ and the norm is approximated by
\begin{equation*}
    \norm{Y}_{L^\infty} \approx \max_{\omega \in \Omega_N} \max_{1 \le i \le n} |Y_i(\omega)|.
\end{equation*}
The interpolation is then defined to be exact on a subsample $\{\omega_j^\star\}_{j=1}^r \subset
\Omega_N$, and on a restricted number of indices $\{i_j\}_{j=1}^r \subset \{1,\hdots,n\}$. The
interpolation points $((\omega_j^\star, i_j))_{j=1}^r$ are selected in a greedy fashion using the
supremum norm of the error, such that
\begin{equation*}
    (\omega_r^\star, i_r) = \arg \max_{\omega\in\Omega_N}\max_{1 \le i \le n} \left|Y_i -
\sum_{j=1}^{r-1}V_{ij} W_j\right|.
\end{equation*}
The approach reduces to the interpolation of the matrix $M_{ij} = Y_i(\omega_j)$ and is therefore
very close to the cross approximation method~\cite{Bebendorf2000} for the low-rank interpolation of
a matrix. A detailed comparison between the singular value decomposition, the empirical
interpolation method and the cross approximation is provided in~\cite{Bebendorf2014}. A weighted
variant of the empirical interpolation method was introduced in~\cite{chen_quarteroni_rozza_2014}
in order to take into account probability measures. Given a positive weight $w:\Omega\to\bR_+$, the
supremum norm is modified such that
\begin{equation*}
    \norm{Y}^w_{L^\infty} = \sup_{\omega \in \Omega} \max_{1 \le i \le n} w(\omega)|Y_i(\omega)|,
\end{equation*}
yielding different interpolation points. While the approximations based on the $L^2$ and $L^\infty$ norms
are widely used, they are optimal with respect to the output of the model and are not directly
related to the distribution of the parameter to infer. An exception to this rule
is~\cite{Giraldi2017} where a weighted singular value decomposition was used to accommodate uniform
priors.

Geppert et al.~\cite{Geppert2017} propose to reduce the number of observations in a Bayesian
regression framework using random projections. The methodology relies on $\varepsilon$-subspace
embeddings that are particular maps of the form $\Pi\in\bR^{r\times n}$ satisfying
\begin{equation*}
    (1-\varepsilon) \norm{B x}^2 \le \norm{\Pi B x}^2 \le (1 + \varepsilon) \norm{B x}^2,
\end{equation*}
for a particular $B \in \bR^{n\times q}$, and for any $x\in\bR^q$ with a probability $1-\alpha$. In
order to obtain an error on the posterior distribution of order $\varepsilon$ in terms of the
Wasserstein metric, it is shown that the number of observations required is $\cO((q +
\log(1/\alpha))/\varepsilon^2)$, $ \cO(q\log(q/\alpha)/\varepsilon^2)$ or
$\cO(q^2/(\alpha\varepsilon^2))$, depending on the embedding. Even though the dimension $r$ can be
drastically smaller than the original number of observations, $n$, it can still be relatively large for a
small risk $\alpha$ and a small error $\varepsilon$.

Another related technique is introduced in a series of papers~\cite{cui2014, Spantini2017,
Spantini2015}. Given a Gaussian linear model, the goal is to directly compute an approximate
posterior covariance matrix as a low-rank update of the prior covariance. Given a particular loss
function depending only on the covariance, it is shown that an optimal low-rank update can be
derived from a generalized eigenvalue problem. The resulting distribution is then optimal in terms
of the Hellinger distance and Kullback-Leibler divergence under the assumption that the mean is exactly
recovered.  An optimal mean is also derived as a linear projection of the data by minimizing the
Bayes risk defined as the expected Mahalanobis distance between the parameter of interest $X$ and
the approximate mean. Regarding this methodology, two disadvantages are notable: there is an
inconsistency between the optimality criteria of the mean and the covariance, and the computation
of the different matrices requires the inversion of the covariance matrix of the noise, which is of
large dimension $n\times n$.

In this work, the approximate mean and covariance are defined similarly, namely as an affine function of
the data and as a low-rank update of the prior covariance, respectively. However, they result from
the optimal projections of the statistical model using criteria from information theory, namely
the Kullback-Leibler divergence, the expected Kullback-Leibler divergence, and the mutual
information, which is the first contribution of this paper. The second contribution concerns the
choice of the practical numerical algorithm for the minimization of the expected
Kullback-Leibler divergence between the posterior distributions of the full and reduced model. It
is moreover shown that a solution to the optimization problem defined as the maximization of the
mutual information $\cI(X,W)$ between the parameter $X$ and the reduced observations $W$ is given
by the solution of a generalized eigenvalue problem that does not require the inversion of a large
matrix. We can moreover estimate the loss $\cI(X,Y)-\cI(X,W)$ and show that no more than $q$
projections are required to recover the full mutual information, where $q$ is the size of $X$. The
last contribution of this work concerns the illustration of the method on linear and nonlinear examples.

This paper is organized as follows. In Section~\ref{sec:linear_model}, the full and reduced linear
models are presented in the Gaussian case, as well as other required definitions.
The posterior distributions are then provided explicitly in closed form.
Section~\ref{sec:optimality_criteria} introduces the three different optimization problems that are
used to define the alternative optimal projections of the
observations.
The analysis of the corresponding optimal subspace and the numerical algorithms for their computation are then provided.
The methodologies are finally applied and illustrated to a Bayesian linear regression problem in
Section~\ref{sec:bayesian_linear_regression} and to a nonlinear problem in Section~\ref{sec:application_to_nonlinear_problems}.
Major conclusions are summarized
in Section~\ref{sec:conclusions}.

%% file: linear.tex
%!TEX root = main.tex

\section{Linear Gaussian model}
\label{sec:linear_model}

\subsection{Models}
\label{sub:models}
We consider an abstract probability space $(\Omega, \cF, \bP)$, where $\Omega$ is the sample space,
$\cF$ is a $\sigma$-algebra and $\bP$ a probability measure. Given an $\bR^n$-valued random vector
$Z$, we denote by $P(Z)$ the pushforward probability measure such that $P(A) = \bP(Z^{-1}(A))$ for
any set $A$ in the Borel algebra of $\bR^n$, and $f_Z$ the probability density function defined with
respect to the Lebesgue measure.

We consider the following linear model
\begin{equation}
  \label{eq:linearmodel}
  Y = BX + \e,
\end{equation}
where $B\in\bR^{n\times q}$ is the design matrix, $X$ is the random parameter we want to infer and
$\e$ is the random noise. The random vector $X$ (resp.\ $\e$) is supposed to follow the multivariate
normal distribution $\cN(\mx,\cx)$ (resp.\ $\cN(\me, \ce)$) with mean $\mx\in\bR^q$ (resp.\
$\me\in\bR^n$) and covariance $\cx\in\bR^{q\times q}$ (resp.\ $\ce\in\bR^{n\times n}$). The input
parameter $X$ and the noise $\e$ are assumed to be independent.

In order to compress the amount of data used for the inference, we introduce $V = (v_i)_{i=1}^r \in
\bR^{n\times r}$, a reduced basis of observations. In the following, the term reduced space may
be used for $V$, as we look for the projection of the observations on the space spanned by the
columns of $V$. The linear model expressed in the reduced coordinates is therefore
\begin{equation}
  \label{eq:reducedlinearmodel}
  W = V^T BX + V^T \e,
\end{equation}
and the reduction is efficient if, for $r \ll n$, the posterior distribution $P(X\mid W)$ is
close to $P(X\mid Y)$ in some sense defined in Section~\ref{sec:optimality_criteria}. Our main goal
is to compute a suitable matrix $V$ which satisfies this condition.

In order to subsequently apply the different methodologies to nonlinear models of the form
\begin{equation*}
    Y= A(X) + \e,
\end{equation*}
it is beneficial to consider the random vector $A(X)=BX$. In the linear case, this random vector
follows the distribution $\cN(\ma, \ca)$ where
\begin{equation*}
    \ma = B\mx \quad \text{and} \quad \ca = B\cx B^T.
\end{equation*}
We also denote by $\cax$ the covariance between $A$ and $X$, that is
\begin{equation*}
    \cax = \bE\left(\left(A(X) - \ma\right)\left(X - \mx\right)^T\right) = B\cx.
\end{equation*}
Given the structure of the problem, the random vectors $Y$ and $W$ are also distributed according to the
multivariate normal distribution $\cN(\my, \cy)$ and $\cN(\mw,\cw)$, respectively, with $\my = \ma +
\me$, $\cy = \ca + \ce$, $\mw = V^T\my$, and $\cw = V^T\cy V$. The extension of the approach
to nonlinear problems is based on the three quantities $\ma$, $\ca$ and $\cax$.

\subsection{Posterior distributions}
\label{sub:posterior_distributions}

Given the linear Gaussian structure of Equation~\eqref{eq:linearmodel}, an observation $y$ of $Y$,
and a reduced basis $V$, the posterior distributions $P(X\mid Y=y)$ and $P(X\mid W=V^T y)$ can
be analytically derived. The result is summarized in Proposition~\ref{prop:posteriordistrib}.

\begin{proposition}
\label{prop:posteriordistrib}
    The posterior distribution $P(X\mid Y=y)$ follows the multivariate normal distribution
    $\cN(\ms, \cs)$, where
    \begin{equation}
        \label{eq:covfullposterior}
        \cs = \cx(\cx + \cax^T\ce^{-1}\cax)^{-1} \cx = \cx - \cax^T\cy^{-1}\cax,
    \end{equation}
    and
    \begin{equation}
        \label{eq:meanfullposterior}
        \ms = \gs (y - \my) + \hs,
    \end{equation}
    with $\gs = \cax^T\cy^{-1}$ and $\hs = \cs \cx^{-1}\mx + \gs\ma$.

    \bigskip
    Regarding the posterior distribution of the reduced model, if the matrix $V\in\bR^{n\times r}$
    is full-rank, the distribution $P(X\mid W=V^T y)$ follows the multivariate normal
    distribution $\cN(\mv, \cv)$, where
    \begin{align}
    \label{eq:covreducedposterior}
    \cv &= \cx\left(\cx + \cax^T V \left(V^T\ce V\right)^{-1}V^T \cax\right)^{-1}\cx \\
        & = \cx - \cax^T V\left(V^T\cy V\right)^{-1}V^T\cax, \notag
    \end{align}
    and
    \begin{equation}
    \label{eq:meanreducedposterior}
    \mv = \gv V^T(y - \my) + \hv,
    \end{equation}
    with
    \begin{equation*}
        \gv = \cax^T V (V^T \cy V)^{-1},
    \end{equation*}
    and
    \begin{equation*}
        \hv = \cv \cx^{-1}\mx + \gv V^T\ma.
    \end{equation*}
\end{proposition}
\begin{proof}
    See Appendix~\ref{proof:posteriordistrib}.
\end{proof}

Regarding Proposition~\ref{prop:posteriordistrib}, we can first notice that the two expressions
\begin{align*}
    \cs &= \cx (\cx + \cax^T\ce^{-1}\cax)^{-1}\cx \\
    \text{and} \quad \cv &= \cx\left(\cx + \cax^T V \left(V^T\ce V\right)^{-1}V^T
    \cax\right)^{-1}\cx
\end{align*}
show that the matrices $\cs$ and $C_V$ are always symmetric positive definite, even for a nonlinear
model. In the following, we denote by $\glr$ the set of invertible matrices in $\bR^{r\times r}$,
and obtain an invariance property expressed in Proposition~\ref{prop:invariancemeancov}.
\begin{proposition}
    \label{prop:invariancemeancov}
    For all matrices $M\in\glr$, we have
    \begin{equation*}
        m_{VM} = \mv \quad \text{and} \quad C_{VM} = \cv.
    \end{equation*}
    Therefore, the posterior distribution $P(X\mid W=V^T y)\sim \cN(\mv, \cv)$ is invariant under
    invertible linear transformation of the matrix $V$ on the right.
\end{proposition}
\begin{proof}
    See Appendix~\ref{proof:invariancemeancov}.
\end{proof}
In practice, this proposition means that $V$ is less important than $\ran(V)$ in the determination
of the posterior distribution. Indeed, rescaling, rotating or permuting the observations in
Equation~\eqref{eq:reducedlinearmodel} does not affect the posterior distribution $P(X\mid W=V^T
y)$.

Formally, the Grassmann manifold $\Grn$ defined as the set of $r$ dimensional subspace of $\bR^n$ is
therefore the set of interest to determine the optimal reduced observations. In this work, we
identify $\Grn$ with the quotient manifold $\bR^{n\times r}_* / \glr$ following~\cite{Absil2004},
where $\bR^{n\times r}_*$ is the set of full rank matrices of $\bR^{n\times r}$ and the quotient
space is defined by
\begin{equation*}
    \Grn = \bR^{n\times r}_*/\glr = \left\{ [V];\ V \in \bR^{n\times r}_*\right\}, \text{ where }
    [V] = \left\{VM;\ M \in \glr\right\}.
\end{equation*}
Finally, Proposition~\eqref{prop:invariancemeancov} means that it is more important to identify the
equivalence class $[V]\in\bR^{n\times r}_*$ than a particular matrix $V\in\bR^{n\times r}_*$. 

The next section presents the different proposed optimization problems, where the Grassmann
manifold $\Grn$ has an important role.

%% file: functionals.tex
%!TEX root = main.tex

\section{Optimality criteria for the definition of the reduced basis}
\label{sec:optimality_criteria}

\subsection{Kullback-Leibler divergence minimization}
\label{sub:kullback_leibler_minimization}

Given two distributions $P(Z_0)$ and $P(Z_1)$, the Kullback-Leibler divergence between them is
defined by
\begin{equation}
    \label{eq:defkld}
    \kld{P(Z_0)}{P(Z_1)} = \bE_{Z_0}\left(\log \frac{f_{Z_0}}{f_{Z_1}}\right).
\end{equation}
This divergence quantifies the ``information lost when [$P(Z_1)$] is used to approximate
[$P(Z_0)$]'' according to~\cite[Section 2.1]{burnham2003model}. The Kullback-Leibler divergence is
always positive and null if and only if the two distributions are identical, therefore
defining a generalized distance between distributions.

This interpretation of the Kullback-Leibler divergence leads us to consider the following
functional $\sJ_0:\bR^{n\times r}_* \to \bR$ defined by
\begin{equation*}
    \sJ_0(V) = \kld{P(X\mid Y=y)}{P(X\mid W=V^T y)}.
\end{equation*}
The domain definition of $\sJ_0$ must be restricted to the set of full rank matrices $\bR^{n\times
r}_*$ in order to comply with Proposition~\ref{prop:posteriordistrib} characterizing the posterior
distributions. Given that we are working with Gaussian distributions, the computation of the
Kullback-Leibler divergence is always well-posed (i.e. $f_{Z_1}$ is always stricly positive in
Equation~\eqref{eq:defkld}). The general expression of the Kullback-Leibler divergence between two
Gaussian distribution is given in Proposition~\ref{prop:kld}.
\begin{proposition}
    \label{prop:kld}
    Assuming that $Z_0\sim \cN(m_0, C_0)$ and $Z_1 \sim \cN(m_1, C_1)$ are $\bR^q$-valued random
    variables, the Kullback-Leibler divergence between $P(Z_0)$ and $P(Z_1)$ is expressed by
    \begin{equation*}
        \kld{P(Z_0)}{P(Z_1)} = \frac{1}{2} \left(\ldd{C_0}{C_1} + \md{C_1}{m_0}{m_1}\right),
    \end{equation*}
    where $\lddt{C_0}{C_1}$ is the Bregman $\log\det$ divergence between $C_0$ and $C_1$ defined by
    \begin{equation*}
        \ldd{C_0}{C_1} = \trace\left(C_0C_1^{-1}\right) - \log\det\left(C_0C_1^{-1}\right) - q,
    \end{equation*}
    and $\mdt{C_1}{m_0}{m_1}$ is the Mahalanobis divergence defined by
    \begin{equation*}
        \md{C_1}{m_0}{m_1} = (m_0 - m_1)^T C_1^{-1}(m_0 - m_1).
    \end{equation*}
\end{proposition}
\begin{proof}
    See Appendix~\ref{proof:kld}.
\end{proof}
As a consequence of Proposition~\ref{prop:kld}, the functional $\sJ_0$ has a closed form depending on $\ms$, $\cs$, $\mv$ and $\cv$:
\begin{align}
    \label{eq:kld_as_fun_of_ldd_md}
    \sJ_0(V) &=\kld{P(X\mid Y=y)}{P(X\mid W=V^T y)} \\
             &= \frac{1}{2} \left(\ldd{\cs}{\cv} + \md{\cv}{\ms}{\mv}\right) \notag\\
    &=\frac{1}{2} \left(\trace\left(\cs\cv^{-1}\right) - \log\det\left(\cs\cv^{-1}\right)-q\right.\notag\\ 
    &\qquad\qquad\qquad+ \left.\left(\ms-\mv\right)^T \cv^{-1}\left(\ms - \mv\right)\right).\notag
\end{align}
Given Proposition~\ref{prop:invariancemeancov}, for all $M\in\glr$, we have $\sJ_0(VM)=\sJ_0(V)$.
It means that we are in fact interested in the map defined on $\Grn$ by $[V] \mapsto \sJ_0(V)$.
The minimization problem of interest is therefore
\begin{equation}
    \label{eq:min_kld}
    \min_{[V]\in\Grn} \kld{P(X\mid Y=y)}{P(X\mid W=V^T y)}.
\end{equation}
We can show that the following result holds.
\begin{theorem}
    \label{th:existenceminkld}
    There exists a solution to Problem~\eqref{eq:min_kld}.
\end{theorem}
\begin{proof}
    See Appendix~\ref{proof:existenceminkld}.
\end{proof}

Note that the minimization of the Kullback-Leibler divergence in Problem~\eqref{eq:min_kld} results
in an \emph{a posteriori} reduction in the sense that a realization $y$ of $Y$ is required to evaluate the
cost function. In the following, other functionals are proposed that circumvent this issue.

\subsection{Expected Kullback-Leibler divergence minimization}
\label{sub:expected_kullback_leibler_divergence_minimization}

The first possibility to remove the dependence on the data is to work on the expected
Kullback-Leibler divergence with respect to the observation, where the measurement $Y$ is treated
as a random variable.
Similarly to Section~\ref{sub:kullback_leibler_minimization}, let $\sJ_1:\bR^{n\times r}_*\to \bR$ be defined by
\begin{equation*}
    \sJ_1(V) = \bE_Y \left(\kld{P(X\mid Y)}{P(X\mid W=V^T Y)}\right).
\end{equation*}
The expected Kullback-Leibler divergence admits a closed form as well,
presented in the next proposition.
\begin{proposition}
    \label{prop:expkld}
    We have the following equality
    \begin{equation}
        \label{eq:defexpkld}
        \sJ_1(V)= \frac{1}{2}\left(\ldd{\cs}{\cv} + \bE_Y\left(\md{\cv}{\ms}{\mv}\right)\right),
    \end{equation}
    where
    \begin{multline*}
 \bE_Y\left(\md{\cv}{\ms}{\mv}\right) =\\  \trace \left(\cv^{-1}\left(\gs -\gv V^T\right)\cy\left(\gs-\gv V^T\right)^T\right) + (\hs - \hv)\cv^{-1}(\hs - \hv).
    \end{multline*}
\end{proposition}
\begin{proof}
    See Appendix~\ref{proof:expkld}.
\end{proof}
Using Proposition~\ref{prop:invariancemeancov} and Equation~\eqref{eq:defexpkld}, we can show that
$\sJ_1(V) = \sJ_1(VM)$ for all matrices $M\in\glr$. We are therefore interested in finding the
optimal equivalence class $[V]$ and solving the minimization problem
\begin{equation}
    \label{eq:min_expect_kld}
    \min_{[V]\in\Grn}  \bE_Y \left(\kld{P(X\mid Y)}{P(X\mid W=V^T Y)}\right).
\end{equation}
As in Section~\ref{sub:kullback_leibler_minimization}, we can prove the following result.
\begin{theorem}
    \label{th:existenceminexpkld}
    There exists a solution to Problem~\eqref{eq:min_expect_kld}.
\end{theorem}
\begin{proof}
    The proof is similar to the one in Appendix~\ref{proof:existenceminkld}, replacing $\sJ_0$ by
    $\sJ_1$.
\end{proof}

\begin{remark}
    The minimization of the log det divergence $\lddt{\cs}{\cv}$ has also been considered, being
    the data-free part of the Kullback-Leibler divergence. It has been ignored in the paper as it
    did not bring additional insights on the optimal construction of the reduced observations.
\end{remark}

\subsection{Mutual information maximization and entropy minimization}
\label{sub:entropy_minimization}

In this section the Shannon entropy and the mutual information are introduced. 
The entropy $H(Z)$ (sometimes denoted $H(P(Z))$) quantifies the uncertainty or the amount of information contained in a random variable $Z\sim P(Z)$ and is defined by
\begin{equation*}
    H(Z) = \bE_Z(-\log(f_Z(Z))).
\end{equation*}
The mutual information $\cI(Z_0, Z_1)$ between the two random variables $Z_0\sim P_0(Z_0)$ and $Z_1\sim
P_1(Z_1)$ is a measure of the information that $Z_0$ contains about $Z_1$,
and is defined by
\begin{equation*}
    \cI(Z_0, Z_1) = H(Z_0) + H(Z_1) - H(Z_0, Z_1),
\end{equation*}
where $H(Z_0, Z_1)$ is the entropy of the \emph{joint distribution} of $Z=(Z_0,Z_1)$. 
From this definition, it is clear that the mutual information is symmetric.

The new definition of the reduced basis, introduced in this section, is related to the definition of the mutual
information. We would like the reduced observations $W$ to contain as much information as possible
about $X$. We therefore consider the following maximization problem
\begin{equation}
    \label{eq:maxmi}
    \max_{V\in\bR^{n\times r}_*} \cI(W, X).
\end{equation}
Note that another expression of the mutual information is 
\begin{equation*}
    \cI(W, X) = \bE_W(\kld{P(X\mid W)}{P(X)}),
\end{equation*}
showing that this strategy aims at maximizing the expected information gain between the prior and the
posterior distributions of $X$.

The optimization problem in Equation~\eqref{eq:maxmi} admits a simple
solution presented in Theorem~\ref{th:maxmi}. Moreover, we shall show that the maximization of the mutual
information is equivalent to the minimization of the entropy of the posterior distribution
$P(X\mid W=V^T y)$.
\begin{theorem}
\label{th:maxmi}
    The following equalities hold
    \begin{align*}
        \cI(W,X) &= \frac{1}{2}\log\det\left(\left(V^T \cy V\right)\left(V^T\ce
        V\right)^{-1}\right) \\
        \text{and} \quad H(P(X\mid W=V^T y)) &= -\cI(W, X) +\frac{1}{2}\log\det\cx +
        \frac{q}{2}\log(2\pi e).
    \end{align*}
    As a consequence, the maximization of the mutual information $\cI(W,X)$ and the minimization of
    the entropy of the posterior distribution $H(P(X\mid W=V^T y))$ with respect to $V$ admit the
    same solutions for any realization $y$ of $Y$. We have the equality
    \begin{equation}
        \label{eq:maxmi_result}
        \max_{V\in\bR^{n\times r}_*} \cI(W,X) = \frac{1}{2} \sum_{i=1}^r \log \lambda_i,
    \end{equation}
    where $(\lambda_i)_{i=1}^r$ are the $r$ dominant eigenvalues of the following generalized
    eigenvector problem
    \begin{equation}
        \label{eq:maxmi_eig}
        \cy v = \lambda \ce v, \quad \lambda \in \bR, \ v \in \bR^n.
    \end{equation}
    A solution to the optimization Problem~\eqref{eq:maxmi_result} is given by the matrix $V$ with
    columns being eigenvectors $(v_i)_{i=1}^r$ associated to the dominant eigenvalues of
    Problem~\eqref{eq:maxmi_eig}.
\end{theorem}
\begin{proof}
    See Appendix~\ref{proof:maxmi}.
\end{proof}
Several remarks follow this result. First, the map $\sJ_2:V\mapsto \cI(W,X)$ is also invariant under the
transformation $\sJ_2(VM) = \sJ_2(V)$ for any invertible matrix $M\in\glr$, and therefore the
solution should be searched in the Grassmann manifold. However in the present case, a particular 
solution admits a simple characterization.

The generalized eigenvalue problem in Equation~\eqref{eq:maxmi_eig} is used to define the
optimal mean in~\cite{Spantini2015} to minimize the Bayes risk. It is however unclear how this
optimal mean is related to the mean defined in Equation~\eqref{eq:meanreducedposterior}. Note
moreover that the computation of the optimal mean from~\cite{Spantini2015} requires the inversion of
the matrix $\ce$, which is not needed in the presently developed approach.

Another interesting feature  of Theorem~\ref{th:maxmi} is that it provides an a priori
estimate for the reduction error, based on the mutual information, summarized in the following Corollary~\ref{coro:maxmi}.
\begin{corollary}
    \label{coro:maxmi}
    Let $V\in\bR^{n\times r}$ be a particuler solution to Problem~\eqref{eq:maxmi_result} and $W$ be the
    reduced model associated to $V$, and $(\lambda_i)_{i=1}^n$ be the eigenvalues associated to
    Problem~\eqref{eq:maxmi_eig} sorted in a decreasing order. Then, the relative error on the mutual
    information is given by
    \begin{equation*}
       \frac{\cI(Y,X) - \cI(W, X)}{\cI(Y,X)} = 1 -
        \frac{\sum_{i=1}^r\log\lambda_i}{\sum_{i=1}^n\log\lambda_i}.
    \end{equation*}
\end{corollary}
In fact, the entire spectrum of $\cy$ is not required to estimate the error. In
practice, we only need to determine the eigenvalues $\nu_i$ associated to the following problem
\begin{equation}
    \label{eq:maxmi_eig_ca}
    \ca v = \nu \ce v, \quad \nu \ge 0, \ v\in\bR^n.
\end{equation}
If $\lambda$ is an eigenvalue associated to Problem~\eqref{eq:maxmi_eig}, then $\nu=\lambda - 1$ is
an eigenvalue associated to Problem~\eqref{eq:maxmi_eig_ca}. Considering Problem~\eqref{eq:maxmi_eig_ca} 
is beneficial in practice because $\ca$ is at most a rank-$q$ matrix. 
This remark leads to the following important result on the number of required projections to get the same mutual information between the observations and the parameter of interest, for the full and the reduced model.
\begin{corollary}
    \label{coro:maxmi_nu}
    Let $V\in\bR^{n\times r}$ be a solution to Problem~\eqref{eq:maxmi_result} and $W$ be the
    reduced model associated to $V$, and $(\lambda_i)_{i=1}^n$ be the eigenvalues associated to
    Problem~\eqref{eq:maxmi_eig}, sorted in a decreasing order. Let $(\nu_i)_{i=1}^n$ be the eigenvalues
    associated to Problem~\eqref{eq:maxmi_eig_ca} (i.e.\ $\lambda_i = 1 + \nu_i$), and let $m \le q
    \ll n$ be the rank of $B\in\bR^{n\times q}$ (see Equation~\eqref{eq:linearmodel}). Then
    $\ca=B\cx B^T$ is a rank-$m$ matrix, and the relative error on the mutual information is given by
    \begin{equation}
        \label{eq:apriorierrorestimate}
       \frac{\cI(Y,X) - \cI(W, X)}{\cI(Y,X)} = 1 -
       \frac{\sum_{i=1}^r\log(1 + \nu_i)}{\sum_{i=1}^m\log(1 + \nu_i)}.
   \end{equation}
   The condition $r\ge m$ implies that $\cI(W,X)=\cI(Y,X)$ and the mutual informations
   between the observations and the parameter of interest are the same for the full and the reduced
   model. In particular, the condition is satisfied for $m=q$.
\end{corollary}

As a side note, the principal component analysis of the observations $Y$ yields a reduced basis defined
as the dominant eigenvectors of $\cy$. Therefore, the resulting reduced space is optimal with
respect to the mutual information in the case of a white noise, i.e.\ $\ce = \sigma^2 \rmI_n$.
However, denoting $(\chi_i)_{i=1}^n$ the eigenvalues of $\cy$ sorted in a decreasing order, the
corresponding estimate of the relative reduction error on the mutual information is given by
\begin{equation*}
    \frac{\cI(Y,X) - \cI(W, X)}{\cI(Y,X)} = 1 -
    \frac{\sum_{i=1}^r\log\left(\frac{\chi_i}{\sigma^2}\right)}{\sum_{i=1}^n\log\left(\frac{\chi_i}{\sigma^2}\right)}.
\end{equation*}
Note that the usual error criteria used in the principal component analysis between the random variable
$Y$ and its rank-$r$ truncated version $Y_r$ controls the $L_2$ norm and is given by (see e.g.\
\cite[Proposition 2.1]{Bebendorf2014})
\begin{equation*}
    \frac{\bE\left(\norm{Y - Y_r}^2_2\right)}{\bE\left(\norm{Y}_2^2\right)} = 1 -
    \frac{\sum_{i=1}^r\chi_i^2}{\sum_{i=1}^n\chi_i^2}.
\end{equation*}

\subsection{Numerical solution to the optimization problems}
\label{sub:numerical_solution_of_the_optimization_problems}

For any functional $\sJ\in\{\sJ_0, \sJ_1, \sJ_2\}$ involved in the optimization problems
presented in Section~\ref{sec:optimality_criteria}, the following property holds
\begin{equation*}
    \sJ(V) = \sJ(VM), \quad \forall M\in\glr.
\end{equation*}
As a consequence of this invariance, there exists an infinite number of solutions to the
optimization problems and the Hessian of the functional $\nabla^2 \sJ$ is ill-conditioned in a
neighbourhood of a solution. The main consequence is that we cannot use a standard Newton algorithm
to solve these nonlinear problems without regularizing the optimization problem first.

In order to circumvent this issue, we consider here the restriction of the optimization problem to
the Grassmann manifold $\Grn$, replacing the search for a $n\times r$ matrix by the search of a
$r$-dimensional linear subspace of $\bR^{n\times n}$. In order to solve Problems~\eqref{eq:min_kld}
and~\eqref{eq:min_expect_kld} we choose to use a specific algorithm exploiting the smooth manifold
structure of $\Grn$, that is the Riemannian trust-region algorithm~\cite{Absil2007} implemented in
the \emph{Pymanopt} library~\cite{Townsend2016}. The derivatives of the cost functions are computed
by automatic differentiation with the \emph{autograd}
library~\cite{autograd}.

Given a finite dimensional vector space $\cV$ equipped with the inner product
$\innert{\cdot}{\cdot}_\cV$ and the associated norm $\normt{\cdot}_\cV$, the trust-region algorithm
consists in correcting the current iterate $V\in\cV$ with $W\in\cV$ using a quadratic approximation
of the functional $\sJ$. $W$ is
defined as the solution to
\begin{equation*}
    \min_{W\in\cV} m(W) = \sJ(V) + \inner{\nabla \sJ(V)}{W}_\cV +
    \frac{1}{2}\inner{\nabla^2\sJ(V)W}{W}_\cV,
\end{equation*}
such that $\normt{W}^2_\cV \le \Delta^2$,
where $\nabla \sJ$ (resp.\ $\nabla^2 \sJ$) is the gradient (resp.\ Hessian) of $\sJ$. The trust-region
radius $\Delta$ is adapted at each iteration of the algorithm based on the quantity
\begin{equation*}
    \rho = \frac{\sJ(V) - \sJ(V + W)}{m(0) - m(W)}.
\end{equation*}
If $\rho$ is close to 1, the quadratic approximation is good and the radius $\Delta$ can be
expanded. Otherwise, $\Delta$ is shrinked.

The Riemannian version of the algorithm consists in considering the tangent space to the manifold
for the search space, which is locally mapped to the smooth manifold. Formally, let $\cM$ be a
smooth manifold equipped with the Riemannian metric $\innert{\cdot}{\cdot}_{\cM, V}$ and the
associated norm $\normt{\cdot}_{\cM, V}$ defined on the tangent space $T_V\cM$ to $\cM$ at $V$. We
denote by $R_V :T_V\cM \to \cM$ a retraction which is a first-order approximation of the
exponential map that maps locally the tangent space to the manifold. The retraction is such that
$R_V(0) = V$. The correction is now defined as
\begin{equation}
    \label{eq:rtr_subproblem}
    \min_{W\in T_V\cM} m(W) = \sJ(V) + \inner{\nabla \sJ(V)}{W}_{\cM,V} +
    \frac{1}{2}\inner{\nabla^2\sJ(V)W}{W}_{\cM,V},
\end{equation}
such that $\normt{W}^2_\cM \le \Delta^2$,
where $\nabla \sJ$ (resp.\ $\nabla^2\sJ$) is the Riemannian gradient (resp.\ Riemannian Hessian) of
$\sJ$. The correction that belongs to the tangent space is mapped to the manifold using the
retraction, such that the new iterate is defined by $R_V(W)$. The trust-region radius is now adapted
according to the ratio
\begin{equation*}
    \rho = \frac{\sJ(V) - \sJ(R_V(W))}{m(0) - m(W)}.
\end{equation*}
The quadratic subproblem presented in Equation~\eqref{eq:rtr_subproblem} is solved with a truncated
conjugate gradient method. We refer the reader to~\cite{Absil2007} for an exhaustive description and
analysis of the algorithm as well as its application on the Grassmann manifold.

%% file: numerical.tex
%!TEX root = main.tex

\section{Application to Bayesian linear regression}
\label{sec:bayesian_linear_regression}

\subsection{Inference problem}
\label{sub:inference_problem}
\newcommand\Yt{\widetilde Y}
\newcommand\At{A_{\text{ref}}}
\newcommand\mt{m_{\text{ref}}}
\newcommand\Ct{C_{\text{ref}}}
\newcommand\mut{\mu_{F}}
\newcommand\mue{\mu_{E}}
\newcommand\kt{k_{F}}
\newcommand\ke{k_{E}}
\newcommand\st{\sigma_{F}}
\newcommand\lt{\ell_{F}}
\newcommand\gp{\cG\cP}

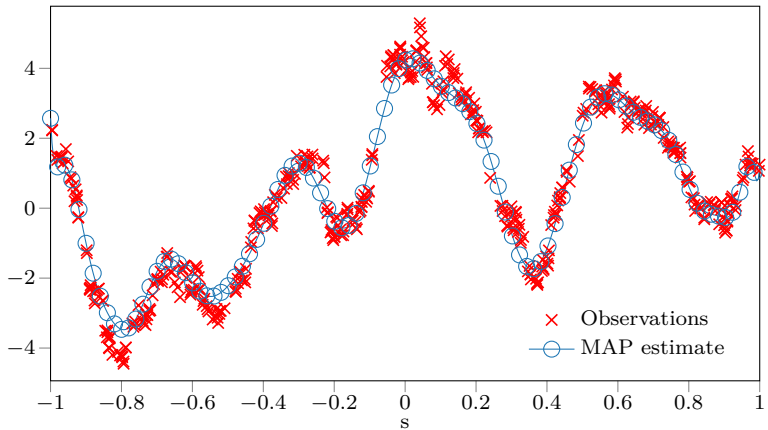
\begin{figure}
    \centering
    \setlength\figurewidth{0.9\textwidth}
    \setlength\figureheight{0.6\figurewidth}
    \input{figs/regression1d_obs_posterior.tikz}
    \caption{Comparison between the observations $y$ and the MAP estimate $B\ms + \me$.}
    \label{fig:regression1d_observations}
\end{figure}

\revis{
    The goal of this section is to illustrate the results of
    Section~\ref{sec:optimality_criteria}, and numerically assess the methods
    in the case of a Bayesian polynomial regression. Given a uniformly
    distributed sample $(s_i)_{i=1}^n$ in $(-1,1)$, we want to infer the 
    random variable, $X$, from the following linear model
    \begin{equation}
        \label{eq:linearapproxmodel}
        Y_i = \sum_{j=0}^{q-1} T_j(s_i)X_j + E(s_i), \quad \forall
        i\in\{1,\hdots,n\},
    \end{equation}
    where $T_j$ is the Chebyshev polynomial of the first kind~\cite{Gradshteyn} of
    order $j$ and $q=30$. The two moments of the prior distribution of
    $X\sim\cN(\mx,\cx)$ are defined by
    \begin{equation*}
        (\mx)_i = -1 + 2\frac{i-1}{q-1}
    \end{equation*}
    and
    \begin{equation*}
        (\cx)_{ij} = \sigma_\x^2 \left(1 + \sqrt{1200} \frac{\vert i - j\vert
        }{(q-1)}\right)\exp\left(-\sqrt{1200}\frac{\vert i - j\vert
        }{(q-1)}\right),
    \end{equation*}
    with $\sigma_\x = 1$.
    Note that the covariance $\cx$ is a Mat\'ern 3/2 covariance matrix, prescribing that polynomial coefficients associated to Chebyshev polynomials with distant degrees are less correlated than close ones.
    The noise $E$ is a stationary Gaussian process with mean and covariance
    functions defined respectively by $\mue(s) = \cos(4\pi s)$ and
    \begin{equation*}
        \ke(s,s')= \sigma_{E, 1}^2
            \exp\left(-\frac{\vert s - s'\vert}{\ell_E}\right) + \sigma_{E,2}^2
            \delta(s-s'),
    \end{equation*}
    with $\sigma_{E,1} = 0.6$, $\ell_E=0.05$, and $\sigma_{E,2}=10^{-3}$.
    The model presented in Equation~\eqref{eq:linearapproxmodel} is
    equivalent to the linear model from Equation~\eqref{eq:linearmodel} with
    $B_{ij}=T_{j-1}(s_i)$.

    Figure~\ref{fig:regression1d_observations} illustrates the data $y$ used
    for the observations and the maximum a posteriori fit $B\ms + \me$.
    The rest of Section~\ref{sec:bayesian_linear_regression} is dedicated to
    the optimal estimation of the posterior distribution $P(X\mid Y=y)$ using
    $P(X\mid W=V^T y)$, where $V$ has been computed according to the criteria
    introduced in Section~\ref{sec:optimality_criteria}.
}

\subsection{Numerical results}
\label{sub:numerical_results}

We first consider three types of approaches for the computation
of the reduced space $V$ based on the principal component analysis. They are denoted PCA-A, PCA-Y and
PCA-YN and are respectively computed as the dominant eigenvectors of the three following eigenvalue
problems
\begin{equation}
    \label{eq:kleeigen}
    \ca v = \lambda v, \quad \cy v = \lambda v \quad \text{and} \quad \cy\ce^{-1}v = \lambda v.
\end{equation}
PCA-A corresponds to the principal component analysis of $A(X)\revis{=BX}$, PCA-Y to the analysis of
$Y$, and PCA-YN to the analysis of $Y$ using the Mahalanobis distance, the metric induced by the
inverse of the noise covariance $\ce^{-1}$. The latter has been successfully used in a Bayesian
inference context in~\cite{Giraldi2017}, where the metric is directly involved in the posterior
distribution due to uniform priors.

We denote by KLD (resp.\ EKLD, MI) the solutions obtained using the minimization of the
Kullback-Leibler divergence (resp.\ minimization of the expected Kullback-Leibler divergence,
maximization of the mutual information).

\begin{figure}[htbp]
    \centering
    \setlength\figurewidth{0.8\textwidth}
    \setlength\figureheight{0.6\figurewidth}
    \input{figs/regression1d_kldiv.tikz}
    \caption{Kullback-Leibler divergence versus the dimension of the reduced space for the
    different numerical methods.}
    \label{fig:regression1d_kldiv}
\end{figure}
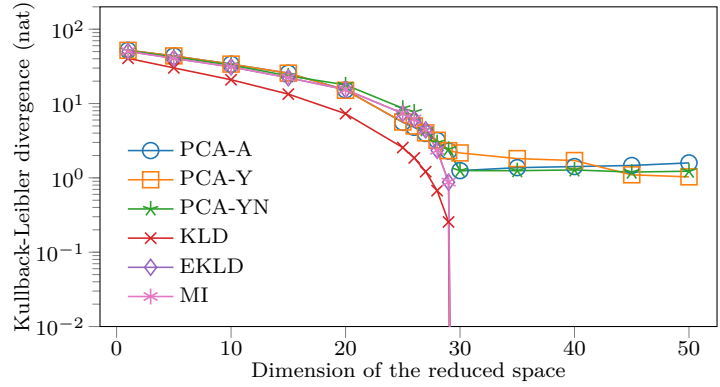

\revis{
    For a particular realization $y$, we compute the Kullback-Leibler divergence
    between the posterior distribution $P(X\mid W=V^T y)$ and $P(X\mid Y=y)$, and
    analyze its dependence on the dimension of the reduced space, $r$.  The results
    are plotted in Figure~\ref{fig:regression1d_kldiv} for the
    different dimensionality reduction methods. We conclude that the
    information theoretic based methods (KLD, EKLD, MI)
    with $r=q$ dimensions yield the exact posterior distribution within
    machine accuracy, and outperform the PCA-based approaches. Given that we
    are measuring the error using the Kullback-Leibler divergence,
    the KLD method performs better than the others. We can however note that
    the EKLD and MI techniques are robust to the realization $y$. We will
    observe below that the PCA methods require a dimension of the order of the
    total number of observations, $n$, to achieve a similar accuracy.
}

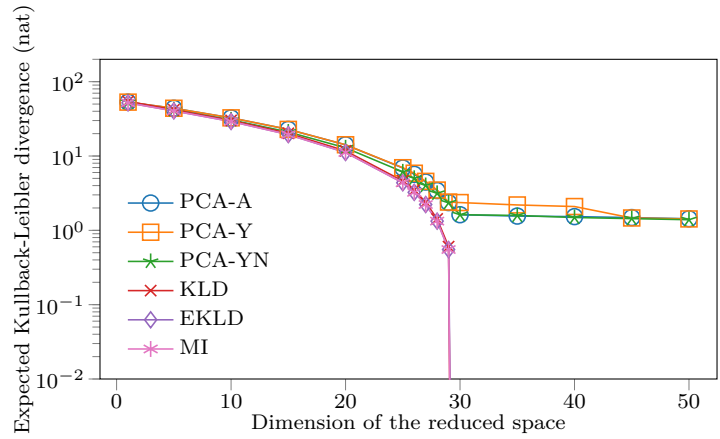
\begin{figure}[htbp]
    \centering
    \setlength\figurewidth{0.8\textwidth}
    \setlength\figureheight{0.6\figurewidth}
    \input{figs/regression1d_exp_kldiv.tikz}
    \caption{Expected Kullback-Leibler divergence versus the dimension of the reduced
        space for the different numerical methods.}
    \label{fig:regression1d_exp_kldiv}
\end{figure}

\revis{
    Figure~\ref{fig:regression1d_exp_kldiv} depicts the dependence
    the {\it expected} Kullback-Leibler divergence between the posterior
    distributions of the reduced and the full models on the
    dimension of the reduced spaces; plotted are results obtained using the 
    different projection techniques.
    Similar to Figure~\ref{fig:regression1d_kldiv}, the information
    theoretic approaches converge to the posterior distribution with
    subspaces of dimension $r=q$, which is not the case for the PCA methods.
    We also note that even when the expected Kullback-Leibler
    divergence is used as error criterion, the EKLD method does not really
    improve the speed of convergence of the distributions compared to the
    other information theoretic approaches. 
}

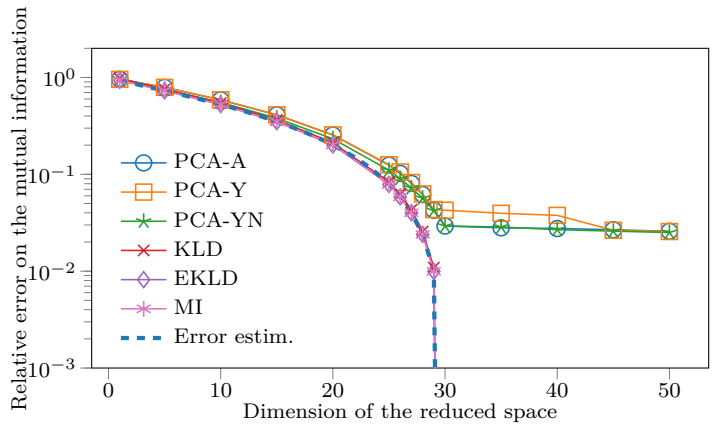
\begin{figure}[htbp]
    \centering
    \setlength\figurewidth{0.8\textwidth}
    \setlength\figureheight{0.6\figurewidth}
    \input{figs/regression1d_mi.tikz}
    \caption{Relative error between the mutual
        information of the observations and the parameter of interest of the full model and the one of
        the reduced model and error estimator~\eqref{eq:apriorierrorestimate} versus the dimension of the
    reduced space for the different numerical methods.}
    \label{fig:regression1d_entropy}
\end{figure}

\revis{
    Figure~\ref{fig:regression1d_entropy} illustrates the relative error with respect to the dimension of the reduced space between the mutual information of:
    \begin{itemize}
       \item the observations and the parameter of interest, $\cI(Y,X)$; and,
       \item the projected observations and the parameter of interest.
    \end{itemize}
    We are in fact looking at the criterion introduced in
    Corrolaries~\ref{coro:maxmi} and~\ref{coro:maxmi_nu}. We note that for
    $r\ge q=30$, all the information theoretic methods converge to the minimal
    value of the relative error. This behavior is predicted by
    Corollary~\ref{coro:maxmi_nu} for the MI approach, as illustrated in the
    figure by the fact that the error
    estimator~\eqref{eq:apriorierrorestimate} overlaps with the error of the MI
    approach. Again, the PCA based methods perform poorly when compared
    to the information theoretic approaches.
    }

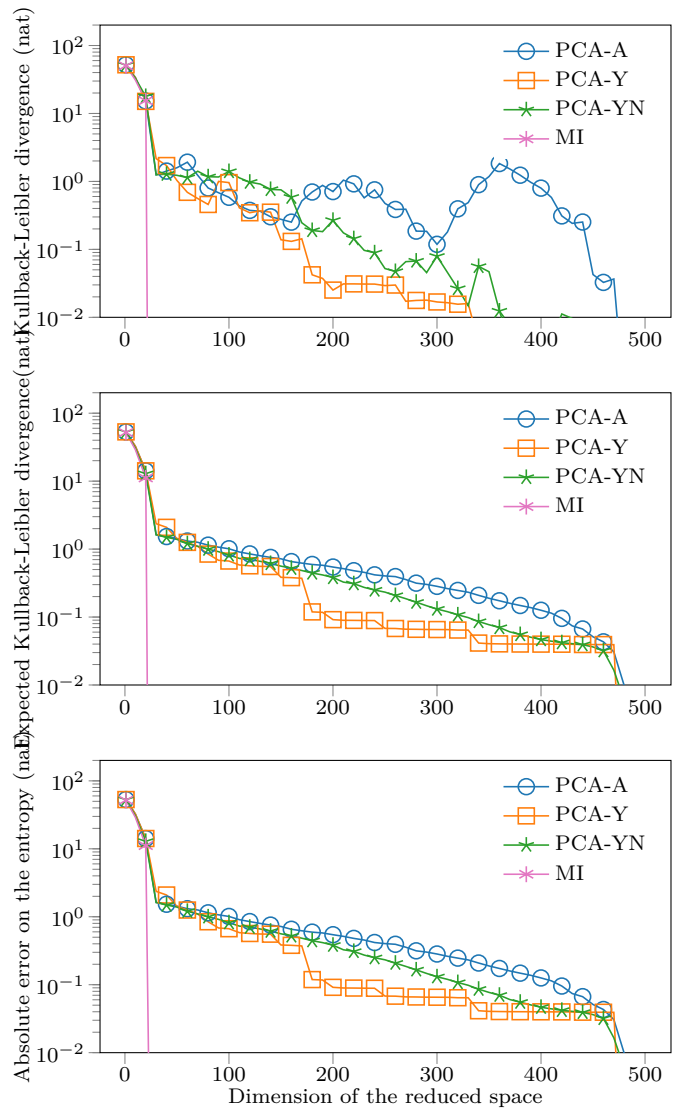
\begin{figure}[htbp]
    \centering
    \setlength\figurewidth{0.75\textwidth}
    \setlength\figureheight{0.6\figurewidth}
    \input{figs/regression1d_kle.tikz}
    \caption{Kullback-Leibler divergence (top), expected Kullback-Leibler divergence (middle), and
        entropy (bottom) versus the dimension of the reduced space for PCA-A, PCA-Y, PCA-YN and MI
    methods.}
    \label{fig:regression1d_kle}
\end{figure}

In Figure~\ref{fig:regression1d_kle}, the different divergences and the absolute error on the entropy with respect to the
dimension of the reduced space are illustrated for the \revis{PCA}-based methods for larger
values of the dimension and compared to the MI approach. \revis{The absolute
error on the entropy is equivalent to the error on the mutual information up
to a constant according to Theorem~\ref{th:maxmi}.} One can see that the
dimension of the
reduced space must be an order of magnitude larger compared to the MI technique to reach the same
accuracy for all three convergence criteria.

\begin{figure}[htbp]
    \centering
    \setlength\figurewidth{0.8\textwidth}
    \setlength\figureheight{0.6\figurewidth}
    \input{figs/regression1d_kle_singval.tikz}
    \caption{Normalized singular values for the PCA-A, PCA-Y, PCA-YN methods.}
    \label{fig:regression1d_kle_singval}
\end{figure}
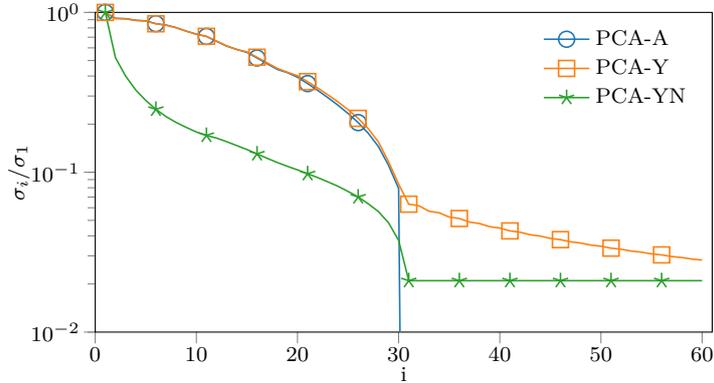

In Figure~\ref{fig:regression1d_kle_singval}, the normalized singular values $(\sigma_i/\sigma_1)_{i=1}^n$
computed for the PCA methods are illustrated. The singular values are defined by
\begin{equation*}
    \sigma_i = \sqrt{\lambda_i}, \quad \text{such that} \quad \sigma_1 \ge \sigma_2 \ge \hdots \ge
    \sigma_n,
\end{equation*}
where $(\lambda_i)_{i=1}^n$ are the eigenvalues involved in Equation~\eqref{eq:kleeigen}.  It is
shown that the spectrum resulting from the PCA-YN method decays faster than the other approaches.
Moreover note that the eigenvalues involved in the MI approach (i.e.\ eigenvalues of
Problem~\eqref{eq:maxmi_eig}) are strictly equal to the eigenvalues of the PCA-YN technique,
see~\cite{Giraldi2017} for more details.

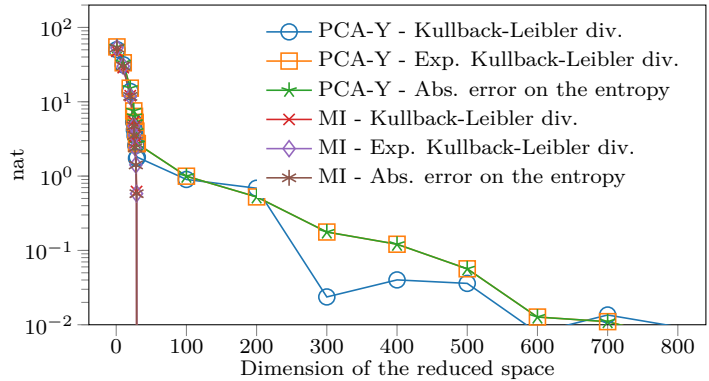
\begin{figure}[htpb]
    \centering
    \setlength\figurewidth{0.8\textwidth}
    \setlength\figureheight{0.6\figurewidth}
    \input{figs/regression1d_largenobs.tikz}
    \caption{Kullback-Leibler divergence, expected Kullback-Leibler divergence,
    and the entropy for the PCA-Y and MI methods versus the dimension of the reduced space
for a larger number of observations ($n=2000$).}
    \label{fig:regression1d_largenobs}
\end{figure}

For the last experiment, we only consider the PCA-Y and MI approaches.  The convergence of the
Kullback-Leibler divergence, the expected Kullback-Leibler divergence, and the entropy with respect to the dimension of the reduced space is plotted in Figure~\ref{fig:regression1d_largenobs} for the PCA-Y
and MI methods using an even larger number of observations ($n=2000$). 
For a dimension $r=30$, the three quantities of interest are null within machine precision for the MI method, whereas the PCA-Y approaches needs a dimension $r=700$ to get a value of $10^{-2}$
nat. 
This highlights that the accuracy of the MI method is more related to the number of parameters
($q=30$) than the number of observations ($n=500$ in
Figures~\ref{fig:regression1d_kldiv},~\ref{fig:regression1d_exp_kldiv},
and~\ref{fig:regression1d_entropy}, and $n=2000$ in Figure~\ref{fig:regression1d_largenobs})
as predicted by Corollary~\ref{coro:maxmi_nu} for the relative error on the mutual information.

\revis{
\subsection{Summary}
}

We have seen that regarding the accuracy on the posterior, the KLD, EKLD, and MI approaches perform
much better in terms of Kullback-Leibler divergence than the PCA approaches. On
the other hand, in terms of computational efforts to evaluate the basis, the PCA methods and the MI
technique only require the solution to an eigenvalue problem, whereas the others need
more advanced strategies like the Riemannian optimization algorithm presented in
Section~\ref{sub:numerical_solution_of_the_optimization_problems}. The maximization of the mutual
information therefore exhibits a good balance between posterior distribution accuracy and
computational difficulty, further providing an a priori error estimate as well as an upper
bound on the number of required projections.

%% file: figs/regression1d_obs_posterior.tikz
% This file was created by matplotlib2tikz v0.6.11.
\begin{tikzpicture}

\definecolor{color0}{rgb}{0.12156862745098,0.466666666666667,0.705882352941177}

\begin{axis}[
xlabel={s},
xmin=-1, xmax=1,
ymin=-4.93666728758486, ymax=5.775534599079,
width=\figurewidth,
height=\figureheight,
tick align=outside,
tick pos=left,
x grid style={lightgray!92.026143790849673!black},
y grid style={lightgray!92.026143790849673!black},
legend style={at={(0.97,0.03)}, anchor=south east, draw=none},
legend cell align={left},
legend entries={{Observations},{MAP estimate}}
]
\addplot [semithick, red, mark=x, mark size=3, mark options={solid,fill opacity=0}, only marks]
table {%
-0.410669994625781 -0.138946020600415
0.0611735112105882 4.02651566636267
-0.61695842610501 -1.99350200496081
-0.864199283617417 -2.30161116707613
0.573970919999827 3.14784251734261
0.312667043551711 -0.125350871713626
0.275041792087272 -0.0789338760814668
0.151205787506068 3.29061367981734
-0.921874167622267 -0.272910849084419
-0.284372791032902 1.52909891381358
0.891366373636815 -0.0371351899982812
-0.879910639369229 -2.37635807456035
0.728084207092422 2.3219944854837
0.754581052296061 1.73573174377421
-0.897612668759093 -1.28600638795025
0.30483723093131 -0.566653015951328
0.103502737269751 3.30321145814298
0.195026506478102 2.99558455736037
-0.0329427513799958 4.37582517929709
-0.434023678286018 -0.963338988820431
-0.404548563396457 -0.124687380041599
0.123017810652722 4.25465459786673
-0.207905127905947 -0.667795422087362
0.577401419357799 2.90092318223903
-0.163031229308246 -0.160868812807661
-0.712192158449344 -2.49348856126571
-0.698186611158646 -2.02465470008567
-0.889517299568352 -2.2941420768393
0.436074387407291 0.27760576327186
-0.415365290991723 -0.573354084925982
-0.602452255599314 -2.38745498002002
0.662727844621042 2.82404578068093
0.135982239449262 3.72290776219139
-0.835320501575038 -3.98539478683994
0.0899971478614408 2.82689274458616
-0.682082510865449 -1.90617796375762
0.353524771405921 -1.89969930129068
-0.763055470241847 -3.34123006022032
-0.110007867759862 0.182339518936577
0.775964868314398 1.57783086060484
0.59453545999062 3.6220233014937
-0.864042543472747 -2.30111472919654
0.921575512933556 0.267738523638609
0.318410784752043 -0.67896077612395
0.437552166448869 0.432287495130988
0.487153530268193 1.49084569748278
0.775472703771274 1.73359586363775
-0.732108664848414 -3.34365018883247
0.55394916121116 2.7318067194236
0.67598302439065 2.72359024745868
0.658609527816856 2.80323646210879
-0.941686006960382 0.842670334398486
-0.219251017446831 -0.139932933099826
-0.125496015663918 -0.0719323946730311
0.0412701270226636 5.28861633150337
-0.0142546516686317 4.391704424252
-0.739141435582762 -3.09896902811432
0.968767061988944 1.47742826730526
-0.868595356163829 -2.33715722376277
0.15856815548814 3.17824396564308
0.161549650173456 3.11615943189872
0.908675731113848 -0.353525901744368
-0.605112444972649 -2.40714423653319
-0.83334919873919 -3.72436485955486
0.767728435466474 1.6787642522778
-0.968607614098113 1.47541606696489
0.591189294707049 3.71981392101898
0.536380878856624 3.34426369977545
0.744087571906341 1.78330550171809
0.832220172779111 -0.00805343510153733
-0.452373614164431 -2.03607299608747
-0.527237747234801 -3.2834028068846
-0.378960312709225 -0.100553516812359
0.431286517836948 0.217099316644453
-0.834529021493193 -4.0072318101295
0.350972886008853 -1.65154460867044
0.717049337755054 2.59044114714737
0.385793768233012 -1.79667051743711
0.570167382527859 3.09967050279207
-0.293437496438761 1.19788346189466
-0.559615861912525 -2.80177170233123
-0.5092552190805 -2.86236922648099
0.454825031148462 1.08219233234151
-0.871744454195582 -2.36621803464763
-0.0925250761155862 1.54959463877156
-0.0246453410629033 4.12440511866747
0.28174965888899 -0.22936046656185
-0.787662295273323 -3.98150047130443
0.147219875602075 3.70372761999996
0.920691588048089 0.278792089845953
0.0313611603950799 4.49849918238339
0.139650091789694 3.979589683537
0.976013323142853 1.28174029870579
0.665741870138125 2.72197674615638
-0.00371263579200765 4.18118359915137
-0.83969732427194 -3.49497277133514
0.742093988129925 1.90863802109104
0.426172817974849 -0.185675074259022
0.464224087069268 0.757274190589388
0.647239974040031 2.875869685911
0.489958422433613 1.95129739353696
-0.173874743808285 -0.669017426250496
0.0294106677512318 4.18351506173718
0.6930829902341 2.39160438454389
0.717724113693793 2.58637224898804
-0.378992772113617 -0.083019194041825
-0.761227212042893 -3.24523821719232
0.0419189626610015 5.16908994731361
-0.102099863405943 0.469741285741274
0.794574549863618 0.900569975538754
-0.0137953228269467 4.6256939712665
0.00895503598318226 3.68891440904054
-0.269752701219188 1.15093531673556
0.566012242184281 3.28592254484472
0.895846038417131 -0.263385626115889
-0.739969973852491 -3.22613394222271
0.322225347052141 -0.476112017909365
-0.672535383019665 -1.3887978611972
-0.382122471185104 -0.275891310859017
0.933260308545161 0.0662962989899303
-0.0342906506691074 4.27152554334953
0.0814948903885861 2.94658559822511
0.422346270137638 -0.101194474030324
0.374863488132597 -2.11911489005338
-0.926061704869924 0.324365118910394
0.776341131911714 1.57468318704386
-0.622831583933737 -1.6356538110507
0.564055516056089 3.37074384638808
0.287852243944277 -0.496621527651893
-0.0292020277962057 3.86162760752881
0.858089890198391 0.109893381573188
-0.460237174609675 -2.0632099474421
0.628357971949269 2.40703709062896
0.390084146808238 -1.72716212854525
-0.528937616211175 -3.13409757554084
0.906605383682072 -0.63480439952063
0.170871417079452 3.26338430289525
0.166132937350707 3.07320643533311
-0.685580915132111 -2.03145726448656
-0.149519071826828 -0.301862669914479
-0.529956171357519 -3.03421018745418
0.201014653239473 2.50486447574607
-0.278783154991592 1.30123860104765
-0.662253103294113 -1.81529466416755
-0.305507101965536 1.49733214625302
-0.57096596062845 -2.76358786426397
0.589056585993052 3.51094203371073
0.583418636192907 3.39440930297515
-0.670919286271143 -1.28130225038687
0.500728097024525 2.64296441288207
-0.888206762307973 -2.33673170133187
-0.53505914879656 -3.10379432245766
0.300192862989738 -0.0916348540760907
0.436483527849729 0.313984122528375
-0.0399010729988909 4.30611266402603
-0.324770664685247 0.89683030082916
-0.220587421356537 0.0222852953035115
0.717594387690007 2.58385710000757
0.0227712222575993 3.75481710820614
0.749380028631645 1.76204712737978
0.892877108309104 0.161901786817749
0.871467422809647 0.0726391784274811
-0.013924136125794 4.58994636224811
0.201036551115646 2.50938596343627
-0.585659799719637 -1.78714966103929
0.912376544090228 -0.143184899971416
0.973620918362694 1.22122542485515
0.658757686706906 2.77021893540151
0.0757789193836036 3.55683087748089
-0.832924759227976 -3.68850478259932
-0.372225021996141 0.172442107791182
-0.787773101540839 -3.98666282986586
-0.59018270841226 -1.65570781179815
-0.34896320737544 0.291060356962718
0.00321797451023742 3.93579588988421
0.756637507621828 1.58530152675498
-0.477982113277742 -2.3024648937909
0.0602735451057168 4.15379682702951
0.687725284233761 2.2900613711657
-0.726871027190767 -3.05855284467544
-0.199551979421401 -0.815212080938934
0.136970573111595 3.86646435054368
-0.627917574879888 -1.72813879331827
-0.972568133413937 1.45886433913081
0.185578111495259 2.68419919946983
-0.968569535820369 1.44431672766127
0.171987447144243 3.31759653613066
-0.182674197368551 -0.610000228380583
-0.54275211624259 -3.18479515407862
-0.620232524258557 -1.65584638398165
0.826605839929452 0.0861245290269833
0.180440743770536 3.15746026269286
-0.657611305437818 -1.82736626742202
0.428948955481286 0.0781163386584864
-0.523107882271393 -2.9391613042272
0.84361028649879 -0.27284187541687
0.442547733499124 0.595827975801031
0.353845958175338 -1.88654821443604
0.973621112686098 1.22334330841273
0.0794903929949811 3.03757335992202
-0.00538584455499924 4.16620379868629
0.796254756010246 0.875055836739996
0.923053702434934 0.338146080841036
-0.995589813028686 2.23960036753757
0.662203656008339 2.82008021186455
0.11452223557287 4.3463533071604
0.870033391256072 0.0918510289617346
0.436882482993901 0.306768705856068
0.91616961772445 0.0938758032157821
0.518300980398082 3.42013047834429
-0.47221154821935 -2.02037744255514
-0.72849693546465 -3.04514737619923
0.803178891568455 0.772171529785623
0.0321464958333091 4.51948692017879
0.698872977754414 2.75311855783139
-0.889405685049886 -2.31839523194242
-0.140096004154132 -0.387150763368926
-0.588865863652135 -1.70921759648114
0.741182598426402 1.72787433914695
0.38806798328864 -1.47584653558931
-0.875141597670941 -2.7239522130512
0.229526122381339 2.1703293660509
0.678569519479552 2.79367225185257
-0.932435387784019 0.198412937489319
0.225230198488522 2.21922667405682
0.851634186052453 0.173967435097716
-0.697078412955228 -2.08536346805293
0.45587273292868 1.03761289343157
-0.69835863108557 -2.03382081110745
0.679903680323166 2.92442477913539
-0.634764478767152 -2.54481410333535
-0.217793226713742 -0.114556976080536
-0.279593846448136 1.2413184582547
-0.204115674000835 -0.858406626062584
0.980433710991078 1.26928126759732
-0.60251483252729 -2.33647487810805
-0.261493379658972 1.53357848993847
0.538788449302666 3.36366785109777
-0.109059085554458 0.316298826815653
0.112540409965786 4.0692769305917
0.527601462970915 3.40793996204339
0.999566767299233 1.24517722239598
0.635850193245846 3.04992083975667
-0.122083716545092 0.177238819619523
-0.316046787831467 0.812042619751878
-0.327976625689 1.00441923459459
0.758925956877943 1.63583842564088
0.503121210848882 2.73214821140238
-0.845759086668627 -3.50557767893958
0.052489263265534 4.5881583312032
0.346134845408951 -1.80754263860979
-0.0247581436220847 4.10232425191444
0.0479252363702221 4.60666278564507
0.717734077069096 2.57025646911913
-0.328702651540655 0.917064980471782
-0.729689136623047 -3.17824868617093
0.219394360152616 2.16993498846943
0.320439915116482 -0.588794187380937
0.954064132238937 0.875788094710464
0.644443314350072 3.02659840666856
0.55731444016483 3.06970303360963
0.88570448494802 -0.460650413331786
-0.202566698461057 -0.914792317583038
-0.395523775261513 -0.119405564099677
0.386634417864964 -1.67101318169479
0.118801096976927 4.12103036070976
0.031693439349626 4.49386712657556
-0.554601536810371 -2.95549864705409
0.299719851662232 0.0538269844414466
0.354753312693679 -2.03876395867782
-0.227304152174576 0.8794800974803
-0.477175531775665 -2.39622715629912
0.997522900580435 1.06613873615339
-0.471168593308261 -2.14867883475679
-0.670496518493535 -1.37937763590908
0.566395933056352 3.21849565941245
-0.0395361742964495 4.34666278034768
0.906636126550173 -0.646455161809987
-0.274323798676561 1.20453465084154
0.206071690813353 2.16966288316711
-0.980977033849658 1.38125511983662
0.536359041492025 3.35667118444378
0.931773659042579 -0.0533500807150118
0.371596422417038 -2.16551055322768
0.674628540296516 2.3863287240373
-0.793956030299199 -4.44974902000923
0.576523446314897 2.84498580768029
-0.138801789686956 -0.563214499001919
-0.649484086591137 -2.12095289719293
0.800755333124209 0.753626634412851
0.326038785096809 -0.990987157794914
-0.568528064499754 -2.65559687347732
-0.465508200998512 -1.84125254837646
-0.689423731878017 -1.81565464009381
0.306174813210452 -0.630921519523017
0.418701776118032 -0.271884607652718
-0.875894559834903 -2.64223011613961
0.952058425765947 0.802609436219262
-0.396494456208046 -0.195114848782878
0.313759217589291 -0.329551352956808
0.590886598092975 3.65806356819164
0.750187888266413 1.88868585113063
0.411909273096044 -0.482086867169404
0.326580665648037 -0.91019248811056
-0.041834394478877 4.05353160110376
-0.282043168794912 1.34304939450297
0.908735283708336 -0.385543028206518
0.218341000611219 2.20486252391987
0.713950660517935 2.38006509762667
0.760853720712857 1.60505894146578
-0.641588458236289 -2.26232284484361
-0.566304655105117 -2.59066795943235
-0.929270767104632 0.539769298510766
-0.699022942943709 -2.05211011664491
-0.769235450203571 -3.29348689523126
0.807369623127034 0.309631503989934
0.96024492878055 1.24464631389693
0.0198237323454729 3.92414559535524
-0.107142922054198 0.277177431565635
0.369396438885013 -2.02052206481327
-0.75187974859819 -2.77736449327309
0.397721824932787 -1.37700398511639
-0.121208220435942 0.274331676844329
0.311574038349934 -0.103197570986536
0.59390075571071 3.68704912207035
0.0967879269169201 2.94932857644938
-0.389497626936111 -0.164959219234596
-0.341578870989156 0.556596834126786
-0.579207851740531 -2.32412493898046
0.57491680011325 2.97680874634955
-0.191635531573705 -0.277384315505734
0.419732754821911 -0.196364755533721
-0.424769897128919 -0.348118884479308
-0.427749194218725 -0.577016803328421
-0.33336797770989 0.68878411852235
0.0108519412578014 3.72936703062613
0.84363963715131 -0.318904720834673
0.443409959699199 0.5868911360896
-0.46608861601503 -1.80567290925427
-0.925120769723364 0.233606361092311
-0.115450741793187 0.379894177667353
-0.420479825253155 -0.397075743724674
0.543742687742555 2.81899469258992
-0.198803937593422 -0.821900628866771
-0.627547121928147 -1.7762629150689
0.776632933714388 1.44575417529912
0.618681003554264 2.74863968662527
-0.31570315527247 0.689988664720822
-0.424617491146014 -0.364941043042957
0.239784934363712 0.86610315086987
0.542437382479013 3.14136123527198
0.404223872891808 -1.43427652226558
0.758310334346032 1.72848255390327
0.175607486951379 2.90831294861446
-0.435137147903934 -1.11945506599142
0.2087954658054 2.27630940141064
0.288689315031129 -0.528388101257209
-0.852512992625788 -2.63900491637089
-0.80393093116625 -4.21323368892968
-0.175006035232822 -0.773475252097411
-0.324393745782642 0.890279380825948
-0.310636765360326 0.978928801925615
-0.0135666005777009 4.56255583987456
-0.367409750136763 0.380765040972016
-0.605299745267629 -2.47294291230693
-0.123179670108359 0.0351441255946968
-0.620367152526619 -1.6609722515335
0.372474658546039 -2.13579363532569
-0.207784331173356 -0.62431800621642
-0.268718793038976 1.16793862659273
0.190556525837167 2.69874995870051
-0.718422670233507 -2.91918459273438
0.298144699872617 -0.2454846892885
-0.766070893209926 -3.39366997746758
0.626420197415851 2.30416742471897
-0.0524701492050101 3.75161617382974
-0.851765227960799 -2.65055734927034
-0.128761055119511 -0.107608725094178
-0.00837325044028647 4.23183769014358
-0.147935253756444 -0.355392547622943
-0.798469655007747 -4.23574979332375
0.874268761750639 -0.0764430617516075
-0.0941735235838106 1.3840863417868
0.829560429933741 0.0705834789897133
0.966648750311788 1.63361348521563
-0.164143924150996 -0.202346616749658
0.901812615008827 -0.589625316716968
-0.210480055696361 -0.694121052980928
-0.194515587353975 -0.247000848494621
-0.809985009683008 -4.11135287000773
0.644929391654487 3.05705378108841
-0.397994366808246 -0.0259879293352441
0.199398091281003 2.7585620576225
0.192029292652811 2.80499472544611
-0.0328981502854002 4.37911488465935
-0.944296717059328 0.888305006718696
-0.407469973195568 -0.0896881287155228
0.816178294451865 0.102247906258914
0.0778860720717687 3.02785448716841
-0.707906205852431 -2.17168822325967
-0.696796343565563 -2.05627996490349
-0.189379185236218 -0.489578833166149
-0.736444614933442 -3.30902442540013
0.614109859687137 3.19136196738561
-0.122711694956284 0.0386251486646022
-0.596086804821536 -2.33510333687602
-0.582035153120867 -1.96982615923038
-0.935977016342992 0.644315766205019
0.40304900822143 -1.64394982849444
0.40552695291902 -1.20675848023417
0.545096118407953 2.97150018016922
0.892076448575368 0.166511498228095
-0.0239992229306349 4.19452963406551
0.902557325993711 -0.701263175206186
-0.557224913499615 -2.89109408109546
0.795450333054014 0.814518925472434
-0.839576131323526 -3.5459492115217
-0.995466933596108 2.22305641610327
-0.825587892860259 -4.19615235551518
0.0625354260045028 3.95457921732542
-0.375807117597847 -0.511254284084065
0.539037949958335 3.3734522430353
-0.162392411981928 -0.0341427097316159
0.106131657927855 3.38316919959852
-0.102859872203269 0.49796552326394
-0.921677652054453 -0.254793794002348
-0.478609727639744 -2.37674845306579
-0.898085065349719 -1.20318337368363
-0.794462916416744 -4.38725179004066
0.670212687252997 2.54641458644157
0.559447072531226 3.13844737477605
0.510696890295155 2.84013549067826
0.703343101859606 2.4626642121168
-0.450844719355934 -2.12160305110829
-0.261044325280416 1.42161507481223
0.817989061177029 -0.0805886277965977
-0.982558653245747 1.49347134003575
-0.385406240452132 -0.206212329341193
-0.448896263959396 -1.73410338672308
0.207704775180667 2.27887832349805
0.629693715707534 2.59188361456024
-0.519930476487301 -2.93662419307516
0.899768553386638 -0.486939654984531
-0.335952874655152 0.901751390561661
-0.00870373023030613 4.33092394100357
-0.841154717266996 -3.69230610132897
-0.361985102846052 0.291966702250354
-0.2755652372299 1.47993556769856
0.372738552730616 -2.19495409070628
-0.374571760286871 -0.434855221995799
0.636132886811773 2.9578245909812
-0.579483735460576 -2.29636652097262
-0.180465393529018 -0.774869377653891
-0.523795115574647 -3.02408753952416
0.0454084034581881 4.91861827717361
0.862047904989825 0.183531050597685
0.977050883067857 1.15307145608063
0.912198906344814 -0.140477424944159
-0.231110677325687 1.36002482188629
0.0156257711107926 3.74886911076526
0.0165411217513571 3.87682398939319
0.158122262553729 3.09028920725138
-0.0926546712303258 1.51054389707513
-0.531873319704185 -2.88344148286973
0.268718192401536 0.0643279282603353
-0.7213372006078 -2.88119667568453
0.203788343807249 2.49478098097554
-0.7215277655336 -2.9360370414541
-0.800681695520803 -4.19119843972688
0.401085719845246 -1.58431083448332
-0.0509568082727545 4.07401196970506
0.914884295225476 0.0830979031228437
-0.136080366249374 -0.60315218979555
-0.886971725789768 -2.26708468405353
-0.966077916809321 1.282527922672
-0.229554105730165 1.28559129955834
-0.972587620612339 1.43953584064461
0.0588605538000992 4.14662908830824
0.10784604140624 3.62986579635669
0.708614740793614 2.43052584541873
0.390471197451458 -1.76199674729506
0.310640462091239 -0.271564033928769
0.919379851560795 0.25475449251716
-0.686082206403881 -2.11807052879189
0.814174916684824 0.288062732252608
0.0541884089971019 4.32168642545242
-0.950209491239014 1.33868683349518
-0.162483196706799 0.0153756651124328
0.274340850411777 0.0689624289783893
0.518098092393368 3.47766143337526
-0.797981364358858 -4.26446515343843
0.972653374199468 1.29437551762849
-0.886896617971328 -2.20124685785669
0.136802887492569 3.79739141840372
0.697215163072812 2.62928440368819
0.522805719163475 3.46707304550462
0.333031290308359 -0.819450357057945
-0.92613891374316 0.365732793612752
0.729215888773994 2.00856581655926
-0.956819867014323 1.69076027333879
};
\addplot [color0, mark=*, mark size=3, mark repeat=5, mark options={solid,fill opacity=0}]
table {%
-1 2.57088425392006
-0.995991983967936 1.8730820667173
-0.991983967935872 1.47807480607026
-0.987975951903808 1.27727895631383
-0.983967935871743 1.19518826085465
-0.979959919839679 1.18001521396908
-0.975951903807615 1.19681437416769
-0.971943887775551 1.22248452056319
-0.967935871743487 1.24217800527089
-0.963927855711423 1.24675164206695
-0.959919839679359 1.23097843165431
-0.955911823647295 1.1923070247125
-0.95190380761523 1.13000917016667
-0.947895791583166 1.04459710577561
-0.943887775551102 0.937425119562043
-0.939879759519038 0.810414181477524
-0.935871743486974 0.665857143699286
-0.93186372745491 0.506275799952749
-0.927855711422846 0.33431111965751
-0.923847695390782 0.152635081743781
-0.919839679358717 -0.0361225845355458
-0.915831663326653 -0.229436233529291
-0.911823647294589 -0.424936211719992
-0.907815631262525 -0.620451239038701
-0.903807615230461 -0.814041167918294
-0.899799599198397 -1.00401993193113
-0.895791583166333 -1.18896864599903
-0.891783567134269 -1.36773909822695
-0.887775551102204 -1.53944819719464
-0.88376753507014 -1.7034642519507
-0.879759519038076 -1.85938622916981
-0.875751503006012 -2.00701733259549
-0.871743486973948 -2.14633437515286
-0.867735470941884 -2.27745446345502
-0.86372745490982 -2.40060049307557
-0.859719438877755 -2.51606686987986
-0.855711422845691 -2.62418673899898
-0.851703406813627 -2.72530183073239
-0.847695390781563 -2.81973583384071
-0.843687374749499 -2.90777199279235
-0.839679358717435 -2.98963540693973
-0.835671342685371 -3.06548029537121
-0.831663326653307 -3.13538228886467
-0.827655310621243 -3.19933562597262
-0.823647294589178 -3.25725496829337
-0.819639278557114 -3.30898141349536
-0.81563126252505 -3.35429217541802
-0.811623246492986 -3.39291331916293
-0.807615230460922 -3.42453488510377
-0.803607214428858 -3.44882770793216
-0.799599198396794 -3.46546123329929
-0.795591182364729 -3.47412165287528
-0.791583166332665 -3.47452971592926
-0.787575150300601 -3.466457628802
-0.783567134268537 -3.44974451975883
-0.779559118236473 -3.42431002252695
-0.775551102204409 -3.39016561427607
-0.771543086172345 -3.3474234299896
-0.767535070140281 -3.29630236240737
-0.763527054108216 -3.23713134257203
-0.759519038076152 -3.17034977834442
-0.755511022044088 -3.09650520524555
-0.751503006012024 -3.01624827413101
-0.74749498997996 -2.93032526232807
-0.743486973947896 -2.83956834809845
-0.739478957915832 -2.74488393206359
-0.735470941883767 -2.64723932325626
-0.731462925851703 -2.54764813171176
-0.727454909819639 -2.44715472418304
-0.723446893787575 -2.34681810505636
-0.719438877755511 -2.24769558143099
-0.715430861723447 -2.15082656032116
-0.711422845691383 -2.05721680786803
-0.707414829659319 -1.96782347622233
-0.703406813627255 -1.88354117434111
-0.69939879759519 -1.80518932533014
-0.695390781563126 -1.7335010161593
-0.691382765531062 -1.6691135065709
-0.687374749498998 -1.612560523745
-0.683366733466934 -1.56426642868711
-0.67935871743487 -1.52454230020604
-0.675350701402806 -1.49358394352097
-0.671342685370741 -1.47147179366361
-0.667334669338677 -1.45817264952161
-0.663326653306613 -1.45354314310614
-0.659318637274549 -1.45733482083121
-0.655310621242485 -1.46920068957843
-0.651302605210421 -1.48870306031091
-0.647294589178357 -1.51532250612653
-0.643286573146293 -1.54846773994958
-0.639278557114229 -1.58748620951933
-0.635270541082164 -1.63167520383842
-0.6312625250501 -1.68029326562158
-0.627254509018036 -1.7325717083076
-0.623246492985972 -1.78772604358579
-0.619238476953908 -1.84496713582221
-0.615230460921844 -1.90351191289858
-0.61122244488978 -1.96259347841963
-0.607214428857715 -2.02147048760974
-0.603206412825651 -2.07943566810366
-0.599198396793587 -2.13582338683455
-0.595190380761523 -2.19001618493681
-0.591182364729459 -2.24145022362314
-0.587174348697395 -2.2896196049941
-0.583166332665331 -2.33407955234599
-0.579158316633267 -2.37444845443505
-0.575150300601202 -2.41040879704131
-0.571142284569138 -2.44170702279138
-0.567134268537074 -2.46815237632142
-0.56312625250501 -2.48961480629819
-0.559118236472946 -2.50602200841505
-0.555110220440882 -2.51735570412422
-0.551102204408818 -2.52364725847764
-0.547094188376753 -2.52497274698207
-0.543086172344689 -2.52144758582183
-0.539078156312625 -2.51322084218691
-0.535070140280561 -2.50046934182096
-0.531062124248497 -2.48339168935409
-0.527054108216433 -2.46220231361836
-0.523046092184369 -2.43712564508922
-0.519038076152305 -2.40839052600412
-0.51503006012024 -2.37622494574578
-0.511022044088176 -2.34085118492236
-0.507014028056112 -2.30248144141791
-0.503006012024048 -2.26131400072028
-0.498997995991984 -2.21753000125953
-0.49498997995992 -2.17129083350637
-0.490981963927856 -2.12273619938638
-0.486973947895792 -2.07198284635205
-0.482965931863727 -2.01912397840743
-0.478957915831663 -1.96422933467375
-0.474949899799599 -1.90734591488253
-0.470941883767535 -1.84849932063501
-0.466933867735471 -1.78769567150865
-0.462925851703407 -1.72492404624164
-0.458917835671343 -1.66015939138638
-0.454909819639279 -1.59336583307695
-0.450901803607214 -1.52450032196946
-0.44689378757515 -1.45351653703451
-0.442885771543086 -1.38036897073833
-0.438877755511022 -1.30501711625256
-0.434869739478958 -1.22742967667809
-0.430861723446894 -1.14758871683019
-0.42685370741483 -1.06549367987409
-0.422845691382766 -0.981165193966306
-0.418837675350701 -0.894648597981755
-0.414829659318637 -0.806017120309839
-0.410821643286573 -0.715374650494391
-0.406813627254509 -0.622858050073392
-0.402805611222445 -0.528638956237337
-0.398797595190381 -0.432925039757318
-0.394789579158317 -0.33596068691733
-0.390781563126253 -0.238027083799957
-0.386773547094188 -0.139441690098495
-0.382765531062124 -0.0405570985404816
-0.37875751503006 0.0582407151119272
-0.374749498997996 0.156534737786921
-0.370741482965932 0.253879838382737
-0.366733466933868 0.349806230324476
-0.362725450901804 0.443823439478243
-0.35871743486974 0.535424732416578
-0.354709418837675 0.624091953581325
-0.350701402805611 0.709300714153272
-0.346693386773547 0.790525870474841
-0.342685370741483 0.867247225734546
-0.338677354709419 0.938955385350872
-0.334669338677355 1.00515769411808
-0.330661322645291 1.06538418171568
-0.326653306613227 1.11919344264377
-0.322645290581162 1.16617837702415
-0.318637274549098 1.20597171998699
-0.314629258517034 1.23825128952005
-0.31062124248497 1.26274488565588
-0.306613226452906 1.27923477766915
-0.302605210420842 1.28756172049658
-0.298597194388778 1.28762844681675
-0.294589178356713 1.2794025870671
-0.290581162324649 1.26291897605868
-0.286573146292585 1.23828131169509
-0.282565130260521 1.20566313852889
-0.278557114228457 1.16530813641028
-0.274549098196393 1.11752970221009
-0.270541082164329 1.06270982044411
-0.266533066132265 1.00129722649807
-0.2625250501002 0.933804873964224
-0.258517034068136 0.860806725264821
-0.254509018036072 0.782933892170698
-0.250501002004008 0.700870159944518
-0.246492985971944 0.615346935572145
-0.24248496993988 0.527137666821136
-0.238476953907816 0.43705178461776
-0.234468937875752 0.345928226404674
-0.230460921843687 0.254628602678741
-0.226452905811623 0.164030072768479
-0.222444889779559 0.0750179990567077
-0.218436873747495 -0.011521548742363
-0.214428857715431 -0.094709366001147
-0.210420841683367 -0.173680688174803
-0.206412825651303 -0.24759269417007
-0.202404809619239 -0.315631864566768
-0.198396793587174 -0.377021122618362
-0.19438877755511 -0.431026688145147
-0.190380761523046 -0.476964577333593
-0.186372745490982 -0.514206685033178
-0.182364729458918 -0.542186390356012
-0.178356713426854 -0.560403631186784
-0.17434869739479 -0.568429398547811
-0.170340681362726 -0.565909607577703
-0.166332665330661 -0.552568308109774
-0.162324649298597 -0.528210204411512
-0.158316633266533 -0.492722460499814
-0.154308617234469 -0.446075774507047
-0.150300601202405 -0.388324712767242
-0.146292585170341 -0.319607301546434
-0.142284569138277 -0.240143881582756
-0.138276553106212 -0.15023523775759
-0.134268537074148 -0.0502600232177339
-0.130260521042084 0.0593284959589565
-0.12625250501002 0.178006342981901
-0.122244488977956 0.305183340752379
-0.118236472945892 0.440208186045251
-0.114228456913828 0.582374037815071
-0.110220440881764 0.730924565246765
-0.1062124248497 0.885060397004173
-0.102204408817635 1.04394590956045
-0.0981963927855712 1.20671628955394
-0.0941883767535071 1.37248480282655
-0.090180360721443 1.54035020118699
-0.0861723446893788 1.70940419701014
-0.0821643286573147 1.87873893554052
-0.0781563126252506 2.0474543952109
-0.0741482965931864 2.21466564740642
-0.0701402805611223 2.37950990888656
-0.0661322645290582 2.54115332249779
-0.0621242484969941 2.69879740484333
-0.0581162324649299 2.85168510318701
-0.0541082164328658 2.99910640801886
-0.0501002004008017 3.14040347235478
-0.0460921843687375 3.27497519393461
-0.0420841683366734 3.40228122196762
-0.0380761523046093 3.5218453558979
-0.0340681362725451 3.63325830976322
-0.030060120240481 3.7361798220402
-0.0260521042084169 3.83034009734181
-0.0220440881763528 3.91554057289668
-0.0180360721442887 3.99165400932912
-0.0140280561122246 4.05862391180917
-0.0100200400801604 4.11646329409036
-0.00601202404809631 4.16525280423661
-0.00200400801603218 4.20513823689893
0.00200400801603196 4.23632746277909
0.00601202404809609 4.25908681135736
0.0100200400801602 4.27373694801325
0.0140280561122244 4.28064829128522
0.0180360721442885 4.28023602015512
0.0220440881763526 4.27295472486812
0.0260521042084167 4.25929275787692
0.0300601202404809 4.23976634400316
0.034068136272545 4.21491351081776
0.0380761523046091 4.18528790154108
0.0420841683366733 4.15145253344317
0.0460921843687374 4.11397356478121
0.0501002004008015 4.07341413274906
0.0541082164328657 4.03032832374038
0.0581162324649298 3.98525533545806
0.0621242484969939 3.93871388805889
0.0661322645290581 3.89119693862933
0.0701402805611222 3.84316674987756
0.0741482965931863 3.79505036003497
0.0781563126252505 3.74723549662741
0.0821643286573146 3.70006697204755
0.0861723446893787 3.65384359378414
0.0901803607214429 3.60881561679151
0.094188376753507 3.56518275987048
0.0981963927855711 3.52309280213328
0.102204408817635 3.48264076970122
0.106212424849699 3.44386871679102
0.110220440881764 3.40676609934582
0.114228456913828 3.37127073341703
0.118236472945892 3.33727032466386
0.122244488977956 3.3046045496641
0.12625250501002 3.27306766427952
0.130260521042084 3.24241160914336
0.134268537074148 3.21234957748756
0.138276553106212 3.18256000604872
0.142284569138276 3.15269094572827
0.14629258517034 3.12236476507123
0.150300601202405 3.09118313650489
0.154308617234469 3.05873225267
0.158316633266533 3.02458821810934
0.162324649298597 2.98832256006741
0.166332665330661 2.94950780121424
0.170340681362725 2.9077230367432
0.174348697394789 2.86255945850634
0.178356713426854 2.81362576963903
0.182364729458918 2.76055343447521
0.186372745490982 2.70300171045151
0.190380761523046 2.64066241111963
0.19438877755511 2.57326435230561
0.198396793587174 2.50057743684032
0.202404809619238 2.42241633709994
0.206412825651302 2.33864373879997
0.210420841683367 2.24917311403415
0.214428857715431 2.1539709963954
0.218436873747495 2.05305873610668
0.222444889779559 1.94651371837379
0.226452905811623 1.83447003359392
0.230460921843687 1.71711859355793
0.234468937875751 1.59470669331309
0.238476953907816 1.46753702385068
0.24248496993988 1.3359661461926
0.246492985971944 1.20040244271772
0.250501002004008 1.06130356663928
0.254509018036072 0.919173415367888
0.258517034068136 0.774558658021573
0.2625250501002 0.628044851531356
0.266533066132264 0.480252183595168
0.270541082164329 0.331830884118781
0.274549098196393 0.183456349716907
0.278557114228457 0.0358240283034595
0.282565130260521 -0.110355887244891
0.286573146292585 -0.254363905929898
0.290581162324649 -0.395476774646149
0.294589178356713 -0.532973051661027
0.298597194388778 -0.666138666009006
0.302605210420842 -0.794272411931674
0.306613226452906 -0.91669132869128
0.31062124248497 -1.03273591777785
0.314629258517034 -1.14177515169478
0.318637274549098 -1.24321123111654
0.322645290581162 -1.33648405023208
0.326653306613226 -1.42107533347989
0.330661322645291 -1.49651241060515
0.334669338677354 -1.56237160097944
0.338677354709419 -1.61828118237205
0.342685370741483 -1.66392392379778
0.346693386773547 -1.69903916663693
0.350701402805611 -1.72342444287646
0.354709418837675 -1.73693662400095
0.358717434869739 -1.73949259871704
0.362725450901803 -1.73106948226932
0.366733466933868 -1.71170436455034
0.370741482965932 -1.68149360847091
0.374749498997996 -1.64059171409408
0.37875751503006 -1.58920976780206
0.382765531062124 -1.52761349922239
0.386773547094188 -1.45612097175014
0.390781563126252 -1.37509993523965
0.394789579158316 -1.2849648717739
0.398797595190381 -1.18617376733611
0.402805611222445 -1.07922464368978
0.406813627254509 -0.964651885815469
0.410821643286573 -0.843022400851669
0.414829659318637 -0.714931644649647
0.418837675350701 -0.580999551788048
0.422845691382765 -0.441866404220219
0.42685370741483 -0.298188672667867
0.430861723446894 -0.150634863457449
0.434869739478958 0.000118598246333201
0.438877755511022 0.15339141988439
0.442885771543086 0.308503400972809
0.44689378757515 0.464778276420136
0.450901803607214 0.62154742214008
0.454909819639278 0.77815340413512
0.458917835671343 0.933953355005733
0.462925851703407 1.0883221646193
0.466933867735471 1.24065547440189
0.470941883767535 1.39037246734597
0.474949899799599 1.53691844831095
0.478957915831663 1.67976721148805
0.482965931863727 1.8184231939673
0.486973947895792 1.95242341614708
0.490981963927856 2.08133921123866
0.49498997995992 2.20477774731513
0.498997995991984 2.32238334622232
0.503006012024048 2.4338386041998
0.507014028056112 2.53886531925291
0.511022044088176 2.63722523018011
0.51503006012024 2.72872057170995
0.519038076152305 2.81319444946281
0.523046092184369 2.89053103745749
0.527054108216433 2.96065559967125
0.531062124248497 3.02353433578247
0.535070140280561 3.07917404973154
0.539078156312625 3.12762163818848
0.543086172344689 3.16896339447986
0.547094188376753 3.20332412207216
0.551102204408817 3.23086605040527
0.555110220440882 3.25178754479026
0.559118236472946 3.26632160130298
0.56312625250501 3.27473411718922
0.567134268537074 3.2773219273141
0.571142284569138 3.27441059769986
0.575150300601202 3.26635196825549
0.579158316633266 3.25352143845272
0.58316633266533 3.23631499198099
0.587174348697395 3.21514595933866
0.591182364729459 3.19044152089678
0.595190380761523 3.16263895719688
0.599198396793587 3.13218165808907
0.603206412825651 3.09951490773859
0.607214428857715 3.06508146846536
0.611222444889779 3.02931699275036
0.615230460921844 2.99264529944461
0.619238476953908 2.95547355713013
0.623246492985972 2.91818742457024
0.627254509018036 2.88114620509139
0.6312625250501 2.84467807838993
0.635270541082164 2.80907547946767
0.639278557114228 2.77459069997349
0.643286573146292 2.74143179195852
0.647294589178357 2.70975885772943
0.651302605210421 2.67968081189739
0.655310621242485 2.65125270266125
0.659318637274549 2.62447367863547
0.663326653306613 2.5992856849501
0.667334669338677 2.57557296774912
0.671342685370741 2.55316245945456
0.675350701402806 2.53182510814095
0.67935871743487 2.51127820300925
0.683366733466934 2.4911887342367
0.687374749498998 2.47117780943559
0.691382765531062 2.45082613066162
0.695390781563126 2.42968051551446
0.69939879759519 2.40726142357916
0.703406813627254 2.3830714255449
0.707414829659319 2.35660452715957
0.711422845691383 2.32735623415841
0.715430861723447 2.29483421794202
0.719438877755511 2.25856941564406
0.723446893787575 2.21812737296139
0.727454909819639 2.17311961442292
0.731462925851703 2.12321480440532
0.735470941883767 2.06814944396939
0.739478957915832 2.0077378343243
0.743486973947896 1.94188102828282
0.74749498997996 1.8705744872998
0.751503006012024 1.7939141644162
0.755511022044088 1.71210074344543
0.759519038076152 1.62544178274315
0.763527054108216 1.53435153850239
0.76753507014028 1.43934827817888
0.771543086172344 1.34104893967887
0.775551102204409 1.24016104642346
0.779559118236473 1.13747185219325
0.783567134268537 1.03383476232402
0.787575150300601 0.93015315862307
0.791583166332665 0.827361843208586
0.795591182364729 0.726406409857589
0.799599198396793 0.628220948484087
0.803607214428858 0.533704586723887
0.807615230460922 0.443697469488347
0.811623246492986 0.35895686951964
0.81563126252505 0.280134205734321
0.819639278557114 0.207753817346076
0.823647294589178 0.142194395881299
0.827655310621242 0.0836740093719085
0.831663326653306 0.0322396580991388
0.835671342685371 -0.0122377260107466
0.839679358717435 -0.0500619891653848
0.843687374749499 -0.0817034351998858
0.847695390781563 -0.107782890580289
0.851703406813627 -0.129047777634115
0.855711422845691 -0.146340080501364
0.859719438877755 -0.160556406822211
0.86372745490982 -0.172600723629277
0.867735470941884 -0.183330776285447
0.871743486973948 -0.193499679861612
0.875751503006012 -0.20369469499191
0.879759519038076 -0.214275752845948
0.88376753507014 -0.225316859478068
0.887775551102204 -0.236554065855148
0.891783567134268 -0.247344207139439
0.895791583166333 -0.256639056521056
0.899799599198397 -0.262979859568119
0.903807615230461 -0.264517359331161
0.907815631262525 -0.259062323824532
0.911823647294589 -0.244171167101907
0.915831663326653 -0.217270420237764
0.919839679358717 -0.175822451206068
0.923847695390781 -0.11753282835606
0.927855711422846 -0.0405969282970916
0.93186372745491 0.0560203565790474
0.935871743486974 0.172282835899493
0.939879759519038 0.306776409146314
0.943887775551102 0.456422655621738
0.947895791583166 0.616247265707423
0.95190380761523 0.779249655943142
0.955911823647294 0.93643421145995
0.959919839679358 1.07708042381229
0.963927855711423 1.18934908976766
0.967935871743487 1.261345078745
0.971943887775551 1.28278435467365
0.975951903807615 1.24744435970003
0.979959919839679 1.15661295169443
0.983967935871743 1.02379225980927
0.987975951903807 0.880960505182693
0.991983967935872 0.786747438204277
0.995991983967936 0.836937956688997
1 1.17778403954744
};
\end{axis}

\end{tikzpicture}

%% file: figs/regression1d_kldiv.tikz
% This file was created by matplotlib2tikz v0.6.11.
\begin{tikzpicture}

\definecolor{color0}{rgb}{0.12156862745098,0.466666666666667,0.705882352941177}
\definecolor{color1}{rgb}{1,0.498039215686275,0.0549019607843137}
\definecolor{color2}{rgb}{0.172549019607843,0.627450980392157,0.172549019607843}
\definecolor{color3}{rgb}{0.83921568627451,0.152941176470588,0.156862745098039}
\definecolor{color4}{rgb}{0.580392156862745,0.403921568627451,0.741176470588235}
\definecolor{color5}{rgb}{0.890196078431372,0.466666666666667,0.76078431372549}

\begin{axis}[
xlabel={Dimension of the reduced space},
ylabel={Kullback-Leibler divergence (nat)},
xmin=-1.45, xmax=52.45,
ymin=0.01, ymax=200,
ymode=log,
width=\figurewidth,
height=\figureheight,
tick align=outside,
tick pos=left,
x grid style={lightgray!92.026143790849673!black},
y grid style={lightgray!92.026143790849673!black},
legend style={at={(0.03,0.03)}, anchor=south west, draw=none},
legend cell align={left},
legend entries={{PCA-A},{PCA-Y},{PCA-YN},{KLD},{EKLD},{MI}}
]
\addplot [semithick, color0, mark=*, mark size=3, mark repeat=1, mark options={solid,fill opacity=0}]
table {%
1 52.0358386309258
5 43.7298426898949
10 33.9242449410232
15 25.6312491689055
20 15.2370329056399
25 5.5847119150925
26 4.78575995003329
27 3.954590746106
28 3.18701113902789
29 2.29276075016026
30 1.25142855196469
35 1.37192557227499
40 1.41934711188582
45 1.4679717993701
50 1.59058985110306
};
\addplot [semithick, color1, mark=square, mark size=3, mark repeat=1, mark options={solid,fill opacity=0}]
table {%
1 52.039548774305
5 43.7461094649654
10 33.900896289634
15 25.6414485413432
20 15.0642673107981
25 5.64282067221071
26 4.99811218830575
27 4.03090101320679
28 3.20507330232947
29 2.32295811039825
30 2.16562890669619
35 1.81665723088621
40 1.71215706291876
45 1.09881348756021
50 1.03209873647133
};
\addplot [semithick, color2, mark=star, mark size=3, mark repeat=1, mark options={solid,fill opacity=0}]
table {%
1 51.4986356715735
5 42.3283400250594
10 32.9920141004972
15 23.590573125355
20 17.8267368252751
25 8.63901375700808
26 7.71444725266219
27 4.44530878858151
28 2.97559765568452
29 2.43076853429167
30 1.25142855196279
35 1.24953837811446
40 1.27857359545475
45 1.19638562492608
50 1.22945395422156
};
\addplot [semithick, color3, mark=x, mark size=3, mark repeat=1, mark options={solid,fill opacity=0}]
table {%
1 40.5859244649867
5 30.1709507451799
10 20.8899622176879
15 13.3955910910769
20 7.30272581011737
25 2.558476164995
26 1.85693354371982
27 1.21183653266253
28 0.672801481080481
29 0.25525181827991
30 3.72908321437435e-24
35 4.57362223191873e-23
40 1.7763567900051e-15
45 3.55271372145608e-15
50 1.77635687785735e-15
};
\addplot [semithick, color4, mark=diamond, mark size=3, mark repeat=1, mark options={solid,fill opacity=0}]
table {%
1 50.2721481595732
5 40.5104817303461
10 31.1236930072281
15 22.2396011058365
20 15.3916643239805
25 7.35810517200525
26 6.18139082027454
27 4.4212599200976
28 2.34805732057005
29 0.876041326782407
30 1.7763569062705e-15
35 1.12489995994254e-22
40 1.02685807437731e-22
45 3.55271376968763e-15
50 9.77472114368791e-23
};
\addplot [semithick, color5, mark=asterisk, mark size=3, mark repeat=1, mark options={solid,fill opacity=0}]
table {%
1 50.2721481595728
5 40.5104817303478
10 31.1236930072233
15 22.2396011058234
20 15.3916643239794
25 7.35810517200504
26 6.18139082027638
27 4.42125992010101
28 2.34805732057232
29 0.876041326781001
30 1.09921453376584e-21
35 1.17491973826026e-21
40 1.77635798987309e-15
45 1.77635798044768e-15
50 1.77635798352567e-15
};
\end{axis}

\end{tikzpicture}

%% file: figs/regression1d_exp_kldiv.tikz
% This file was created by matplotlib2tikz v0.6.11.
\begin{tikzpicture}

\definecolor{color0}{rgb}{0.12156862745098,0.466666666666667,0.705882352941177}
\definecolor{color1}{rgb}{1,0.498039215686275,0.0549019607843137}
\definecolor{color2}{rgb}{0.172549019607843,0.627450980392157,0.172549019607843}
\definecolor{color3}{rgb}{0.83921568627451,0.152941176470588,0.156862745098039}
\definecolor{color4}{rgb}{0.580392156862745,0.403921568627451,0.741176470588235}
\definecolor{color5}{rgb}{0.890196078431372,0.466666666666667,0.76078431372549}

\begin{axis}[
xlabel={Dimension of the reduced space},
ylabel={Expected Kullback-Leibler divergence (nat)},
xmin=-1.45, xmax=52.45,
ymin=0.01, ymax=200,
ymode=log,
width=\figurewidth,
height=\figureheight,
tick align=outside,
tick pos=left,
x grid style={lightgray!92.026143790849673!black},
y grid style={lightgray!92.026143790849673!black},
legend style={at={(0.03,0.03)}, anchor=south west, draw=none},
legend cell align={left},
legend entries={{PCA-A},{PCA-Y},{PCA-YN},{KLD},{EKLD},{MI}}
]
\addplot [semithick, color0, mark=*, mark size=3, mark repeat=1, mark options={solid,fill opacity=0}]
table {%
1 53.078639858062
5 43.8832459719521
10 32.5714378947297
15 22.8088781506219
20 14.1471781933401
25 6.90314055050326
26 5.6761263133566
27 4.46476802638121
28 3.42775777046691
29 2.36208883571927
30 1.62208039224254
35 1.55582494658662
40 1.52971351506177
45 1.48574308319776
50 1.43204953943753
};
\addplot [semithick, color1, mark=square, mark size=3, mark repeat=1, mark options={solid,fill opacity=0}]
table {%
1 53.083716153486
5 43.9061783611184
10 32.6314961144896
15 22.821227554135
20 14.1802710554842
25 6.96103077954502
26 5.86628449841655
27 4.5625952537149
28 3.48341468438266
29 2.40112738208661
30 2.36462787202628
35 2.19598029325216
40 2.09102963767975
45 1.46712454393009
50 1.42418495766648
};
\addplot [semithick, color2, mark=star, mark size=3, mark repeat=1, mark options={solid,fill opacity=0}]
table {%
1 52.8574240325398
5 41.8987595380283
10 30.6405757036471
15 21.0380170666636
20 12.89003308012
25 6.12359998772266
26 5.01394254192112
27 4.01439652278095
28 3.1472749671088
29 2.34295738142966
30 1.62208039224256
35 1.58615371962671
40 1.48872395413579
45 1.43731323822989
50 1.38900492879602
};
\addplot [semithick, color3, mark=x, mark size=3, mark repeat=1, mark options={solid,fill opacity=0}]
table {%
1 54.3399635211417
5 41.9312906928936
10 30.1744134286089
15 20.229569295286
20 11.7226662947568
25 4.6276623036421
26 3.47686426909713
27 2.38955569743304
28 1.42396005865524
29 0.604157313565674
30 1.51518809709854e-20
35 2.09057565732528e-20
40 1.77634220594683e-15
45 3.55273483822448e-15
50 1.77637338122167e-15
};
\addplot [semithick, color4, mark=diamond, mark size=3, mark repeat=1, mark options={solid,fill opacity=0}]
table {%
1 51.8759130427854
5 40.444951704132
10 29.1564370116286
15 19.4204621667652
20 11.1564163991571
25 4.36609426540769
26 3.24379873802834
27 2.22487642110508
28 1.33043010594222
29 0.551847249678435
30 1.77636997241235e-15
35 2.10792348856889e-20
40 1.50117086723566e-20
45 3.55273379045367e-15
50 1.35191978226287e-20
};
\addplot [semithick, color5, mark=asterisk, mark size=3, mark repeat=1, mark options={solid,fill opacity=0}]
table {%
1 51.8759130427853
5 40.444951704132
10 29.1564370116285
15 19.4204621667652
20 11.1564163991571
25 4.36609426540769
26 3.24379873802834
27 2.22487642110508
28 1.33043010594222
29 0.55184724967843
30 2.0016709652439e-20
35 2.00128129609992e-20
40 1.77637697607348e-15
45 1.77637697630709e-15
50 1.77637697661899e-15
};
\end{axis}

\end{tikzpicture}

%% file: figs/regression1d_mi.tikz
% This file was created by matplotlib2tikz v0.6.11.
\begin{tikzpicture}

\definecolor{color0}{rgb}{0.12156862745098,0.466666666666667,0.705882352941177}
\definecolor{color1}{rgb}{1,0.498039215686275,0.0549019607843137}
\definecolor{color2}{rgb}{0.172549019607843,0.627450980392157,0.172549019607843}
\definecolor{color3}{rgb}{0.83921568627451,0.152941176470588,0.156862745098039}
\definecolor{color4}{rgb}{0.580392156862745,0.403921568627451,0.741176470588235}
\definecolor{color5}{rgb}{0.890196078431372,0.466666666666667,0.76078431372549}

\begin{axis}[
xlabel={Dimension of the reduced space},
ylabel={Relative error on the mutual information},
xmin=-1.5, xmax=53.5,
ymin=0.001, ymax=2,
ymode=log,
width=\figurewidth,
height=\figureheight,
tick align=outside,
tick pos=left,
x grid style={lightgray!92.026143790849673!black},
y grid style={lightgray!92.026143790849673!black},
legend style={at={(0.03,0.03)}, anchor=south west, draw=none},
legend cell align={left},
legend entries={{PCA-A},{PCA-Y},{PCA-YN},{KLD},{EKLD},{MI},{Error estim.}}
]
\addplot [semithick, color0, mark=*, mark size=3, mark repeat=1, mark options={solid,fill opacity=0}]
table {%
1 0.956940183176236
5 0.791158958689562
10 0.5872214854909
15 0.411214983915468
20 0.255055580322829
25 0.124454820254148
26 0.102333318422236
27 0.0804940734052062
28 0.0617981010346945
29 0.0425854492344578
30 0.0292440407630548
35 0.0280495395762876
40 0.0275787837669235
45 0.02678605295779
50 0.0258180268414707
};
\addplot [semithick, color1, mark=square, mark size=3, mark repeat=1, mark options={solid,fill opacity=0}]
table {%
1 0.957031702308712
5 0.791572400419572
10 0.588304258598349
15 0.411437627911084
20 0.255652202422542
25 0.12549850725388
26 0.105761627981981
27 0.0822577735507769
28 0.0628015242109464
29 0.0432892644379277
30 0.0426312248209301
35 0.0395907240591179
40 0.0376985975872145
45 0.0264503844398127
50 0.025676238462174
};
\addplot [semithick, color2, mark=star, mark size=3, mark repeat=1, mark options={solid,fill opacity=0}]
table {%
1 0.952951943968103
5 0.755381199186516
10 0.552410502696952
15 0.379288617026705
20 0.232390857222515
25 0.110400698088694
26 0.0903949895345904
27 0.0723744495733694
28 0.056741353801887
29 0.0422405335127809
30 0.0292440407630603
35 0.0285963286746189
40 0.0268397942591642
45 0.0259129246176671
50 0.0250419874082332
};
\addplot [semithick, color3, mark=x, mark size=3, mark repeat=1, mark options={solid,fill opacity=0}]
table {%
1 0.979680240201431
5 0.755967694420282
10 0.544006191394666
15 0.364713334756801
20 0.211344722875274
25 0.0834308494781247
26 0.0626834285773278
27 0.0430806417214551
28 0.0256721838197994
29 0.0108921858554042
30 4.91910886661858e-13
35 4.69621174609992e-13
40 4.87811629273009e-13
45 4.88324036446615e-13
50 4.954977368771e-13
};
\addplot [semithick, color4, mark=diamond, mark size=3, mark repeat=1, mark options={solid,fill opacity=0}]
table {%
1 0.935256552585799
5 0.72917089804474
10 0.525653374865448
15 0.350126165118582
20 0.20113595942036
25 0.0787151113376749
26 0.0584815542906213
27 0.0401116843917738
28 0.0239859580551218
29 0.00994910211677979
30 4.79356910908508e-13
35 4.84480982644569e-13
40 4.9293570100907e-13
45 4.75513857106462e-13
50 4.72951821238432e-13
};
\addplot [semithick, color5, mark=asterisk, mark size=3, mark repeat=1, mark options={solid,fill opacity=0}]
table {%
1 0.935256552585811
5 0.729170898044739
10 0.525653374865472
15 0.350126165118585
20 0.201135959420371
25 0.0787151113376885
26 0.0584815542906207
27 0.040111684391793
28 0.0239859580551263
29 0.00994910211677261
30 4.80637928842523e-13
35 4.81406539602932e-13
40 4.80125521668917e-13
45 4.8038172525572e-13
50 4.79869318082114e-13
};
\addplot [ultra thick, color0, dashed]
table {%
1 0.935256552585736
2 0.877819861109529
3 0.824926016022053
4 0.77567069996916
5 0.729170898044415
6 0.685003718422537
7 0.642857014894188
8 0.602419930691374
9 0.563481874971819
10 0.525653374864905
11 0.488518028590322
12 0.452342884400706
13 0.41722122930501
14 0.383086668659192
15 0.350126165117816
16 0.318334452514397
17 0.287544751071998
18 0.2577800382252
19 0.228942333241949
20 0.201135959419424
21 0.174316572608321
22 0.148592987815483
23 0.12396158788655
24 0.100653311390432
25 0.0787151113365959
26 0.0584815542895043
27 0.0401116843906554
28 0.0239859580539695
29 0.00994910211559985
30 1.66489044772788e-12
31 3.548272786702e-12
32 4.05409039672122e-12
33 4.40047998040427e-12
34 4.52882176205094e-12
35 4.6349590832051e-12
36 4.72843986187854e-12
37 4.8123727225402e-12
38 4.88342699611621e-12
39 4.94959628838387e-12
40 5.01043651013333e-12
41 5.06505948294489e-12
42 5.1183501881269e-12
43 5.17008658107443e-12
44 5.21738208192346e-12
45 5.26134691369862e-12
46 5.30420152244915e-12
47 5.34505772975535e-12
48 5.38458166943201e-12
49 5.42343947529389e-12
50 5.46007683510652e-12
51 5.4944937488699e-12
};
\end{axis}

\end{tikzpicture}

%% file: figs/regression1d_kle.tikz
% This file was created by matplotlib2tikz v0.6.11.
\begin{tikzpicture}

\definecolor{color0}{rgb}{0.12156862745098,0.466666666666667,0.705882352941177}
\definecolor{color1}{rgb}{1,0.498039215686275,0.0549019607843137}
\definecolor{color2}{rgb}{0.172549019607843,0.627450980392157,0.172549019607843}
\definecolor{color3}{rgb}{0.890196078431372,0.466666666666667,0.76078431372549}

\begin{groupplot}[group style={group size=1 by 3}]
\nextgroupplot[
ylabel={Kullback-Leibler divergence (nat)},
xmin=-23.95, xmax=524.95,
ymin=0.01, ymax=200,
ymode=log,
width=\figurewidth,
height=\figureheight,
tick align=outside,
tick pos=left,
x grid style={lightgray!92.026143790849673!black},
y grid style={lightgray!92.026143790849673!black},
legend style={draw=none},
legend entries={{PCA-A},{PCA-Y},{PCA-YN},{MI}},
legend cell align={left}
]
\addplot [semithick, color0, mark=*, mark size=3, mark repeat=2, mark options={solid,fill opacity=0}]
table {%
1 52.0358386309258
10 33.9242449410232
20 15.2370329056399
30 1.25142855196469
40 1.41934711188582
50 1.59058985110306
60 1.91918203658925
70 1.22806492717682
80 0.800517914521256
90 0.681159351622375
100 0.58406400802341
110 0.434921381262573
120 0.374938387409697
130 0.364002438951144
140 0.301766221675787
150 0.276926082200592
160 0.251577206049176
170 0.514152731029046
180 0.698474922428877
190 0.874719083546259
200 0.708548625582065
210 1.043663314075
220 0.920506402328247
230 0.572984936545356
240 0.751795835261793
250 0.474749790884025
260 0.386679206217727
270 0.38978690616886
280 0.185606818364441
290 0.18434819629734
300 0.11867580790808
310 0.179719218282847
320 0.395020554236375
330 0.487802618745304
340 0.895952384258498
350 1.17960037823494
360 1.832483274896
370 1.5254684480669
380 1.23393868639128
390 0.973021257810495
400 0.796839548293433
410 0.581122311424138
420 0.310779330939079
430 0.241335746006401
440 0.252489619659613
450 0.0426015856771766
460 0.0326542903493064
470 0.036628255545369
480 0.000607444481751443
490 0.000537421115572072
500 1.77635567453502e-15
};
\addplot [semithick, color1, mark=square, mark size=3, mark repeat=2, mark options={solid,fill opacity=0}]
table {%
1 52.039548774305
10 33.900896289634
20 15.0642673107981
30 2.16562890669619
40 1.71215706291876
50 1.03209873647133
60 0.691177747247622
70 0.559437752359517
80 0.458091195084912
90 1.00043331463835
100 0.960563216136441
110 0.407997673205888
120 0.345998263944622
130 0.369234607706684
140 0.352549753861303
150 0.138458865310467
160 0.131133272452354
170 0.143020453687925
180 0.0422218174833614
190 0.0373783681756879
200 0.0251645875597607
210 0.0310131877573026
220 0.0310426877615362
230 0.0308745739445044
240 0.030828627010866
250 0.0292977323318661
260 0.029849468300031
270 0.0172226141895207
280 0.0177248982636202
290 0.0178896626685702
300 0.0168295397860263
310 0.01641308631789
320 0.0155847192366205
330 0.0159521601424188
340 0.0042146077500937
350 0.00275606990567076
360 0.00244797303406431
370 0.00227929176785605
380 0.00225028433061583
390 0.00231806198231094
400 0.00224357234258867
410 0.00208313940329516
420 0.0021163379888128
430 0.00211203781359992
440 0.00172076992583501
450 0.00170783749478335
460 0.00172458515396674
470 0.00160369135662565
480 2.2544959727873e-05
490 6.45685843624761e-06
500 1.77635514275177e-15
};
\addplot [semithick, color2, mark=star, mark size=3, mark repeat=2, mark options={solid,fill opacity=0}]
table {%
1 51.4986356715735
10 32.9920141004972
20 17.8267368252751
30 1.25142855196279
40 1.27857359545475
50 1.22945395422156
60 1.13696258868259
70 1.41318102269377
80 1.17219109987865
90 1.16214669132409
100 1.39997761117858
110 1.12862661969095
120 0.975636836026936
130 0.915156969456559
140 0.762864794577053
150 0.751712066557618
160 0.585767976006883
170 0.246346559224935
180 0.189831939764716
190 0.182323159059952
200 0.26781381292907
210 0.175057777112323
220 0.143929270768837
230 0.0967118422792841
240 0.0891067292654855
250 0.052194015241382
260 0.0473181110568655
270 0.0659276149908775
280 0.0673960201079976
290 0.045490305787694
300 0.0792869241869611
310 0.0458620092374733
320 0.0265484216241576
330 0.0148114540202929
340 0.0561552455479284
350 0.046885315300868
360 0.0123944476063035
370 0.0042745834902162
380 0.00663490788010387
390 0.00354719730339702
400 0.00738911607024798
410 0.00460732889784044
420 0.0112000190439903
430 0.00960915122545832
440 0.00682757118080655
450 0.00303409578334343
460 0.00203316837466631
470 0.00825641138093096
480 0.00215993335247927
490 2.51959137837261e-05
500 9.03098109347466e-15
};
\addplot [semithick, color3, mark=asterisk, mark size=3, mark repeat=2, mark options={solid,fill opacity=0}]
table {%
1 50.2721481595728
10 31.1236930072233
20 15.3916643239794
30 1.09921453376584e-21
40 1.77635798987309e-15
50 1.77635798352567e-15
60 1.77635799176536e-15
70 1.77635799345484e-15
80 1.77635798587777e-15
90 1.13889860631921e-21
100 1.7763557069716e-15
110 1.77635797577786e-15
120 3.5527148178597e-15
130 1.77635796833466e-15
140 1.77635796723751e-15
150 1.134597456866e-21
160 1.13123745532235e-21
170 1.12630718617805e-21
180 3.55271480801311e-15
190 1.77635571092764e-15
200 1.77635795941178e-15
210 1.12630314986481e-21
220 1.77635797324803e-15
230 1.1320265665026e-21
240 1.77635797011219e-15
250 1.77635796887517e-15
260 1.77635570556131e-15
270 1.77635796802172e-15
280 1.77635796506036e-15
290 1.77635796896536e-15
300 1.13046163436852e-21
310 1.77635797593404e-15
320 1.77635797035237e-15
330 1.77635571614861e-15
340 1.77635570194846e-15
350 1.77635797596677e-15
360 1.14581556552132e-21
370 1.77635799180074e-15
380 1.77635798258122e-15
390 1.77635568886684e-15
400 1.77635569841779e-15
410 1.77635797967779e-15
420 3.55271480894777e-15
430 1.77635570742124e-15
440 1.776357973559e-15
450 1.13139878734454e-21
460 1.77635794726692e-15
470 1.13835414541771e-21
480 1.77635798131412e-15
490 1.77635794537091e-15
500 1.77635501019176e-15
};
\nextgroupplot[
ylabel={Expected Kullback-Leibler divergence(nat)},
xmin=-23.95, xmax=524.95,
ymin=0.01, ymax=200,
ymode=log,
width=\figurewidth,
height=\figureheight,
tick align=outside,
tick pos=left,
x grid style={lightgray!92.026143790849673!black},
y grid style={lightgray!92.026143790849673!black},
legend style={draw=none},
legend entries={{PCA-A},{PCA-Y},{PCA-YN},{MI}},
legend cell align={left}
]
\addplot [semithick, color0, mark=*, mark size=3, mark repeat=2, mark options={solid,fill opacity=0}]
table {%
1 53.078639858062
10 32.5714378947297
20 14.1471781933401
30 1.62208039224254
40 1.52971351506177
50 1.43204953943753
60 1.31393203325197
70 1.25858099718159
80 1.14256960248123
90 1.06769127336039
100 1.01088207855983
110 0.906330434599143
120 0.849110033784864
130 0.798376842919811
140 0.753852648844379
150 0.719797746300175
160 0.65449713373579
170 0.621965729741488
180 0.593924770653357
190 0.577783214348596
200 0.544384719199595
210 0.515339581689923
220 0.480468471118182
230 0.454771460676912
240 0.417830454737426
250 0.402754220323054
260 0.393623508080477
270 0.356260310330645
280 0.315376349441795
290 0.302147580461026
300 0.283144167888962
310 0.262170557069792
320 0.246385311475523
330 0.231763182838168
340 0.208637727277316
350 0.188959469310745
360 0.173498884676824
370 0.159422643138777
380 0.147908027102844
390 0.137743022008211
400 0.125840920900986
410 0.114050125987795
420 0.0956243468754603
430 0.0741400715017861
440 0.0667005354317694
450 0.0511718678832028
460 0.0431623510647149
470 0.0307236994822012
480 0.00975589179810281
490 0.0029413318760716
500 1.77633698082959e-15
};
\addplot [semithick, color1, mark=square, mark size=3, mark repeat=2, mark options={solid,fill opacity=0}]
table {%
1 53.083716153486
10 32.6314961144896
20 14.1802710554842
30 2.36462787202628
40 2.09102963767975
50 1.42418495766648
60 1.25581880467298
70 1.09554197926079
80 0.844215930611212
90 0.683732003359081
100 0.667116770301505
110 0.584172287247303
120 0.567575686713831
130 0.557563830987182
140 0.554505611092866
150 0.384281855125399
160 0.380095164638403
170 0.372201792469691
180 0.1190251626267
190 0.117371692513669
200 0.0915520841316985
210 0.0895186956547191
220 0.0893006086631121
230 0.0889291005456608
240 0.0887317251130658
250 0.0680591000523737
260 0.0678555466758455
270 0.0662255165530703
280 0.0660202437022851
290 0.0655742451130578
300 0.0654778095181733
310 0.0653289509210375
320 0.0643561205759832
330 0.0642883820740013
340 0.0414045801439627
350 0.0405347941978045
360 0.0402822123984636
370 0.0401358916605348
380 0.0400586626537627
390 0.0400232288694681
400 0.0399481979812876
410 0.0398511418103802
420 0.0398255565280888
430 0.0397340265301265
440 0.0394408479278861
450 0.0393685195216917
460 0.0393556228687808
470 0.0392286599598473
480 3.02619045462414e-05
490 9.66735137660971e-06
500 1.77633923680494e-15
};
\addplot [semithick, color2, mark=star, mark size=3, mark repeat=2, mark options={solid,fill opacity=0}]
table {%
1 52.8574240325398
10 30.6405757036471
20 12.89003308012
30 1.62208039224256
40 1.48872395413579
50 1.38900492879602
60 1.20491698121956
70 1.09903896926993
80 0.997613795977128
90 0.915694398974256
100 0.818764906759663
110 0.729287146804488
120 0.702438273548823
130 0.661277245540907
140 0.614616591886132
150 0.561220546025067
160 0.522430077141516
170 0.483958042440858
180 0.444499346852045
190 0.416506393530272
200 0.381715195198357
210 0.326740366401233
220 0.309497547245522
230 0.272960781000374
240 0.249835608969364
250 0.232620526897811
260 0.207719738609772
270 0.186188104341983
280 0.165484770700772
290 0.146647596001851
300 0.130992208726981
310 0.120046919919971
320 0.107718189954179
330 0.0981530685830798
340 0.0860174197737027
350 0.0770760736255345
360 0.0701099549905769
370 0.060515765632837
380 0.0553311442931116
390 0.05050984742868
400 0.0465815829156196
410 0.0441485126449466
420 0.0419870177094811
430 0.040540109498515
440 0.0391703330234702
450 0.0362963424448033
460 0.0319493984010061
470 0.0165409929976936
480 0.00592086048722051
490 4.06796523993081e-05
500 4.15710136527008e-12
};
\addplot [semithick, color3, mark=asterisk, mark size=3, mark repeat=2, mark options={solid,fill opacity=0}]
table {%
1 51.8759130427853
10 29.1564370116285
20 11.1564163991571
30 2.0016709652439e-20
40 1.77637697607348e-15
50 1.77637697661899e-15
60 1.77637697666052e-15
70 1.77637697667064e-15
80 1.77637697738099e-15
90 2.01378790783347e-20
100 1.776336701689e-15
110 1.77637697737128e-15
120 3.55273381692421e-15
130 1.77637697760323e-15
140 1.77637697951762e-15
150 2.01401657686103e-20
160 2.01402307177912e-20
170 2.01401554904709e-20
180 3.55273381914091e-15
190 1.77633669938156e-15
200 1.77637697985564e-15
210 2.01402109686321e-20
220 1.77637698004029e-15
230 2.01402578256631e-20
240 1.77637697976806e-15
250 1.77637698011767e-15
260 1.77633669868588e-15
270 1.77637697982339e-15
280 1.77637697978004e-15
290 1.77637698003756e-15
300 2.01403816636749e-20
310 1.7763769800381e-15
320 1.77637698009397e-15
330 1.77633669894005e-15
340 1.77633669861303e-15
350 1.77637697995238e-15
360 2.01408158491184e-20
370 1.7763769800676e-15
380 1.77637698026335e-15
390 1.77633669864827e-15
400 1.77633669825868e-15
410 1.7763769800901e-15
420 3.55273381971539e-15
430 1.77633669830426e-15
440 1.77637698043019e-15
450 2.01414473864554e-20
460 1.77637698102465e-15
470 2.01416096253546e-20
480 1.77637698116038e-15
490 1.77637698196133e-15
500 1.77633653841398e-15
};
\nextgroupplot[
xlabel={Dimension of the reduced space},
ylabel={Absolute error on the entropy (nat)},
xmin=-23.95, xmax=524.95,
ymin=0.01, ymax=200,
ymode=log,
width=\figurewidth,
height=\figureheight,
tick align=outside,
tick pos=left,
x grid style={lightgray!92.026143790849673!black},
y grid style={lightgray!92.026143790849673!black},
legend style={draw=none},
legend entries={{PCA-A},{PCA-Y},{PCA-YN},{MI}},
legend cell align={left}
]
\addplot [semithick, color0, mark=*, mark size=3, mark repeat=2, mark options={solid,fill opacity=0}]
table {%
1 53.0786398580363
10 32.5714378947052
20 14.1471781933152
30 1.62208039221825
40 1.52971351503754
50 1.43204953941322
60 1.3139320332277
70 1.25858099715727
80 1.14256960245678
90 1.06769127333591
100 1.01088207853533
110 0.906330434574517
120 0.849110033760279
130 0.79837684289517
140 0.75385264881972
150 0.719797746275457
160 0.654497133711086
170 0.621965729716688
180 0.593924770628529
190 0.577783214323798
200 0.544384719174598
210 0.515339581664929
220 0.480468471093129
230 0.454771460651937
240 0.417830454712455
250 0.402754220297989
260 0.393623508055352
270 0.356260310305554
280 0.31537634941655
290 0.30214758043584
300 0.283144167863551
310 0.262170557044453
320 0.246385311450066
330 0.231763182812667
340 0.208637727251823
350 0.188959469285134
360 0.173498884651249
370 0.159422643113157
380 0.147908027077161
390 0.137743021982455
400 0.125840920875188
410 0.114050125961992
420 0.0956243468496147
430 0.0741400714757781
440 0.0667005354056549
450 0.0511718678569757
460 0.043162351038486
470 0.0307236994561002
480 0.00975589177202707
490 0.00294133184973688
500 2.63078447915177e-11
};
\addplot [semithick, color1, mark=square, mark size=3, mark repeat=2, mark options={solid,fill opacity=0}]
table {%
1 53.0837161534603
10 32.6314961144651
20 14.1802710554593
30 2.36462787200189
40 2.09102963765529
50 1.42418495764202
60 1.25581880464848
70 1.09554197923628
80 0.844215930586689
90 0.683732003334466
100 0.667116770276802
110 0.584172287222568
120 0.567575686689064
130 0.55756383096239
140 0.554505611068091
150 0.384281855100252
160 0.380095164613277
170 0.372201792444553
180 0.119025162600533
190 0.117371692487492
200 0.0915520841053699
210 0.0895186956283389
220 0.0893006086367301
230 0.0889291005191915
240 0.0887317250866992
250 0.0680591000260016
260 0.0678555466494508
270 0.0662255165266821
280 0.0660202436758617
290 0.065574245086605
300 0.065477809491707
310 0.0653289508946209
320 0.0643561205495438
330 0.0642883820475468
340 0.041404580117554
350 0.0405347941714176
360 0.0402822123720448
370 0.040135891634133
380 0.0400586626273132
390 0.0400232288430367
400 0.0399481979548817
410 0.0398511417839948
420 0.0398255565016647
430 0.0397340265036625
440 0.0394408479014281
450 0.0393685194952802
460 0.0393556228423613
470 0.0392286599334106
480 3.02618781908848e-05
490 9.66732495655265e-06
500 2.63931099198089e-11
};
\addplot [semithick, color2, mark=star, mark size=3, mark repeat=2, mark options={solid,fill opacity=0}]
table {%
1 52.8574240325147
10 30.6405757036209
20 12.8900330800954
30 1.62208039221856
40 1.48872395411114
50 1.38900492877128
60 1.20491698119452
70 1.09903896924474
80 0.997613795951377
90 0.915694398948123
100 0.818764906732721
110 0.729287146776723
120 0.702438273521746
130 0.661277245512931
140 0.614616591857921
150 0.56122054599491
160 0.522430077112812
170 0.483958042411015
180 0.444499346821072
190 0.416506393500679
200 0.381715195167967
210 0.32674036636848
220 0.309497547213123
230 0.272960780967818
240 0.24983560893595
250 0.23262052686297
260 0.20771973857563
270 0.186188104304957
280 0.165484770663976
290 0.146647595963554
300 0.130992208687186
310 0.120046919878629
320 0.107718189915101
330 0.098153068541702
340 0.0860174197339738
350 0.0770760735820311
360 0.0701099549528514
370 0.0605157655898942
380 0.0553311442499016
390 0.0505098473856798
400 0.0465815828742926
410 0.0441485126028986
420 0.0419870176691539
430 0.0405401094554882
440 0.0391703329821347
450 0.0362963423989378
460 0.0319493983561721
470 0.0165409929529794
480 0.00592086042038531
490 4.06795829768214e-05
500 8.37481195503642e-11
};
\addplot [semithick, color3, mark=asterisk, mark size=3, mark repeat=2, mark options={solid,fill opacity=0}]
table {%
1 51.8759130427575
10 29.1564370115988
20 11.1564163991307
30 2.66631161593978e-11
40 2.66346944499674e-11
50 2.66204835952522e-11
60 2.66346944499674e-11
70 2.66204835952522e-11
80 2.66346944499674e-11
90 2.66204835952522e-11
100 2.66489053046826e-11
110 2.66204835952522e-11
120 2.65920618858217e-11
130 2.65494293216761e-11
140 2.66346944499674e-11
150 2.6606272740537e-11
160 2.65778510311065e-11
170 2.66204835952522e-11
180 2.65920618858217e-11
190 2.65352184669609e-11
200 2.66346944499674e-11
210 2.66204835952522e-11
220 2.65920618858217e-11
230 2.66346944499674e-11
240 2.66346944499674e-11
250 2.6606272740537e-11
260 2.65778510311065e-11
270 2.66346944499674e-11
280 2.6606272740537e-11
290 2.66346944499674e-11
300 2.65920618858217e-11
310 2.66346944499674e-11
320 2.66346944499674e-11
330 2.66346944499674e-11
340 2.66204835952522e-11
350 2.66346944499674e-11
360 2.65920618858217e-11
370 2.66346944499674e-11
380 2.66346944499674e-11
390 2.66346944499674e-11
400 2.66489053046826e-11
410 2.66346944499674e-11
420 2.65920618858217e-11
430 2.65494293216761e-11
440 2.66346944499674e-11
450 2.66204835952522e-11
460 2.65920618858217e-11
470 2.66346944499674e-11
480 2.66346944499674e-11
490 2.6606272740537e-11
500 2.65778510311065e-11
};
\end{groupplot}

\end{tikzpicture}

%% file: figs/regression1d_kle_singval.tikz
% This file was created by matplotlib2tikz v0.6.11.
\begin{tikzpicture}

\definecolor{color0}{rgb}{0.12156862745098,0.466666666666667,0.705882352941177}
\definecolor{color1}{rgb}{1,0.498039215686275,0.0549019607843137}
\definecolor{color2}{rgb}{0.172549019607843,0.627450980392157,0.172549019607843}

\begin{axis}[
xlabel={i},
ylabel={$\sigma_i/\sigma_1$},
xmin=0, xmax=61,
ymin=0.01, ymax=1,
ymode=log,
width=\figurewidth,
height=\figureheight,
tick align=outside,
tick pos=left,
x grid style={lightgray!92.026143790849673!black},
y grid style={lightgray!92.026143790849673!black},
legend cell align={left},
legend style={draw=none},
legend entries={{PCA-A},{PCA-Y},{PCA-YN}}
]
\addplot [semithick, color0, mark=*, mark size=3, mark repeat=5, mark options={solid,fill opacity=0}]
table {%
1 1
2 0.918560387555081
3 0.914223430903836
4 0.893407547596022
5 0.884899518954346
6 0.850604089620674
7 0.833857066767184
8 0.802897550748354
9 0.76371167096945
10 0.731318829063651
11 0.70762804898171
12 0.659068978715492
13 0.615522022253694
14 0.586603345490822
15 0.562626079221831
16 0.518576585535896
17 0.477264523375852
18 0.440287051707904
19 0.420025116899271
20 0.393278178370994
21 0.35870698971032
22 0.325250726769756
23 0.293166215121735
24 0.262717894150487
25 0.233206555712198
26 0.20446065107026
27 0.175231612455315
28 0.144811638466864
29 0.110531928099896
30 0.079129033719244
31 1.18847751454166e-08
32 7.26816747148608e-09
33 6.30561539969385e-09
34 5.78960236877383e-09
35 5.25638244144193e-09
36 5.24815594316764e-09
37 5.24465533706174e-09
38 5.23443602513854e-09
39 5.18027746727466e-09
40 4.81914034330007e-09
41 4.76783696797228e-09
42 4.14005736359591e-09
43 4.14000995507465e-09
44 3.96606908780528e-09
45 3.95767112145554e-09
46 3.74372338620647e-09
47 3.66440766808217e-09
48 3.55387334195599e-09
49 3.34188614466351e-09
50 3.26969887178466e-09
51 3.21970308402299e-09
52 2.8450705743999e-09
53 2.77270961320301e-09
54 2.7235039502251e-09
55 2.6666207025505e-09
56 2.65366295778743e-09
57 2.5726580350223e-09
58 2.51425902203134e-09
59 2.48438220966721e-09
60 2.4119350973182e-09
};
\addplot [semithick, color1, mark=square, mark size=3, mark repeat=5, mark options={solid,fill opacity=0}]
table {%
1 1
2 0.914327326176776
3 0.910941421379159
4 0.890026232598142
5 0.882535391280873
6 0.848964423333488
7 0.831064102781969
8 0.802251453595132
9 0.763665453140373
10 0.72954761560577
11 0.708869484860946
12 0.661602119984671
13 0.618412209773045
14 0.58537850450045
15 0.566872548053896
16 0.524282624215375
17 0.484001244936811
18 0.446514281445748
19 0.422561830312161
20 0.401512958134344
21 0.368935563431302
22 0.336599660518743
23 0.305408823637516
24 0.275617814276282
25 0.246377255927045
26 0.217386598095035
27 0.186946555426282
28 0.154867547331352
29 0.116486992207093
30 0.0842813284222298
31 0.0631469448902233
32 0.0617499958317237
33 0.0570103426856744
34 0.0557894357520224
35 0.0524446917496063
36 0.0513604863204017
37 0.0487508335877411
38 0.047778666550796
39 0.0456437830234955
40 0.0447657513864084
41 0.0429683640455591
42 0.0421706666231716
43 0.0406271968996187
44 0.0398988376104658
45 0.0385537811671212
46 0.0378858089094734
47 0.0367000913831162
48 0.0360851054784609
49 0.0350300399256257
50 0.0344618468676688
51 0.0335156934280253
52 0.0329890581613967
53 0.0321349262308612
54 0.0316453863299419
55 0.0308698860367999
56 0.0304136096590336
57 0.0297059479672412
58 0.0292796232999477
59 0.028630977204275
60 0.0282317251554081
};
\addplot [semithick, color2, mark=star, mark size=3, mark repeat=5, mark options={solid,fill opacity=0}]
table {%
1 1
2 0.520879061174101
3 0.400424102553347
4 0.327378672953895
5 0.281314549167855
6 0.247601336544235
7 0.223007245203967
8 0.204906304647255
9 0.19085333639291
10 0.178629000277713
11 0.169989645325924
12 0.163613107920213
13 0.155361295567489
14 0.146805004525683
15 0.138748714141159
16 0.130299830056813
17 0.122197430151209
18 0.115109048562145
19 0.10893361655527
20 0.103460554826998
21 0.0978751551378024
22 0.0926113826046246
23 0.0871007548594283
24 0.0817155079824768
25 0.075842558588103
26 0.0700436139118395
27 0.0634659593189023
28 0.0567597916669564
29 0.0483124479652257
30 0.0374831201709075
31 0.0209742081572238
32 0.020974208157221
33 0.0209742081572205
34 0.0209742081572182
35 0.0209742081572157
36 0.0209742081572148
37 0.0209742081572117
38 0.0209742081572102
39 0.02097420815721
40 0.0209742081572097
41 0.0209742081572092
42 0.0209742081572083
43 0.0209742081572064
44 0.0209742081572058
45 0.0209742081572047
46 0.0209742081572043
47 0.020974208157204
48 0.0209742081572037
49 0.0209742081572037
50 0.0209742081572033
51 0.020974208157203
52 0.0209742081572025
53 0.0209742081572018
54 0.0209742081572014
55 0.0209742081572009
56 0.0209742081572006
57 0.0209742081572001
58 0.0209742081572
59 0.0209742081572
60 0.0209742081571999
};
\end{axis}

\end{tikzpicture}

%% file: figs/regression1d_largenobs.tikz
% This file was created by matplotlib2tikz v0.6.11.
\begin{tikzpicture}

\definecolor{color0}{rgb}{0.12156862745098,0.466666666666667,0.705882352941177}
\definecolor{color1}{rgb}{1,0.498039215686275,0.0549019607843137}
\definecolor{color2}{rgb}{0.172549019607843,0.627450980392157,0.172549019607843}
\definecolor{color3}{rgb}{0.83921568627451,0.152941176470588,0.156862745098039}
\definecolor{color4}{rgb}{0.580392156862745,0.403921568627451,0.741176470588235}
\definecolor{color5}{rgb}{0.549019607843137,0.337254901960784,0.294117647058824}

\begin{axis}[
xlabel={Dimension of the reduced space},
ylabel={nat},
xmin=-38.95, xmax=839.95,
ymin=0.01, ymax=200,
ymode=log,
width=\figurewidth,
height=\figureheight,
tick align=outside,
tick pos=left,
x grid style={lightgray!92.026143790849673!black},
y grid style={lightgray!92.026143790849673!black},
legend cell align={left},
legend style={draw=none},
legend entries={{PCA-Y - Kullback-Leibler div.},{PCA-Y - Exp. Kullback-Leibler div.},{PCA-Y - Abs. error on the entropy},{MI - Kullback-Leibler div.},{MI - Exp. Kullback-Leibler div.},{MI - Abs.  error on the entropy}}
]
\addplot [semithick, color0, mark=*, mark size=3, mark options={solid,fill opacity=0}]
table {%
1 50.2271252782069
10 31.4552534561414
20 13.7458351791764
25 5.31380024673023
26 4.17570802197361
27 3.47966400551462
28 2.58708572880649
29 1.73515307143851
30 1.79743357215997
100 0.90103507928527
200 0.686797968432074
300 0.0236652340347313
400 0.0401165273642529
500 0.0360148496740641
600 0.00782658793463377
700 0.0135275941175423
800 0.00936826733067075
};
\addplot [semithick, color1, mark=square, mark size=3, mark options={solid,fill opacity=0}]
table {%
1 54.7068497084344
10 33.4444418288129
20 15.2437525467849
25 7.53955366054026
26 6.36618842310295
27 5.28994943368292
28 4.025392476092
29 2.76935790077259
30 2.72438006033146
100 0.996676348134934
200 0.524737478467543
300 0.176826882793531
400 0.121112617902752
500 0.0565834430686126
600 0.0126981246543656
700 0.0110060094787352
800 0.00632600940074397
};
\addplot [semithick, color2, mark=star, mark size=3, mark options={solid,fill opacity=0}]
table {%
1 54.7068497084437
10 33.4444418288201
20 15.2437525467872
25 7.53955366054144
26 6.36618842310408
27 5.28994943368424
28 4.02539247609334
29 2.76935790077385
30 2.72438006033298
100 0.996676348136432
200 0.52473747846885
300 0.176826882794273
400 0.12111261790335
500 0.0565834430692718
600 0.0126981246542357
700 0.0110060094783364
800 0.00632600940028283
};
\addplot [semithick, color3, mark=x, mark size=3, mark options={solid,fill opacity=0}]
table {%
1 47.8424160215484
10 28.7790733132787
20 12.0559976002055
25 5.40808785267155
26 4.65769312037517
27 3.3937494075008
28 2.82936022559439
29 0.620310215905605
30 1.7763660933584e-15
100 1.1053518903937e-20
200 1.09948089509808e-20
300 1.77636782967369e-15
400 1.10370089760916e-20
500 -3.55270263308494e-15
600 1.77636789860537e-15
700 1.10068997236553e-20
800 -3.55270267745054e-15
};
\addplot [semithick, color4, mark=diamond, mark size=3, mark options={solid,fill opacity=0}]
table {%
1 52.9790265673353
10 30.1526022263702
20 11.8078933784559
25 4.71938105174962
26 3.51457234528469
27 2.4081733487729
28 1.41350524252296
29 0.579752694224889
30 1.77639668715754e-15
100 4.12779096577935e-20
200 4.12851561381321e-20
300 1.77639812931982e-15
400 4.12963291926699e-20
500 -3.55267238103853e-15
600 1.77639813658201e-15
700 4.13004090379291e-20
800 -3.55267237753525e-15
};
\addplot [semithick, color5, mark=asterisk, mark size=3, mark options={solid,fill opacity=0}]
table {%
1 52.9790265673423
10 30.1526022263766
20 11.8078933784545
25 4.71938105174744
26 3.51457234528252
27 2.40817334877146
28 1.41350524252198
29 0.579752694224624
30 3.83693077310454e-13
100 3.41060513164848e-13
200 3.5527136788005e-13
300 3.41060513164848e-13
400 3.41060513164848e-13
500 2.98427949019242e-13
600 3.41060513164848e-13
700 3.41060513164848e-13
800 2.98427949019242e-13
};
\end{axis}

\end{tikzpicture}

%% file: nonlinear.tex
%!TEX root = main.tex
\section{Application to nonlinear problems}
\label{sec:application_to_nonlinear_problems}
\revis{
    We focused in the previous sections on the case of a linear problem where $X$,
    $A(X)$ and $E$ follow the multivariate normal distribution.
    In~\cite{Giraldi2017}, we showed that when $A$ is nonlinear and $X$ is
    drawn according to a uniform distribution, the PCA-YN is an appropriate
    dimensionality reduction method. The reason is that the probability density
    function of the posterior distribution of $X$ is
    \begin{equation*}
        \log f_X(X\mid Y=y) = -\frac{1}{2} \norm{y - A(X)}_{\ce^{-1}}^2 +
        \text{constant},
    \end{equation*}
    and that the PCA-YN aims to approximate the random variable $Y$ with respect to the
    Mahalanobis norm $\norm{\cdot}_{\ce^{-1}}$. This illustrates that the
    appropriate dimensionality reduction method depends on the statistical
    model used for the inference.

    In this section, we assess the benefits of the proposed approaches when the
    normality assumption is violated for $A(X)$, specifically using a log-normal model. In
    particular, we numerically evaluate the robustness of the approach using
    two values of the variance of the underlying Gaussian random vector, i.e.\ 
    a small variance yielding a model that could be well approximated with a
    Gaussian process, and a large variance where the Gaussian assumption
    no longer holds.  The last example also involves dimensionality
    reduction in a large-scale data setting, where the number of observations is
    drastically reduced.
}

\subsection{Inference problem}
\label{sub:inference_problem_nl}

\newcommand\sigt{\sigma_{F}}
\newcommand\ellt{\ell_{F}}

\revis{
    For $\bfs \in (-1,1)^2$, let $F$ be a centered stationary Gaussian process with covariance function given by
    \begin{equation*}
        k_{F}(\bfs,\bfs') = \sigma_{F, 1}^2 \exp\left(
        -\frac{\norm{\bfs-\bfs'}^2}{2\ellt^2}\right) + \sigma_{F, 2}^2
        \delta(\bfs-\bfs').
    \end{equation*}
    The nonlinear regression model of interest is based on the PCA of the
    random vector $(F(\bfs_i))_{i=1}^n$, where $(\bfs_i)_{i=1}^n$ is a
    uniformly distributed sample in $(-1,1)^2$. Let $C_F$ be the covariance
    matrix of the random vector, and its eigenpairs $(\lambda_i, w_i)_{i=1}^n
    \in (\bR^+ \times \bR^n)^n$ ordered such that $\lambda_1 \ge \hdots \ge
    \lambda_n$. The nonlinear model of interest is then given by
    \begin{equation*}
        Y_i = \exp \left(\sum_{j=1}^q B_{ij} X_j\right) + E(\bfs_i), \quad \text{where} \quad B_{ij} = (w_j)_i.
    \end{equation*}
    The prior distribution of $X$ is deduced from the PCA of $F$. Given that
    $F$ is a Gaussian process, $X$ is chosen to follow the multivariate
    normal distribution $\cN(0,C_X)$ such that
    $(C_X)_{ij}=\lambda_i\delta_{ij}$. The noise $E$ is a centered Gaussian
    process independent of $X$ and $F$ with covariance function
    \begin{equation*}
        k_E(\bfs,\bfs') = \sigma_{E,1}^2 \exp\left(-\frac{\norm{\bfs-\bfs'}}{\ell_E}\right) +
        \sigma_{E,2}^2\delta(\bfs-\bfs').
    \end{equation*}
    The synthetic data $y$ are generated using the nonlinear model
    \begin{equation*}
        \Yt_i = \exp\left(F(\bfs_i)\right) + E(\bfs_i), \quad \forall i \in
        \{1,\hdots,n\}.
    \end{equation*}
    We are therefore introducing a model error accounting for the truncation to $q$
    terms of the PCA-based expansion of $F$.
}

In the applications below, the number of parameters is set to $q=20$. Two
different sets of values for the standard deviation parameters
($\sigma_{F,1}$, $\sigma_{F,2}$, $\sigma_{E,1}$ and $\sigma_{E,2}$) will be tested to control the nonlinearity
of the mapping between the predictions and the observations.  The correlation
lengths are set to $\ellt = 0.2$ and $\ell_E = 0.05$. Finally, we shall use
$n=2,000$ observation points.

\subsection{Computation of the bases and error estimation}
\label{sub:error_estimation}

To compute the reduced bases, we rely on the expressions of the linear case which need the
determination of the second moments of the nonlinear model $A(X)= (A_i(X))_{i=1}^n$, where $A_i(X)
= \exp((BX)_i)$.
\revis{
    The analytical expressions of the mean $\ma = \bE(A(X))$ and the
    covariances $\ca = \bE((A(X)-\ma)(A(X)-\ma)^T)$ and $\cax
    =\bE((A(X)-\ma)(X-\mx)^T)$ are given by
    \begin{align*}
        (\ma)_i = \exp\left(\frac{1}{2} D_{ii}\right),\quad (\ca)_{ij} = \exp\left(\frac{1}{2} \left(D_{ii} +
        D_{jj}\right)\right) (\exp\left(D_{ij}\right) - 1),\\ \text{and} \quad
        (\cax)_{ij} = (\ma)_i \left(B\cx\right)_{ij},
    \end{align*}
    where $D = B\cx B^T$.
}

To assess the reduction error, the Kullback-Leibler divergence and the mutual information are not
available in closed form. Their accurate numerical estimation is challenging and would require prohibitive sampling of the posterior distributions, for instance using a Markov-Chain Monte Carlo method, and an estimation of the probability
density function with inherent source of error.
The situation is even more complicated for the expected Kullback-Leibler divergence,
requiring a repetitive sampling of the posterior distribution for the estimation of only one value
of this quantity.
Therefore, we choose to characterize the reduction error by its impact on the MAP value of the parameter.
The MAP is computed by solving
\begin{equation}
    \label{eq:mapf}
    \max_{x} \log f_X(x\mid Y=y) \Leftrightarrow \max_{x} \log f_Y(y\mid X=x) + \log f_X(x),
\end{equation}
for the full (unreduced) approach and, in the case of the reduced models,
\begin{equation}
    \label{eq:mapr}
    \max_{x} \log f_X(x\mid W=V^T y) \Leftrightarrow \max_{x} \log f_W(V^T y\mid X=x) + \log f_X(x).
\end{equation}
These optimization problems are solved with a trust-region Newton method, using automatic differentiation for
the evaluation of the gradient and the Hessian of the log density function.

\newcommand\xmap{x^{\text{MAP}}}
\newcommand\cmap{C^{\text{MAP}}}
\newcommand\xv{x^{\text{MAP}}_V}
\newcommand\cmapv{C^{\text{MAP}}_V}
We denote by $\xmap$ (resp.\ $\xv$) the MAP estimate of the full (resp.\ reduced)
model. Since $\xmap$ is a stationary point of the log density function of the posterior
distribution, the second-order Taylor expansion of $f_X(\cdot\mid Y=y)$ is given by
\begin{multline*}
    \log f_X(x\mid Y=y) \approx \log f_X(\xmap\mid Y=y) \\ + \frac{1}{2}(x - \xmap)^T \nabla^2\log
    f_X(\xmap\mid Y=y) (x-\xmap).
\end{multline*}
Approximating locally the distribution by the multivariate normal distribution $\cN(\xmap, \cmap)$, where
\begin{equation*}
    \cmap =-(\nabla^2\log f_X(\xmap\mid Y=y))^{-1},
\end{equation*}
gives the so-called Laplace approximation of the distribution~\cite{tierney1986accurate}.
Similarly, the posterior distribution of the reduced model will be approximated by the multivariate normal distribution
$\cN(\xv, \cmapv)$ where
\begin{equation*}
    \cmapv=-(\nabla^2\log f_X(\xv\mid W=V^T y))^{-1}.
\end{equation*}

In the following, we monitor the convergence of $\xv$ to $\xmap$ with the dimension of the reduced
space, as well as the convergence of the Hessian $(\cmapv)^{-1}$ to $(\cmap)^{-1}$ in Frobenius norm.
Note that it is empirically checked that the posterior distribution is unimodal by solving 200 times the
Problems~\eqref{eq:mapf} and~\eqref{eq:mapr} with random initial guesses drawn according to the
prior distribution. We denote by $\epsilon$ and $\epsilon_H$ the ($\Yt$-averaged) relative errors on the MAP and Hessian, respectively defined by
\begin{equation}
    \label{eq:error_estimators}
    \epsilon =  \frac{\bE_{\Yt}\left(\norm{\xv - \xmap}\right)}{\bE_{\Yt}\left(\norm{\xmap}\right)} \quad \text{and} \quad
    \epsilon_H = \frac{\bE_{\Yt}\left(\norm{\left(\cmapv\right)^{-1} -
    \left(\cmap\right)^{-1}}_{\text{Fro}}\right)}{\bE_{\Yt}\left(\norm{\left(\cmap\right)^{-1}}_{\text{Fro}}\right)}.
\end{equation}
The expectations appearing in the errors $\epsilon$ and $\epsilon_H$ are estimated by a crude Monte-Carlo method with a sample of size 70. This low sample size was found enough to obtain sufficiently correct error estimates, reflecting the robustness of all approaches which exhibit moderate dependences of the reduction error with the particular realization of $\Yt$.

\subsection{Weak nonlinearity}
\label{sub:weak_nonlinearity}
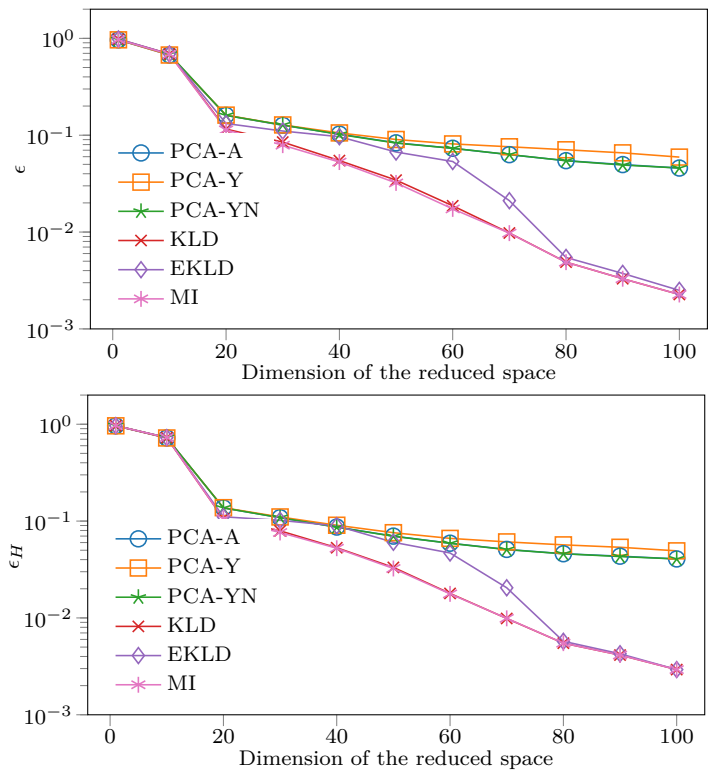
\begin{figure}[htpb]
    \centering
    \setlength\figurewidth{0.8\textwidth}
    \setlength\figureheight{0.6\figurewidth}
    \input{./figs/lognormal2d_quasilin_expectedmaperror.tikz}
    \setlength\figurewidth{0.8\textwidth}
    \setlength\figureheight{0.6\figurewidth}
    \input{./figs/lognormal2d_quasilin_expectedmaperrorh.tikz}
    \caption{Error $\epsilon$ and $\epsilon_H$ versus the dimension of the reduced space for the case of a weak
    nonlinearity.}
    \label{fig:lognormal2d_quasilinerrors}
\end{figure}

In this section, the case of a weak nonlinearity is considered, setting the standard deviations to
\begin{equation*}
    \sigma_{F,1} =0.3, \quad \sigma_{F,2} =0.001, \quad \sigma_{E,1}=0.1, \quad
    \text{and} \quad \sigma_{E,2}=0.001.
\end{equation*}
The error estimates $\epsilon$ and $\epsilon_H$ introduced in
Equation~\eqref{eq:error_estimators} are plotted in Figure~\ref{fig:lognormal2d_quasilinerrors}
against the dimension of the reduced space.
First, we observe that all the methods converge in terms of $\epsilon$ or $\epsilon_H$. All the
principal component analysis based approaches perform poorly compared to the information theoretic
techniques introduced here, with more than one order of magnitude difference when considering a
reduced space of dimension $100$. \revis{As a consequence, the normality
assumption for the computation of the reduced basis is shown to improve the
quality of the posterior distribution even when the statistical model does not
have a Gaussian structure anymore.}

It is interesting to note that the maximization of the mutual information (MI method) yields a
basis that performs slightly better than the KLD or the EKLD approaches regarding the error on the
MAP parameter $\epsilon$. The difference is less significant when considering the error on the
Hessian $\epsilon_H$ but the information theoretic methods converges faster than the principal
component analysis based approaches which tend to stagnate.

\begin{figure}[htpb]
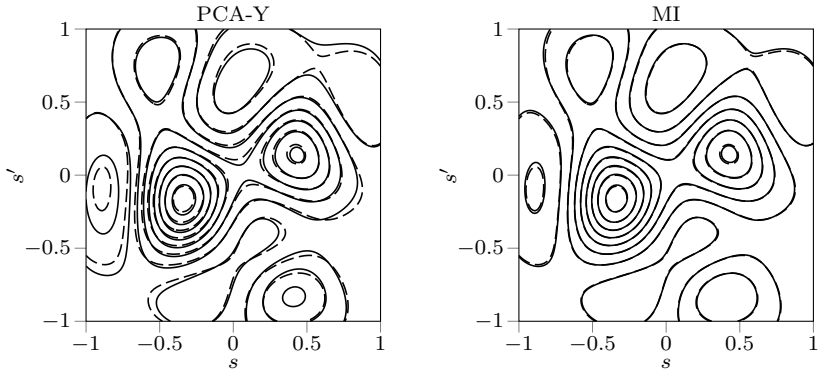

    \centering
    \setlength\figurewidth{0.45\textwidth}
    \setlength\figureheight{\figurewidth}
    \input{figs/lognormal2d_quasilin_map_kley.tikz}
    \setlength\figurewidth{0.45\textwidth}
    \setlength\figureheight{\figurewidth}
    \input{figs/lognormal2d_quasilin_map_mi.tikz}
    \caption{Contour plot of $A(\xmap)$ (dashed lines) and $A(\xv)$ (solid
        lines) for the PCA-Y (left) and the MI (right) methods with a dimension of the reduced
    space $r=60$.}
    \label{fig:lognormal2d_quasilin_map}
\end{figure}

A comparison between the MAP estimates of the field, $A(\xmap)$ and $A(\xv)$, for the same sample
$y$ of $\Yt$ and the PCA-Y and MI methods is provided in Figure~\ref{fig:lognormal2d_quasilin_map}
for the reduction with $r=60$. The plots highlight the better approximation for the MI method.

\subsection{Strong nonlinearity}
\label{sub:strong_nonlinearity}
A strong nonlinearity is considered by considerably increasing $\sigma_{F,1}$ and
$\sigma_{E,1}$ compared to Section~\ref{sub:weak_nonlinearity}. The standard deviations are
now set to
\begin{equation*}
    \sigma_{F,1} =1.5 \quad \text{and} \quad \sigma_{E,1}=0.6,
\end{equation*}
while $\sigma_{F,2}$ and $\sigma_{E,2}$ are identically set to
$0.001$. \revis{We expect now that the Gaussian assumption to be less useful
than in Section~\ref{sub:weak_nonlinearity}.}

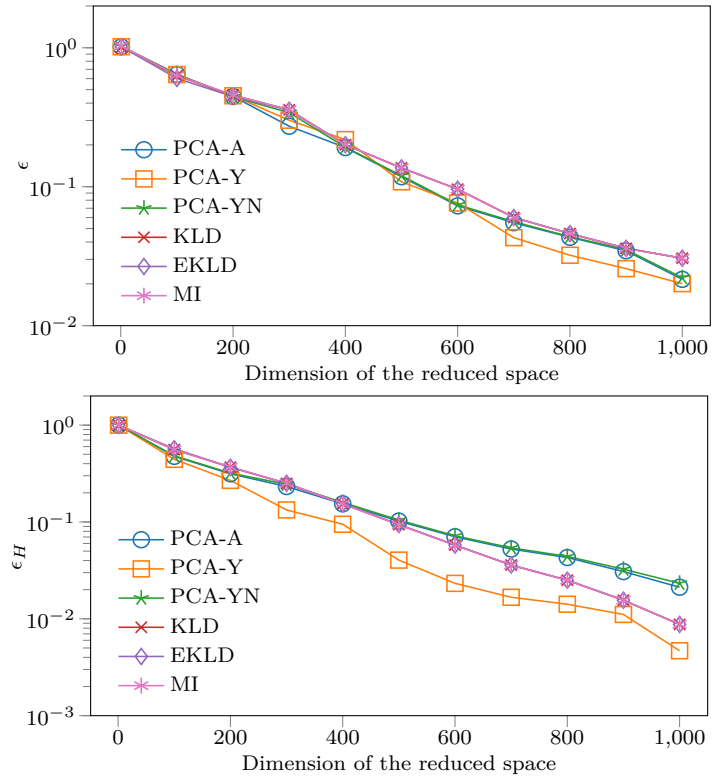
\begin{figure}[htpb]
    \centering
    \setlength\figurewidth{0.8\textwidth}
    \setlength\figureheight{0.6\figurewidth}
    \input{./figs/lognormal2d_nonlin_expectedmaperror.tikz}
    \setlength\figurewidth{0.8\textwidth}
    \setlength\figureheight{0.6\figurewidth}
    \input{./figs/lognormal2d_nonlin_expectedmaperrorh.tikz}
    \caption{Error estimates versus the dimension of the reduced space for the case of a strong
    nonlinearity.}
    \label{fig:lognormal2d_nonlinerrors}
\end{figure}

Figure~\ref{fig:lognormal2d_nonlinerrors} depicts the convergence of the error estimators
$\epsilon$ and $\epsilon_H$ with respect to the dimension of the reduced space for the different methods.
In contrast to Section~\ref{sub:weak_nonlinearity}, all the approaches exhibit
a similar convergence in terms of the error criteria $\epsilon$ and
$\epsilon_H$. Note that the PCA-Y method performs slightly better, especially
for the error on the Hessian matrix. One major difference with the previous
convergence curves reported previously in
Figure~\ref{fig:lognormal2d_quasilinerrors} is the larger dimension of the
reduced space needed to achieve a given relative error.  Indeed, the dimension
of the reduced space varies from 1 to 100 in
Figure~\ref{fig:lognormal2d_quasilinerrors} and from 1 to 1000 in
Figure~\ref{fig:lognormal2d_nonlinerrors}.  It indicates that a larger amount
of observations is required to identify the posterior distribution of the model
parameters, with similar relative accuracy, because of the non-linearities.
\revis{Even if the normality assumption is violated, the information theoretic
approaches are shown to be robust and converge to the original posterior
distribution at the same rate as the PCA based methods.}

\begin{figure}[htpb]
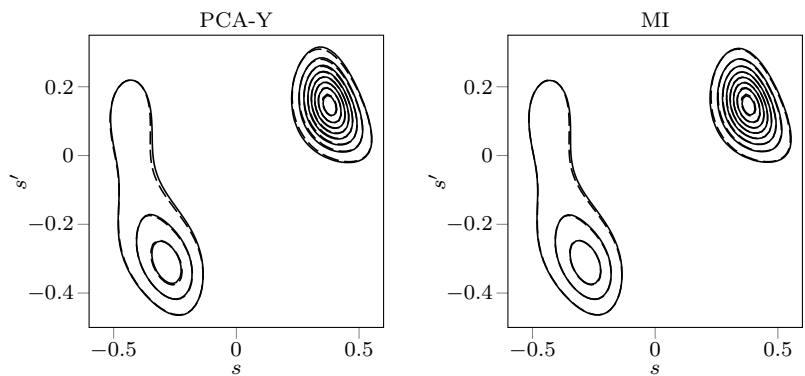

    \centering
    \setlength\figurewidth{0.45\textwidth}
    \setlength\figureheight{\figurewidth}
    \input{figs/lognormal2d_nonlin_map_kley.tikz}
    \setlength\figurewidth{0.45\textwidth}
    \setlength\figureheight{\figurewidth}
    \input{figs/lognormal2d_nonlin_map_mi.tikz}
    \caption{Contour plot of $A(\xmap)$ (dashed lines) and $A(\xv)$ (solid
        lines) for the
    PCA-Y (left) and the MI (right) methods with a dimension of the reduced space $r=400$.}
    \label{fig:lognormal2d_nonlin_map}
\end{figure}

The estimates $A(\xmap)$ and $A(\xv)$ of the field are compared in
Figure~\ref{fig:lognormal2d_nonlin_map}  for the PCA-Y and MI methods and dimension $r=400$ and the
same sample of $\Yt$.
It confirms that for this highly non-linear case and this dimension of the reduced space, the two reduction approaches yield similar accuracy.

\input{largepb}

\revis{
    \subsection{Summary}
    \label{sub:Summary}

    The numerical experiments of Sections~\ref{sub:weak_nonlinearity}--\ref{sub:largescalepb}
    suggest that the information theoretic approaches yield robust reductions
    even though they were developed for linear Gaussian models. We have shown
    in particular that they perform better then the
    PCA-based approaches, except in the strongly
    nonlinear case where all approaches behave similarly.  Moreover, the
    solution to the maximization of the mutual information is significantly
    simpler to compute than in the KLD and EKLD techniques.  Indeed, it only
    requires the solution of an eigenvalue problem and has therefore a
    computational complexity similar to the computation of the principal
    component analysis.

    Moreover, the  proposed approaches are robust to model errors as
    illustrated in Section~\ref{sub:largescalepb}. Indeed, even if we truncate
    the PCA-based expansion of the Gaussian process, $F$, the information
    theoretic reduction methods provide the lowest errors on the posterior
    distribution approximation.
}

%% file: figs/lognormal2d_quasilin_expectedmaperror.tikz
% This file was created by matplotlib2tikz v0.6.13.
\begin{tikzpicture}

\definecolor{color0}{rgb}{0.12156862745098,0.466666666666667,0.705882352941177}
\definecolor{color1}{rgb}{1,0.498039215686275,0.0549019607843137}
\definecolor{color2}{rgb}{0.172549019607843,0.627450980392157,0.172549019607843}
\definecolor{color3}{rgb}{0.83921568627451,0.152941176470588,0.156862745098039}
\definecolor{color4}{rgb}{0.580392156862745,0.403921568627451,0.741176470588235}
\definecolor{color5}{rgb}{0.890196078431372,0.466666666666667,0.76078431372549}

\begin{axis}[
xlabel={Dimension of the reduced space},
ylabel={$\epsilon$},
xmin=-3.95, xmax=104.95,
ymin=0.001, ymax=2,
ymode=log,
width=\figurewidth,
height=\figureheight,
tick align=outside,
tick pos=left,
x grid style={lightgray!92.026143790849673!black},
y grid style={lightgray!92.026143790849673!black},
legend style={at={(0.03,0.03)}, anchor=south west, draw=none},
legend cell align={left},
legend entries={{PCA-A},{PCA-Y},{PCA-YN},{KLD},{EKLD},{MI}}
]
\addplot [semithick, color0, mark=*, mark size=3, mark repeat=1, mark options={solid,fill opacity=0}]
table {%
1 0.965378744068137
10 0.672930241227853
20 0.160219287058233
30 0.126569562663705
40 0.102700234814437
50 0.0834052117768019
60 0.0732672501110935
70 0.0629489389620361
80 0.0545924362900625
90 0.049728534168395
100 0.0457991023405877
};
\addplot [semithick, color1, mark=square, mark size=3, mark repeat=1, mark options={solid,fill opacity=0}]
table {%
1 0.965406038174346
10 0.672932243466218
20 0.161204973679078
30 0.12763945439071
40 0.105568121798133
50 0.0903109170106833
60 0.081293136327848
70 0.0759644570307023
80 0.0707371932534001
90 0.0657219335677222
100 0.059426286998877
};
\addplot [semithick, color2, mark=star, mark size=3, mark repeat=1, mark options={solid,fill opacity=0}]
table {%
1 0.966732516807711
10 0.687259252786253
20 0.160015109785906
30 0.127620065855699
40 0.101852676498126
50 0.0829316440254243
60 0.0731865848168421
70 0.0631177251806346
80 0.0544887837715792
90 0.0491859678533473
100 0.0458124722645654
};
\addplot [semithick, color3, mark=x, mark size=3, mark repeat=1, mark options={solid,fill opacity=0}]
table {%
1 0.979537206364057
10 0.688687436045188
20 0.115591154110539
30 0.0849051283990551
40 0.0547191140492108
50 0.0343741188738183
60 0.0186761388506722
70 0.00977979837815434
80 0.00489071215274291
90 0.0032972454626894
100 0.00226631931783707
};
\addplot [semithick, color4, mark=diamond, mark size=3, mark repeat=1, mark options={solid,fill opacity=0}]
table {%
1 0.984224793710101
10 0.689995054803125
20 0.131992071543711
30 0.110426160462106
40 0.0964747512531145
50 0.0670440761964623
60 0.0535234359539684
70 0.021089686616729
80 0.00546736297881978
90 0.0037417255056988
100 0.00249953904493307
};
\addplot [semithick, color5, mark=asterisk, mark size=3, mark repeat=1, mark options={solid,fill opacity=0}]
table {%
1 0.967573324113702
10 0.684817226891136
20 0.113463961762532
30 0.0789109281104094
40 0.0528001643387715
50 0.0321557749634954
60 0.0173451596402706
70 0.00974349877591166
80 0.00488691011344427
90 0.00329764756484343
100 0.00226593459787866
};
\end{axis}

\end{tikzpicture}

%% file: figs/lognormal2d_quasilin_expectedmaperrorh.tikz
% This file was created by matplotlib2tikz v0.6.13.
\begin{tikzpicture}

\definecolor{color0}{rgb}{0.12156862745098,0.466666666666667,0.705882352941177}
\definecolor{color1}{rgb}{1,0.498039215686275,0.0549019607843137}
\definecolor{color2}{rgb}{0.172549019607843,0.627450980392157,0.172549019607843}
\definecolor{color3}{rgb}{0.83921568627451,0.152941176470588,0.156862745098039}
\definecolor{color4}{rgb}{0.580392156862745,0.403921568627451,0.741176470588235}
\definecolor{color5}{rgb}{0.890196078431372,0.466666666666667,0.76078431372549}

\begin{axis}[
xlabel={Dimension of the reduced space},
ylabel={$\epsilon_H$},
xmin=-3.95, xmax=104.95,
ymin=0.001, ymax=2,
ymode=log,
width=\figurewidth,
height=\figureheight,
tick align=outside,
tick pos=left,
x grid style={lightgray!92.026143790849673!black},
y grid style={lightgray!92.026143790849673!black},
legend style={at={(0.03,0.03)}, anchor=south west, draw=none},
legend cell align={left},
legend entries={{PCA-A},{PCA-Y},{PCA-YN},{KLD},{EKLD},{MI}}
]
\addplot [semithick, color0, mark=*, mark size=3, mark repeat=1, mark options={solid,fill opacity=0}]
table {%
1 0.963101408354071
10 0.721436979076386
20 0.136810170595438
30 0.108421255759059
40 0.0867616132596578
50 0.0699404943285899
60 0.0590645822048906
70 0.0509875160363899
80 0.0459972005877853
90 0.0434218920026148
100 0.0406503091024748
};
\addplot [semithick, color1, mark=square, mark size=3, mark repeat=1, mark options={solid,fill opacity=0}]
table {%
1 0.963063266081227
10 0.721395654190145
20 0.137515032820856
30 0.110338130600718
40 0.0905745694513047
50 0.0758770314909774
60 0.0662825647087259
70 0.0609683179112815
80 0.0568103278737561
90 0.0536548782243697
100 0.0492715893854724
};
\addplot [semithick, color2, mark=star, mark size=3, mark repeat=1, mark options={solid,fill opacity=0}]
table {%
1 0.966336725122194
10 0.728945555095205
20 0.136719341253987
30 0.108579311852623
40 0.086395992353525
50 0.069688564605856
60 0.0589549231505093
70 0.0507274750779708
80 0.0459739901416155
90 0.0430290771328231
100 0.0406359240375786
};
\addplot [semithick, color3, mark=x, mark size=3, mark repeat=1, mark options={solid,fill opacity=0}]
table {%
1 0.966109421566324
10 0.726108351369665
20 0.105825260537214
30 0.0792785618740123
40 0.0530490928949565
50 0.0332440296005863
60 0.0179044594121369
70 0.00987226866040916
80 0.00548682060299225
90 0.00413730715600406
100 0.0029284624447309
};
\addplot [semithick, color4, mark=diamond, mark size=3, mark repeat=1, mark options={solid,fill opacity=0}]
table {%
1 0.966377043764908
10 0.726904503891031
20 0.110720467623172
30 0.101761369430975
40 0.0892946413203527
50 0.0604233567704411
60 0.0468372975851278
70 0.0204840356070627
80 0.00573819170766851
90 0.00426570049712635
100 0.00292976890005688
};
\addplot [semithick, color5, mark=asterisk, mark size=3, mark repeat=1, mark options={solid,fill opacity=0}]
table {%
1 0.966600640470464
10 0.727705873058879
20 0.105755047851189
30 0.0762615128027939
40 0.0524305254209451
50 0.0323946750719034
60 0.0177210175628733
70 0.00990404221597358
80 0.0054879747415325
90 0.00413848328106278
100 0.00292858189562308
};
\end{axis}

\end{tikzpicture}

%% file: figs/lognormal2d_nonlin_expectedmaperror.tikz
% This file was created by matplotlib2tikz v0.6.13.
\begin{tikzpicture}

\definecolor{color0}{rgb}{0.12156862745098,0.466666666666667,0.705882352941177}
\definecolor{color1}{rgb}{1,0.498039215686275,0.0549019607843137}
\definecolor{color2}{rgb}{0.172549019607843,0.627450980392157,0.172549019607843}
\definecolor{color3}{rgb}{0.83921568627451,0.152941176470588,0.156862745098039}
\definecolor{color4}{rgb}{0.580392156862745,0.403921568627451,0.741176470588235}
\definecolor{color5}{rgb}{0.890196078431372,0.466666666666667,0.76078431372549}

\begin{axis}[
xlabel={Dimension of the reduced space},
ylabel={$\epsilon$},
xmin=-48.95, xmax=1049.95,
ymin=0.01, ymax=2,
ymode=log,
width=\figurewidth,
height=\figureheight,
tick align=outside,
tick pos=left,
x grid style={lightgray!92.026143790849673!black},
y grid style={lightgray!92.026143790849673!black},
legend style={at={(0.03,0.03)}, anchor=south west, draw=none},
legend cell align={left},
legend entries={{PCA-A},{PCA-Y},{PCA-YN},{KLD},{EKLD},{MI}}
]
\addplot [semithick, color0, mark=*, mark size=3, mark repeat=1, mark options={solid,fill opacity=0}]
table {%
1 1.01720307579403
100 0.643255246424734
200 0.446929362339516
300 0.272709668606505
400 0.191413822684989
500 0.118020100014628
600 0.0728635987078805
700 0.0550941881397691
800 0.0432669682246356
900 0.0346052797619707
1000 0.021551183340943
};
\addplot [semithick, color1, mark=square, mark size=3, mark repeat=1, mark options={solid,fill opacity=0}]
table {%
1 1.0171674614649
100 0.641705133927769
200 0.451728324887671
300 0.301042433688715
400 0.21846575751022
500 0.108499232684071
600 0.0763892311111142
700 0.0429977328335927
800 0.0321404969416221
900 0.0257362548027022
1000 0.0200599991241258
};
\addplot [semithick, color2, mark=star, mark size=3, mark repeat=1, mark options={solid,fill opacity=0}]
table {%
1 1.02715965464368
100 0.645720946396434
200 0.443057536404114
300 0.341544855290233
400 0.192647316054857
500 0.120160644968234
600 0.0741166624947651
700 0.0560125188354486
800 0.0436914254913548
900 0.0352335867474049
1000 0.0221384637922886
};
\addplot [semithick, color3, mark=x, mark size=3, mark repeat=1, mark options={solid,fill opacity=0}]
table {%
1 1.01013323613026
100 0.628841677957096
200 0.456309560288238
300 0.357818968662822
400 0.200739298364591
500 0.136654104722464
600 0.0960407405060493
700 0.0600876304346434
800 0.0460040366759968
900 0.0360372157213113
1000 0.030586335550864
};
\addplot [semithick, color4, mark=diamond, mark size=3, mark repeat=1, mark options={solid,fill opacity=0}]
table {%
1 1.01891175867237
100 0.601674737232742
200 0.445902067699038
300 0.355573272303224
400 0.200785307639121
500 0.136672823131375
600 0.0960418355911994
700 0.0600878029023608
800 0.0460042726176651
900 0.0360371983532674
1000 0.0305863427523352
};
\addplot [semithick, color5, mark=asterisk, mark size=3, mark repeat=1, mark options={solid,fill opacity=0}]
table {%
1 1.03013982227307
100 0.632056936194767
200 0.456278188850948
300 0.357820551182829
400 0.200739006483786
500 0.136654431336287
600 0.0960407665301159
700 0.0600876007440083
800 0.0460040342907534
900 0.036037233718441
1000 0.0305862752014068
};
\end{axis}

\end{tikzpicture}

%% file: figs/lognormal2d_nonlin_expectedmaperrorh.tikz
% This file was created by matplotlib2tikz v0.6.13.
\begin{tikzpicture}

\definecolor{color0}{rgb}{0.12156862745098,0.466666666666667,0.705882352941177}
\definecolor{color1}{rgb}{1,0.498039215686275,0.0549019607843137}
\definecolor{color2}{rgb}{0.172549019607843,0.627450980392157,0.172549019607843}
\definecolor{color3}{rgb}{0.83921568627451,0.152941176470588,0.156862745098039}
\definecolor{color4}{rgb}{0.580392156862745,0.403921568627451,0.741176470588235}
\definecolor{color5}{rgb}{0.890196078431372,0.466666666666667,0.76078431372549}

\begin{axis}[
xlabel={Dimension of the reduced space},
ylabel={$\epsilon_H$},
xmin=-48.95, xmax=1049.95,
ymin=0.001, ymax=2,
ymode=log,
width=\figurewidth,
height=\figureheight,
tick align=outside,
tick pos=left,
x grid style={lightgray!92.026143790849673!black},
y grid style={lightgray!92.026143790849673!black},
legend style={at={(0.03,0.03)}, anchor=south west, draw=none},
legend cell align={left},
legend entries={{PCA-A},{PCA-Y},{PCA-YN},{KLD},{EKLD},{MI}}
]
\addplot [semithick, color0, mark=*, mark size=3, mark repeat=1, mark options={solid,fill opacity=0}]
table {%
1 0.998860077483421
100 0.476017253431825
200 0.314404834086744
300 0.232293153516048
400 0.153991621797876
500 0.102010798734133
600 0.070442308996698
700 0.0525634564557087
800 0.0429374767880818
900 0.0307484687323474
1000 0.0211963482469341
};
\addplot [semithick, color1, mark=square, mark size=3, mark repeat=1, mark options={solid,fill opacity=0}]
table {%
1 0.99886288382948
100 0.44336172894735
200 0.268642219665078
300 0.133072428258241
400 0.0947870646777252
500 0.0403776650366137
600 0.0232397199216185
700 0.0167030460889129
800 0.0141482255854238
900 0.0110926348459291
1000 0.00468099785088912
};
\addplot [semithick, color2, mark=star, mark size=3, mark repeat=1, mark options={solid,fill opacity=0}]
table {%
1 0.998044804594478
100 0.482655417651996
200 0.31965610332967
300 0.244888578787653
400 0.159327614583039
500 0.104579957296946
600 0.0720463352228339
700 0.0540948639180947
800 0.0440913230859923
900 0.0325221614316721
1000 0.0232712378150932
};
\addplot [semithick, color3, mark=x, mark size=3, mark repeat=1, mark options={solid,fill opacity=0}]
table {%
1 0.996632440609248
100 0.558908184485655
200 0.368041043815553
300 0.250127817918817
400 0.15397924287709
500 0.0938774202981754
600 0.05766918520911
700 0.0359288718449547
800 0.0250384905220332
900 0.0155744254959658
1000 0.00875866529794042
};
\addplot [semithick, color4, mark=diamond, mark size=3, mark repeat=1, mark options={solid,fill opacity=0}]
table {%
1 0.995923289427044
100 0.570111722846242
200 0.367577891155937
300 0.250995934364184
400 0.154229814357358
500 0.0939264418434855
600 0.0576734403927977
700 0.0359295746378309
800 0.0250385505375545
900 0.0155744558293722
1000 0.00875870419861325
};
\addplot [semithick, color5, mark=asterisk, mark size=3, mark repeat=1, mark options={solid,fill opacity=0}]
table {%
1 0.997737738987625
100 0.558593725031164
200 0.36804139003036
300 0.250126119667017
400 0.153979104163643
500 0.0938770590075299
600 0.0576692959193144
700 0.0359288721370981
800 0.0250384951019288
900 0.0155744260592416
1000 0.0087586663817238
};
\end{axis}

\end{tikzpicture}

%% file: largepb.tex
%!TEX root = main.tex
\revis{

\subsection{Large-scale problem}
\label{sub:largescalepb}

The objective of this section is to demonstrate the feasibility, robustness and
efficiency of the proposed information-based reduction method in the context of large-scale simulations and large-dimensional observations.
To this end, we consider the problem of identifying three values $\kappa_{\Omega_{1,2,3}}$ associated with the three subdomains, $\Omega_{1,2,3}$, of the two-dimensional domain $\Omega$ depicted in the left plot of Figure~\ref{fig:add1}. These $\kappa_{\Omega_j}$ are independent and follow a log-normal distribution with parameters $\mu_\kappa$, $\sigma_\kappa$. They are
therefore expressed as
$$
    \kappa_{\Omega_j} = \exp\left[ \mu_\kappa + \sigma_\kappa X_j \right], \quad X_j \sim \cN(0,1).
$$
Thus, the vector of parameters to be inferred is $X\in\mathbb R^q$, $q=3$.
For simplicity, but without loss of generality, we shall use hereafter $\mu_\kappa = 0$ and $\sigma_\kappa=1$. The inference uses a large set of $n\approx 32,000$ observations $Y_i$ modeled as
\begin{equation} \label{eq:modeladd}
    Y_i = A_i(X) + E_i,
\end{equation}
where $A_i(X):= U(x_i)$ is the solution at the observation point $x_i \in \Omega$ of the elliptic partial differential equation with uncertain parameters $\kappa_{\Omega_j}$:
\begin{equation*}
    {\mathbb \nabla} \cdot \left( \kappa(x) {\mathbb \nabla} U(x) \right) = -1, \quad \kappa(x\in \Omega_j) = \kappa_{\Omega_j}.
\end{equation*}
The model equation is equipped with homogeneous Dirichlet (resp.\ Neumann) boundary conditions on the vertical and horizontal (resp.\ oblique) boundaries of $\Omega$.
The model for the $E_i$ is again the independent centered Gaussian model with variance $\sigma_\epsilon^2$.

\begin{figure}[htpb]
    \centering
    \includegraphics[width=0.24\textwidth]{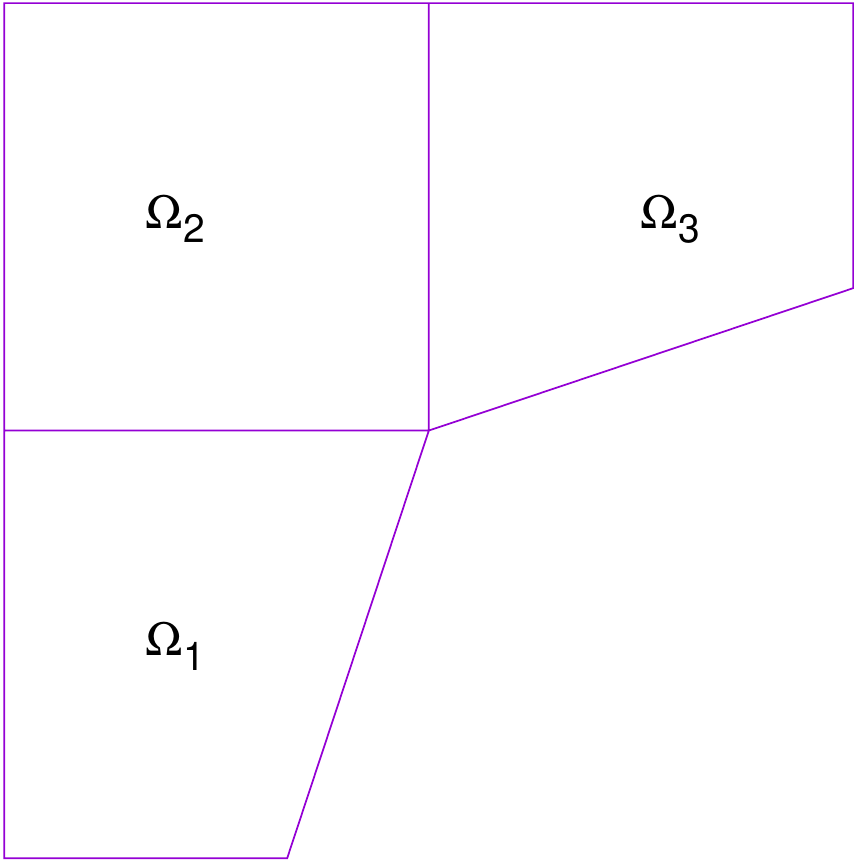}
    \includegraphics[width=0.24\textwidth]{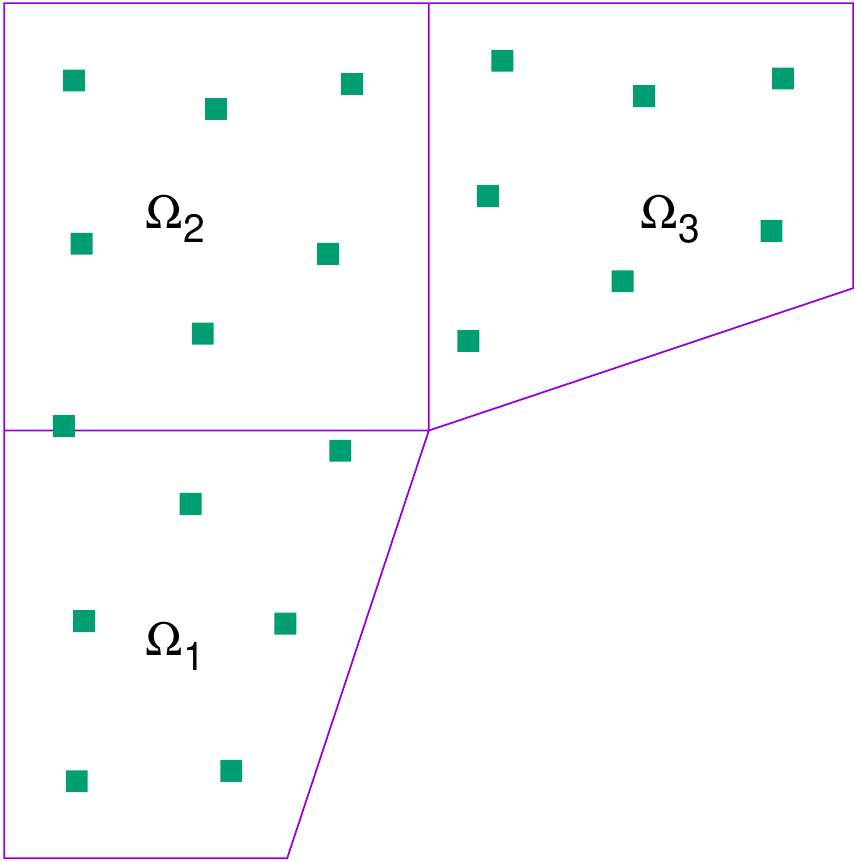}
    \includegraphics[width=0.24\textwidth]{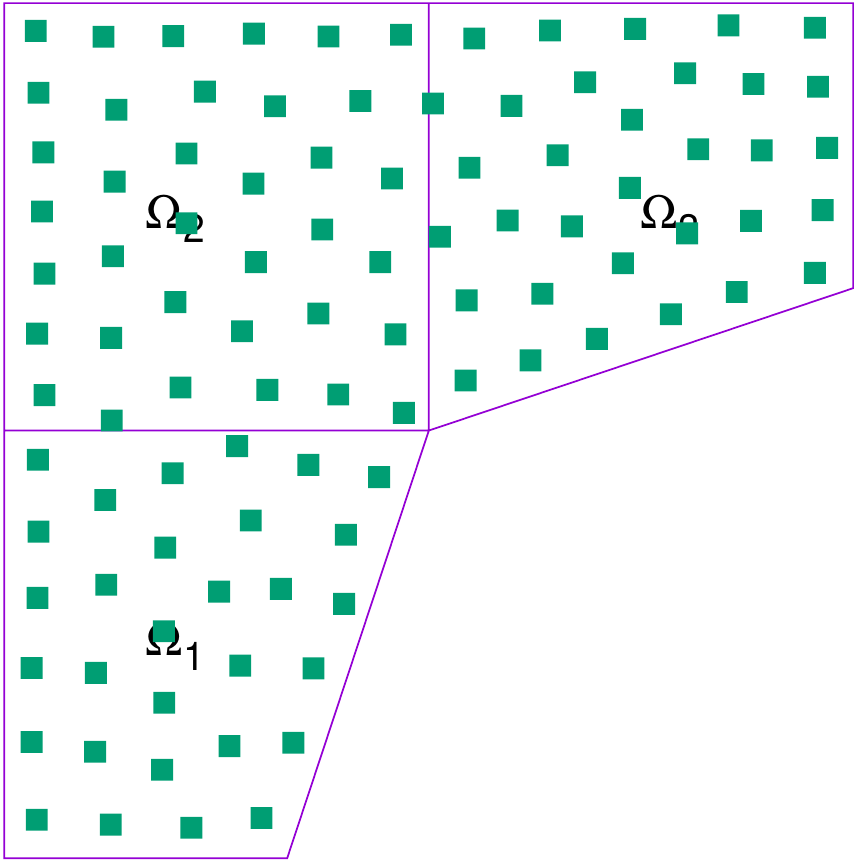}
    \includegraphics[width=0.24\textwidth]{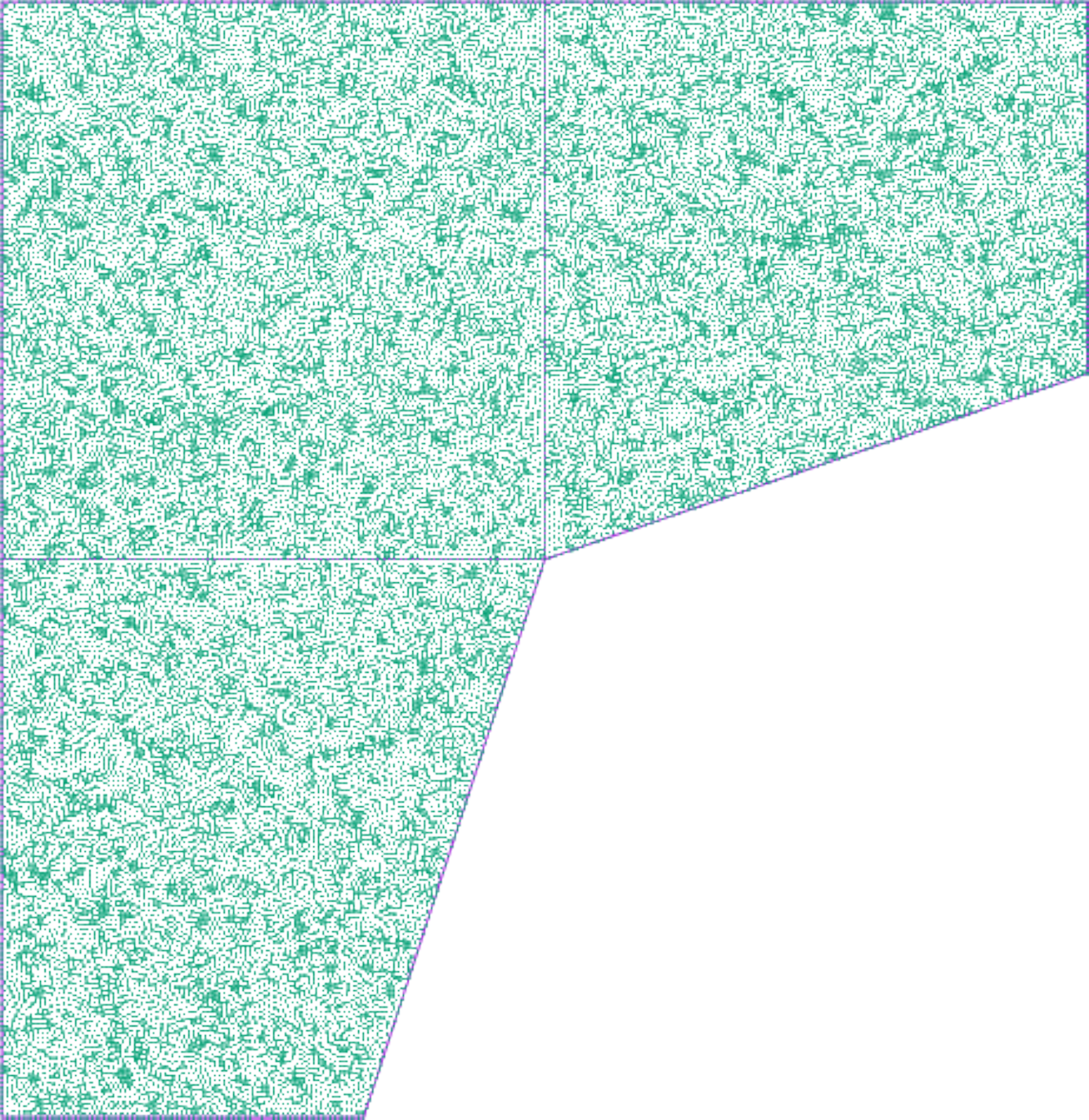}
    \caption{Left plot: Schematic of the problem domain (contained in a $3\times 3$ square) and its three subdomains, $\Omega_{j}$, over which $\kappa = \kappa_j$ is constant. Centre plots: centroids location for 20 and 100 clusters. Right plot: observation points $x_i$.}\label{fig:add1}
\end{figure}

For the reduction, we consider the maximization of the mutual information (MI), requiring the solution of~\eqref{eq:maxmi_eig_ca}.
Since $C_E$ is diagonal, the reduced basis is given by the dominant eigenspace of $C_A$.
Different approaches can be used to estimate $C_A$. Here, we rely on a Polynomial Chaos (PC) method~\cite{LeMaitre2010}, exploiting the low dimensionality of $X$, and a standard,
second-order finite element method for the spatial discretization of the elliptic problem on a very fine mesh.
As expected from the low dimensionality of $X$, the decay of the spectrum of $C_A$ is very fast.
In Figure~\ref{fig:add2} we plot the first five dominant modes of $C_A$ using the observation points shown in the right plot of Figure~\ref{fig:add1}.
Note that these observation points cover well the entire domain $\Omega$.

\begin{figure}[htpb]
    \centering
    \includegraphics[angle=-0,width=0.18\textwidth]{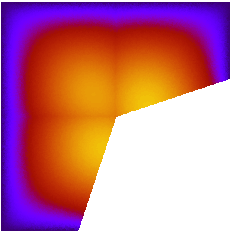}
    \includegraphics[angle=-0,width=0.18\textwidth]{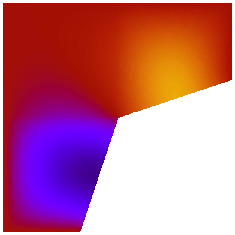}
    \includegraphics[angle=-0,width=0.18\textwidth]{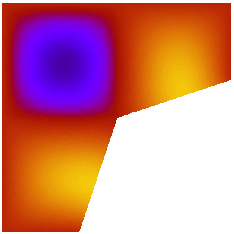}
    \includegraphics[angle=-0,width=0.18\textwidth]{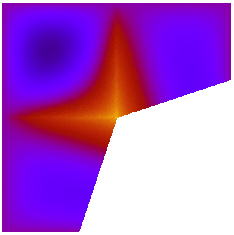}
    \includegraphics[angle=-0,width=0.18\textwidth]{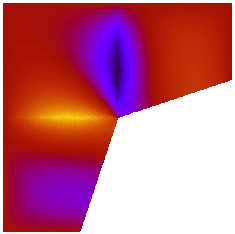}
    \caption{The five leading reduced modes (from left to right) of the MI method plotted against the $n\approx 32,000$ observation points shown in the right plot of Figure~\ref{fig:add1}.}\label{fig:add2}
\end{figure}

For comparison purposes, we also consider more reduction approaches based on observation clustering. Indeed, the amount of observations ($n\approx 32,000$) appears an overkill to learn just $q=3$ parameters. It is consequently tempting to disregard some observations and retain only $k>0$ of them to carry the inference. However, we want to maintain a sufficient coverage of the domain, and so we rely on a clustering method (k-means~\cite{Hartigan1979,MacQueen1967}) to partition the observations set into $k>0$ distinct subsets, minimizing the Euclidean distances between the $x_i$ and their respective cluster's centroids.
The  k-means procedure randomly generates clusters with a roughly equal number of observations.
In each cluster, the position $x_i$ of the selected observation is the one closest to the corresponding cluster centroid. Two examples of selected observation points are depicted in the two center plots of Figure~\ref{fig:add1}, for $k=20$ and 100 clusters respectively.
We shall refer to this reduction approach as ``Centroids.''
Disregarding all observations but the $k$-th closest to the centroids is clearly a brutal reduction approach, which is more susceptible to be affected by the noise compared to an approach
involving the projection of all observations. Consequently, one may prefer to
average (with equal weight) all the observations belonging to a cluster to
define the corresponding reduced observation.  This approach is referred to Cluster Averages (CAv) in the following.

The MI, Centroids and CAv reduction approaches are compared for three noise levels. The measurements $y_i$ are randomly generated from~\eqref{eq:modeladd} and plotted in Figure~\ref{fig:add3} to appreciate the noise to signal ratio.

\begin{figure}[htpb]
    \centering
    \includegraphics[angle=-0,width=0.32\textwidth]{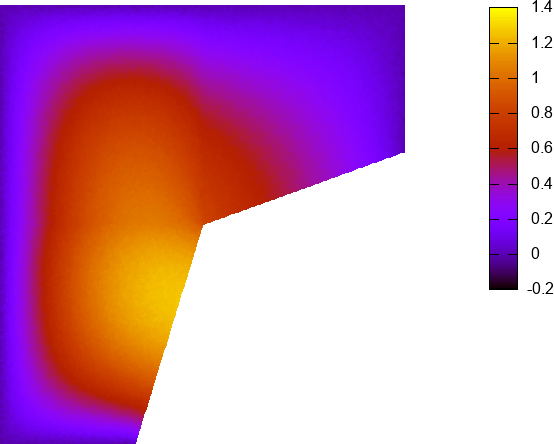}
    \includegraphics[angle=-0,width=0.32\textwidth]{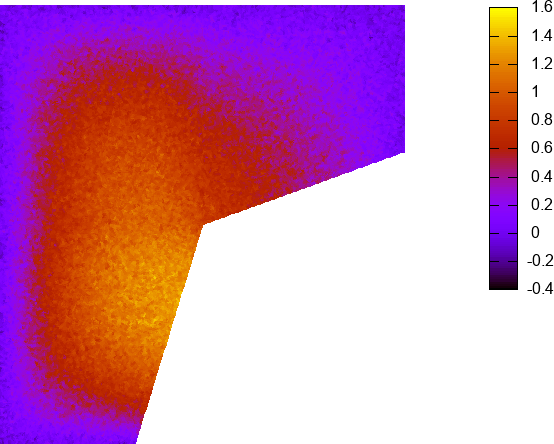}
    \includegraphics[angle=-0,width=0.32\textwidth]{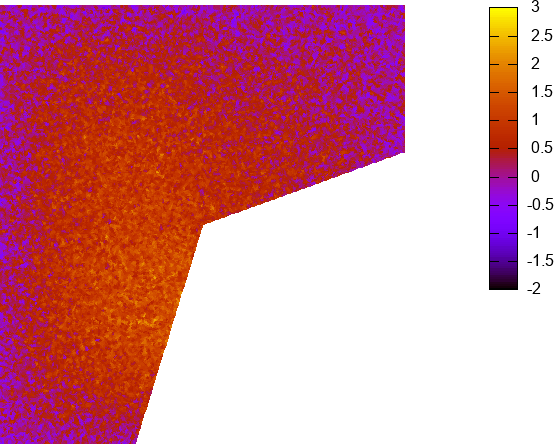}
    \caption{Measurements $y_i$ for noise level $\sigma_\epsilon=0.01$, 0.1 and 0.5 from left to right.}\label{fig:add3}
\end{figure}

To quantify the reduction errors, we consider as before the distance to the unreduced MAP point and Hessian:
\begin{equation*}
    \hat \epsilon(y) =  \frac{\norm{\xv - \xmap}}{\norm{\xmap}} \quad \text{and} \quad
    \hat \epsilon_H (y) = \frac{\norm{\left(\cmapv\right)^{-1} -
    \left(\cmap\right)^{-1}}_{\text{Fro}}}{\norm{\left(\cmap\right)^{-1}}_{\text{Fro}}}.
\end{equation*}
Note that we do not average over random observations $Y$, and restrict the
analysis to a unique measurement $y$, because of the involved computational times.
The convergence of the errors $\hat \epsilon(y)$ and $\hat \epsilon_H (y)$ with
the dimension of the reduced spaces is reported in Figure~\ref{fig:add4}, for
the three approaches and the highest noise level ($\sigma_\epsilon$). It is
seen that the MI reduction converges for roughly 10 reduced modes, and
outperforms the cluster-based reduction methods that converges at a much lower rate.
As one may have expected, the convergence of the errors in the cluster-based methods is also noisier than in MI, with Centroids exhibiting higher sensitivity to noise than CAv.

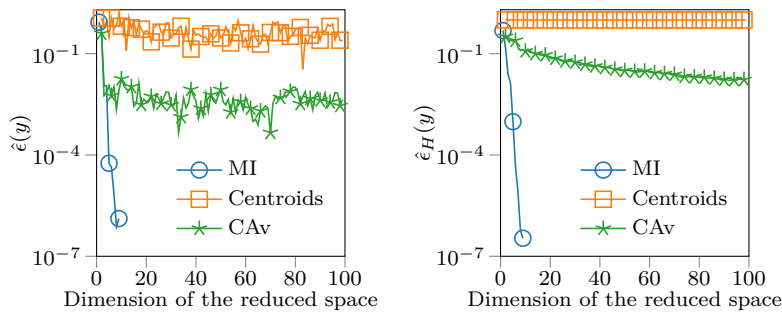
\begin{figure}[htpb]
    \centering
    \setlength\figurewidth{0.4\textwidth}
    \setlength\figureheight{\figurewidth}
    \input{./figs/fig_conv_MAP.tikz}
    \input{./figs/fig_conv_HES.tikz}
    \caption{Convergence with the reduction dimension of the MI, Centroids and Cluster Averages errors on MAP ($\hat \epsilon(y)$, left) and Hessian ($\hat \epsilon_H(y)$, right). Case of high noise level $\sigma_\epsilon=0.5$.}\label{fig:add4}
\end{figure}

However, the slow convergence of the cluster-based methods is due to the large noise in the previous example. This can be appreciated form the results reported in Figure~\ref{fig:add5}, which show that $\hat \epsilon$ and $\hat\epsilon_H$ decrease with the noise level in the CAv method, but that the convergence rate remains the same. Also note that the convergence rate of the MI method appears to be insensitive to the noise level.

% \begin{figure}[htpb]
%     \centering
%     \includegraphics[angle=-90,width=0.4\textwidth]{eps_add/fig_comp_map_MI.eps}
%     \includegraphics[angle=-90,width=0.4\textwidth]{eps_add/fig_comp_map_cav.eps} \\
%     \includegraphics[angle=-90,width=0.4\textwidth]{eps_add/fig_comp_hess_MI.eps}
%     \includegraphics[angle=-90,width=0.4\textwidth]{eps_add/fig_comp_hess_cav.eps}
%     \caption{Convergence with the reduction dimension of the MI (left) and CAv (right) MAP errors ($\hat \epsilon(y)$, top) and Hessian errors ($\hat \epsilon_H(y)$, bottom). Plotted are the errors for different noise intensities as indicated.}\label{fig:add5}
% \end{figure}
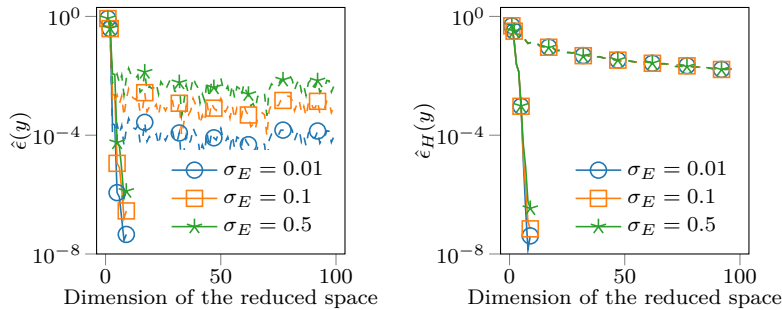
\begin{figure}[htpb]
    \centering
    \setlength\figurewidth{0.4\textwidth}
    \setlength\figureheight{\figurewidth}
    \input{./figs/fig_comp_map.tikz}
    \input{./figs/fig_comp_hess.tikz}
    \caption{Convergence with the reduction dimension of the MI (solid lines)
    and CAv (dashed lines) MAP errors ($\hat \epsilon(y)$, left) and Hessian
errors ($\hat \epsilon_H(y)$, right). Plotted are the errors for different noise intensities as indicated.}\label{fig:add5}
\end{figure}

}

%% file: figs/fig_conv_MAP.tikz
% This file was created by matplotlib2tikz v0.6.15.
\begin{tikzpicture}

\definecolor{color0}{rgb}{0.12156862745098,0.466666666666667,0.705882352941177}
\definecolor{color1}{rgb}{1,0.498039215686275,0.0549019607843137}
\definecolor{color2}{rgb}{0.172549019607843,0.627450980392157,0.172549019607843}

\begin{axis}[
xlabel={Dimension of the reduced space},
ylabel={$\hat \epsilon(y)$},
xmin=0, xmax=100,
ymin=1e-07, ymax=2,
ymode=log,
width=\figurewidth,
height=\figureheight,
tick align=outside,
tick pos=left,
x grid style={white!69.01960784313725!black},
y grid style={white!69.01960784313725!black},
legend entries={{MI},{Centroids},{CAv}},
legend cell align={left},
legend style={at={(0.97,0.03)}, anchor=south east, draw=none}
]
\addplot [semithick, color0, mark=*, mark size=3, mark repeat=4, mark options={solid,fill opacity=0}]
table {%
1 0.8542649
2 0.4175421
3 0.01148963
4 0.01116042
5 5.682254e-05
6 4.247478e-05
7 7.078861e-06
8 7.507172e-07
9 1.307389e-06
};
\addplot [semithick, color1, mark=square, mark size=3, mark repeat=4, mark options={solid,fill opacity=0}]
table {%
2 1.127348
3 0.6255686
4 0.8554357
5 0.3117713
6 1.133473
7 0.5418937
8 0.6783791
9 1.509972
10 0.6513301
11 0.9387993
12 0.122989
13 1.192423
14 0.5873419
15 0.9281147
16 0.3862511
17 0.7919518
18 0.5316503
19 0.4032765
20 0.4836151
21 0.4162074
22 0.2228457
23 0.598275
24 0.6785507
25 0.7664587
26 0.4342696
27 0.3645008
28 0.8025946
29 0.8225288
30 0.290868
31 0.3731102
32 0.3938223
33 0.5898532
34 0.6507032
35 0.2418637
36 0.5947091
37 0.2771269
38 0.1368567
39 0.4922156
40 0.352855
41 0.3842709
42 0.3224446
43 0.3081362
44 0.3133038
45 0.2257305
46 0.386555
47 0.2818606
48 0.3320542
49 0.5539906
50 0.2928341
51 0.3521978
52 0.6349424
53 0.4214505
54 0.2068712
55 0.6553629
56 0.4788875
57 0.9465154
58 0.4151391
59 0.2131707
60 0.2817981
61 0.1373662
62 0.2761942
63 0.6946787
64 0.3017283
65 0.2271422
66 0.1922754
67 0.3565725
68 0.3784862
69 0.1934878
70 0.4260194
71 0.8606372
72 0.7876662
73 0.8693384
74 0.4161687
75 0.692342
76 0.6584677
77 0.2176976
78 0.3264654
79 0.3337054
80 0.529859
81 0.5733066
82 0.5981725
83 0.0347472
84 0.2873746
85 0.1964887
86 0.3511848
87 0.310811
88 0.2940734
89 0.5489955
90 0.2503116
91 0.4894214
92 0.4361343
93 0.7368746
94 0.6460028
95 0.7637334
96 0.2583007
97 0.4177466
98 0.2505691
99 0.2360334
};
\addplot [semithick, color2, mark=star, mark size=3, mark repeat=4, mark options={solid,fill opacity=0}]
table {%
2 0.3979382
3 0.005749611
4 0.007047869
5 0.0113572
6 0.005717293
7 0.0104522
8 0.002252687
9 0.008013365
10 0.01792314
11 0.01131831
12 0.008881371
13 0.008069622
14 0.01051261
15 0.002793174
16 0.003424694
17 0.01341562
18 0.003104044
19 0.004474145
20 0.003752181
21 0.003978798
22 0.005456108
23 0.00241495
24 0.006338312
25 0.004719943
26 0.003330903
27 0.004917197
28 0.00730438
29 0.007207239
30 0.00303622
31 0.002901942
32 0.005932939
33 0.0004695041
34 0.001339258
35 0.004852263
36 0.002420801
37 0.003444069
38 0.00869523
39 0.005502883
40 0.004035447
41 0.001211407
42 0.002412985
43 0.001428613
44 0.006385772
45 0.001976188
46 0.006098633
47 0.004035376
48 0.00341915
49 0.00586523
50 0.008601027
51 0.003670779
52 0.003069632
53 0.004423278
54 0.001861657
55 0.004827204
56 0.001782002
57 0.004286486
58 0.003782953
59 0.00215729
60 0.005218037
61 0.003095254
62 0.00239417
63 0.0009794882
64 0.001158405
65 0.0007674879
66 0.00200878
67 0.002313883
68 0.002821343
69 0.002389719
70 0.0004599475
71 0.003420133
72 0.005357224
73 0.00453183
74 0.004994912
75 0.004917961
76 0.007244538
77 0.007335351
78 0.007620339
79 0.01020141
80 0.01013902
81 0.007238438
82 0.003287335
83 0.005428911
84 0.007513999
85 0.002938193
86 0.004679521
87 0.00218963
88 0.005803132
89 0.007171899
90 0.004146546
91 0.004175762
92 0.006860263
93 0.002680087
94 0.003647463
95 0.003961107
96 0.00699962
97 0.006797732
98 0.002876055
99 0.004348337
};
\end{axis}

\end{tikzpicture}

%% file: figs/fig_conv_HES.tikz
% This file was created by matplotlib2tikz v0.6.15.
\begin{tikzpicture}

\definecolor{color0}{rgb}{0.12156862745098,0.466666666666667,0.705882352941177}
\definecolor{color1}{rgb}{1,0.498039215686275,0.0549019607843137}
\definecolor{color2}{rgb}{0.172549019607843,0.627450980392157,0.172549019607843}

\begin{axis}[
xlabel={Dimension of the reduced space},
ylabel={$\hat \epsilon_H(y)$},
xmin=0, xmax=100,
ymin=1e-07, ymax=2,
ymode=log,
width=\figurewidth,
height=\figureheight,
tick align=outside,
tick pos=left,
x grid style={white!69.01960784313725!black},
y grid style={white!69.01960784313725!black},
legend entries={{MI},{Centroids},{CAv}},
legend cell align={left},
legend style={at={(0.97,0.03)}, anchor=south east, draw=none}
]
\addplot [semithick, color0, mark=*, mark size=3, mark repeat=4, mark options={solid,fill opacity=0}]
table {%
1 0.4837494
2 0.1487955
3 0.02480509
4 0.0134052
5 0.0009685271
6 4.811779e-05
7 8.503044e-06
8 5.939308e-07
9 3.425812e-07
};
\addplot [semithick, color1, mark=square, mark size=3, mark repeat=4, mark options={solid,fill opacity=0}]
table {%
2 0.9998685
3 0.9998775
4 0.999796
5 0.9997357
6 0.9998862
7 0.9997333
8 0.9995946
9 0.9998482
10 0.9997002
11 0.9992185
12 0.9995619
13 0.9998546
14 0.9996408
15 0.9994381
16 0.9995729
17 0.9995257
18 0.9993025
19 0.9992053
20 0.9991957
21 0.9994404
22 0.999072
23 0.9996519
24 0.9990268
25 0.9992117
26 0.9991881
27 0.999254
28 0.9994702
29 0.9996553
30 0.9993057
31 0.9990078
32 0.9990239
33 0.9991564
34 0.9989464
35 0.9987661
36 0.9989878
37 0.9985982
38 0.9987384
39 0.9986899
40 0.9988843
41 0.9978661
42 0.9981326
43 0.9988057
44 0.9989042
45 0.9986659
46 0.998094
47 0.9981601
48 0.9982892
49 0.9988476
50 0.9980413
51 0.9984627
52 0.9983323
53 0.9976969
54 0.9981261
55 0.9986788
56 0.9984987
57 0.9985919
58 0.9984006
59 0.9983221
60 0.9984133
61 0.9982674
62 0.9972063
63 0.9977125
64 0.9973611
65 0.9980193
66 0.9976252
67 0.9973182
68 0.9969453
69 0.9977301
70 0.9973441
71 0.996917
72 0.9979487
73 0.997747
74 0.9975645
75 0.9977116
76 0.9976485
77 0.9968163
78 0.9974698
79 0.9972412
80 0.9968131
81 0.9975339
82 0.9976234
83 0.9973047
84 0.9978744
85 0.997543
86 0.9973198
87 0.9972479
88 0.9970613
89 0.9974192
90 0.9973519
91 0.998074
92 0.9969483
93 0.9975981
94 0.9968334
95 0.9973097
96 0.9966138
97 0.9975141
98 0.9970081
99 0.9964499
};
\addplot [semithick, color2, mark=star, mark size=3, mark repeat=4, mark options={solid,fill opacity=0}]
table {%
2 0.3211986
3 0.2560193
4 0.2517845
5 0.2353417
6 0.2469215
7 0.1559307
8 0.1241276
9 0.1297902
10 0.1187083
11 0.1356853
12 0.1220024
13 0.1070486
14 0.09954392
15 0.0964768
16 0.09132336
17 0.08884447
18 0.09097102
19 0.08706542
20 0.07266862
21 0.07137748
22 0.07271926
23 0.07149274
24 0.06246675
25 0.06074834
26 0.06041555
27 0.05771263
28 0.06136251
29 0.06131479
30 0.05527048
31 0.0554833
32 0.0457112
33 0.04908116
34 0.04832171
35 0.04504976
36 0.04323705
37 0.04099829
38 0.04266054
39 0.04435286
40 0.04148377
41 0.03943054
42 0.0386126
43 0.03837514
44 0.03849186
45 0.03593487
46 0.03415178
47 0.03400593
48 0.03278301
49 0.03129481
50 0.03137487
51 0.03268631
52 0.03180702
53 0.03065241
54 0.02997762
55 0.02827581
56 0.02949137
57 0.02996472
58 0.03030607
59 0.03016689
60 0.0282615
61 0.02877613
62 0.02687294
63 0.02686469
64 0.02731901
65 0.02522593
66 0.02540027
67 0.02599254
68 0.02499079
69 0.02368981
70 0.02430125
71 0.02271913
72 0.02058564
73 0.02008995
74 0.02107333
75 0.02037846
76 0.02094812
77 0.02037336
78 0.0205063
79 0.02014498
80 0.01975614
81 0.01958985
82 0.01988595
83 0.01859134
84 0.01849045
85 0.0191709
86 0.01859602
87 0.01816228
88 0.01729418
89 0.01665487
90 0.01725303
91 0.01691767
92 0.01611763
93 0.01677304
94 0.01684996
95 0.01645364
96 0.01761489
97 0.01736516
98 0.01704926
99 0.01708467
};
\end{axis}

\end{tikzpicture}

%% file: figs/fig_comp_map.tikz
% This file was created by matplotlib2tikz v0.6.15.
\begin{tikzpicture}

\definecolor{color0}{rgb}{0.12156862745098,0.466666666666667,0.705882352941177}
\definecolor{color1}{rgb}{1,0.498039215686275,0.0549019607843137}
\definecolor{color2}{rgb}{0.172549019607843,0.627450980392157,0.172549019607843}

\begin{axis}[
xlabel={Dimension of the reduced space},
ylabel={$\hat \epsilon(y)$},
xmin=-3.9, xmax=103.9,
ymin=1e-08, ymax=2,
ymode=log,
width=\figurewidth,
height=\figureheight,
tick align=outside,
tick pos=left,
x grid style={white!69.01960784313725!black},
y grid style={white!69.01960784313725!black},
legend cell align={left},
legend style={at={(0.97,0.03)}, anchor=south east, draw=none},
legend entries={{$\sigma_E=0.01$},{$\sigma_E=0.1$},{$\sigma_E=0.5$}}
]
\addplot [semithick, color0, mark=*, mark size=3, mark repeat=4, mark options={solid,fill opacity=0}]
table {%
1 0.843511
2 0.4046246
3 0.0002327623
4 0.0002260173
5 1.163535e-06
6 8.209643e-07
7 1.542955e-07
8 2.904158e-08
9 4.555012e-08
};
\addplot [semithick, color1, mark=square, mark size=3, mark repeat=4, mark options={solid,fill opacity=0}]
table {%
1 0.8455469
2 0.4070594
3 0.002322631
4 0.002255223
5 1.129181e-05
6 8.157383e-06
7 1.477167e-06
8 1.571802e-07
9 2.81766e-07
};
\addplot [semithick, color2, mark=star, mark size=3, mark repeat=4, mark options={solid,fill opacity=0}]
table {%
1 0.8542649
2 0.4175421
3 0.01148963
4 0.01116042
5 5.682254e-05
6 4.247478e-05
7 7.078861e-06
8 7.507172e-07
9 1.307389e-06
};
\addplot [semithick, color0, dashed, mark=*, mark size=3, mark repeat=15, mark options={solid,fill opacity=0}, forget plot]
table {%
2 0.3860051
3 0.0001142147
4 0.0001394051
5 0.0002262596
6 0.0001134853
7 0.0002081874
8 4.356675e-05
9 0.0001640521
10 0.0003647024
11 0.0002271844
12 0.0001785863
13 0.0001593492
14 0.0002108529
15 5.585856e-05
16 6.826684e-05
17 0.0002709243
18 6.213534e-05
19 9.05901e-05
20 7.389612e-05
21 8.009407e-05
22 0.0001078277
23 4.81434e-05
24 0.0001278032
25 9.48522e-05
26 6.641537e-05
27 9.986764e-05
28 0.0001476512
29 0.0001455452
30 6.047377e-05
31 5.938834e-05
32 0.0001190821
33 1.015498e-05
34 2.740348e-05
35 9.760588e-05
36 4.919119e-05
37 6.846306e-05
38 0.000173284
39 0.0001099883
40 8.00651e-05
41 2.374729e-05
42 4.703368e-05
43 2.881953e-05
44 0.0001272443
45 4.044195e-05
46 0.0001228696
47 8.162897e-05
48 6.871722e-05
49 0.0001178228
50 0.0001719404
51 7.489959e-05
52 6.045131e-05
53 8.747716e-05
54 3.735589e-05
55 9.605471e-05
56 3.583355e-05
57 8.672815e-05
58 7.694627e-05
59 4.352248e-05
60 0.0001047694
61 6.19008e-05
62 4.795154e-05
63 1.956595e-05
64 2.307658e-05
65 1.553154e-05
66 4.024116e-05
67 4.64168e-05
68 5.686676e-05
69 4.844645e-05
70 9.284075e-06
71 6.926035e-05
72 0.0001076033
73 9.130993e-05
74 0.0001003438
75 9.936058e-05
76 0.0001466429
77 0.0001485296
78 0.0001546044
79 0.0002062749
80 0.0002053832
81 0.0001469246
82 6.650298e-05
83 0.0001093155
84 0.0001519235
85 5.915528e-05
86 9.402673e-05
87 4.404931e-05
88 0.0001160049
89 0.0001444909
90 8.345687e-05
91 8.420234e-05
92 0.0001384739
93 5.362188e-05
94 7.379663e-05
95 8.020217e-05
96 0.0001411819
97 0.0001367066
98 5.753328e-05
99 8.753049e-05
};
\addplot [semithick, color1, dashed, mark=square, mark size=3, mark repeat=15, mark options={solid,fill opacity=0}, forget plot]
table {%
2 0.3882585
3 0.00114354
4 0.001396931
5 0.002264017
6 0.001136166
7 0.002083453
8 0.0004380991
9 0.001634014
10 0.003635979
11 0.002270387
12 0.00178392
13 0.001597673
14 0.002107528
15 0.0005582765
16 0.0006832099
17 0.002704895
18 0.0006214106
19 0.0009038925
20 0.0007412858
21 0.0007999936
22 0.001080998
23 0.0004815925
24 0.001275968
25 0.0009476271
26 0.0006647375
27 0.0009961902
28 0.001473454
29 0.001453048
30 0.0006054989
31 0.0005912064
32 0.001190133
33 9.98274e-05
34 0.0002729582
35 0.0009751259
36 0.0004905374
37 0.0006857318
38 0.001734757
39 0.001100032
40 0.0008021293
41 0.0002386105
42 0.0004729001
43 0.0002876216
44 0.001273713
45 0.0004026668
46 0.001227196
47 0.0008144772
48 0.0006863992
49 0.001176956
50 0.001719696
51 0.0007461798
52 0.0006065154
53 0.0008768235
54 0.0003735837
55 0.0009620815
56 0.0003577011
57 0.0008649155
58 0.0007671011
59 0.0004348661
60 0.00104703
61 0.0006187312
62 0.0004796025
63 0.0001956874
64 0.000231025
65 0.0001551456
66 0.0004021107
67 0.0004638144
68 0.0005677513
69 0.0004833841
70 9.271159e-05
71 0.0006912074
72 0.00107527
73 0.0009121546
74 0.001002881
75 0.0009922097
76 0.001463361
77 0.001482052
78 0.001542133
79 0.002058898
80 0.002049326
81 0.001465139
82 0.000663509
83 0.001092166
84 0.00151637
85 0.0005908183
86 0.0009395063
87 0.0004398592
88 0.001160102
89 0.001443127
90 0.000833717
91 0.0008406089
92 0.001382575
93 0.0005363097
94 0.0007364212
95 0.0008002522
96 0.001409682
97 0.001365947
98 0.0005763151
99 0.0008748695
};
\addplot [semithick, color2, dashed, mark=star, mark size=3, mark repeat=15, mark options={solid,fill opacity=0}, forget plot]
table {%
2 0.3979382
3 0.005749611
4 0.007047869
5 0.0113572
6 0.005717293
7 0.0104522
8 0.002252687
9 0.008013365
10 0.01792314
11 0.01131831
12 0.008881371
13 0.008069622
14 0.01051261
15 0.002793174
16 0.003424694
17 0.01341562
18 0.003104044
19 0.004474145
20 0.003752181
21 0.003978798
22 0.005456108
23 0.00241495
24 0.006338312
25 0.004719943
26 0.003330903
27 0.004917197
28 0.00730438
29 0.007207239
30 0.00303622
31 0.002901942
32 0.005932939
33 0.0004695041
34 0.001339258
35 0.004852263
36 0.002420801
37 0.003444069
38 0.00869523
39 0.005502883
40 0.004035447
41 0.001211407
42 0.002412985
43 0.001428613
44 0.006385772
45 0.001976188
46 0.006098633
47 0.004035376
48 0.00341915
49 0.00586523
50 0.008601027
51 0.003670779
52 0.003069632
53 0.004423278
54 0.001861657
55 0.004827204
56 0.001782002
57 0.004286486
58 0.003782953
59 0.00215729
60 0.005218037
61 0.003095254
62 0.00239417
63 0.0009794882
64 0.001158405
65 0.0007674879
66 0.00200878
67 0.002313883
68 0.002821343
69 0.002389719
70 0.0004599475
71 0.003420133
72 0.005357224
73 0.00453183
74 0.004994912
75 0.004917961
76 0.007244538
77 0.007335351
78 0.007620339
79 0.01020141
80 0.01013902
81 0.007238438
82 0.003287335
83 0.005428911
84 0.007513999
85 0.002938193
86 0.004679521
87 0.00218963
88 0.005803132
89 0.007171899
90 0.004146546
91 0.004175762
92 0.006860263
93 0.002680087
94 0.003647463
95 0.003961107
96 0.00699962
97 0.006797732
98 0.002876055
99 0.004348337
};
\end{axis}

\end{tikzpicture}

%% file: figs/fig_comp_hess.tikz
% This file was created by matplotlib2tikz v0.6.15.
\begin{tikzpicture}

\definecolor{color0}{rgb}{0.12156862745098,0.466666666666667,0.705882352941177}
\definecolor{color1}{rgb}{1,0.498039215686275,0.0549019607843137}
\definecolor{color2}{rgb}{0.172549019607843,0.627450980392157,0.172549019607843}

\begin{axis}[
xlabel={Dimension of the reduced space},
ylabel={$\hat \epsilon_H(y)$},
xmin=-3.9, xmax=103.9,
ymin=1e-08, ymax=2,
ymode=log,
width=\figurewidth,
height=\figureheight,
tick align=outside,
tick pos=left,
x grid style={white!69.01960784313725!black},
y grid style={white!69.01960784313725!black},
legend cell align={left},
legend style={at={(0.97,0.03)}, anchor=south east, draw=none},
legend entries={{$\sigma_E=0.01$},{$\sigma_E=0.1$},{$\sigma_E=0.5$}}
]
\addplot [semithick, color0, mark=*, mark size=3, mark repeat=4, mark options={solid,fill opacity=0}]
table {%
1 0.4828048
2 0.1483177
3 0.02403351
4 0.01401087
5 0.0009348731
6 1.305883e-05
7 3.105342e-07
8 1.401932e-08
9 4.077651e-08
};
\addplot [semithick, color1, mark=square, mark size=3, mark repeat=4, mark options={solid,fill opacity=0}]
table {%
1 0.4829957
2 0.1484014
3 0.02416299
4 0.01387332
5 0.0009409314
6 1.838368e-05
7 1.588083e-06
8 1.213148e-07
9 7.018573e-08
};
\addplot [semithick, color2, mark=star, mark size=3, mark repeat=4, mark options={solid,fill opacity=0}]
table {%
1 0.4837494
2 0.1487955
3 0.02480509
4 0.0134052
5 0.0009685271
6 4.811779e-05
7 8.503044e-06
8 5.939308e-07
9 3.425812e-07
};
\addplot [semithick, color0, dashed, mark=*, mark size=3, mark repeat=15, mark options={solid,fill opacity=0}, forget plot]
table {%
2 0.3178807
3 0.2519024
4 0.2465053
5 0.2288298
6 0.2407515
7 0.1526108
8 0.1223308
9 0.1308115
10 0.1168218
11 0.1319993
12 0.1186995
13 0.1010667
14 0.09523929
15 0.09512129
16 0.09371473
17 0.09227723
18 0.09109167
19 0.08557606
20 0.07229588
21 0.0703872
22 0.07222417
23 0.07215615
24 0.06302634
25 0.06109504
26 0.06077964
27 0.05966097
28 0.0596835
29 0.05975619
30 0.05634049
31 0.05572186
32 0.04803772
33 0.04920342
34 0.04827659
35 0.04445084
36 0.0434115
37 0.04271672
38 0.04355631
39 0.04308115
40 0.04157629
41 0.04009634
42 0.03928307
43 0.03690161
44 0.03644296
45 0.03557539
46 0.03384506
47 0.03356993
48 0.03185884
49 0.03086528
50 0.02955924
51 0.03239006
52 0.03147206
53 0.03036797
54 0.02960247
55 0.02909827
56 0.02952237
57 0.02923096
58 0.02937867
59 0.02924795
60 0.02865246
61 0.02826367
62 0.02636514
63 0.02650309
64 0.02639279
65 0.0254984
66 0.02525171
67 0.02477835
68 0.02535836
69 0.0247398
70 0.02474231
71 0.0243234
72 0.02284094
73 0.02279642
74 0.02286159
75 0.02240095
76 0.02257535
77 0.02194585
78 0.02109216
79 0.02128562
80 0.02074891
81 0.019629
82 0.01991417
83 0.01890323
84 0.0189672
85 0.01864572
86 0.01824194
87 0.01800691
88 0.01772006
89 0.01744592
90 0.01740078
91 0.01751148
92 0.01641294
93 0.01638111
94 0.01634992
95 0.01672672
96 0.0169864
97 0.01688834
98 0.01653608
99 0.01627804
};
\addplot [semithick, color1, dashed, mark=square, mark size=3, mark repeat=15, mark options={solid,fill opacity=0}, forget plot]
table {%
2 0.3184928
3 0.2526613
4 0.2474788
5 0.2300315
6 0.2418904
7 0.1532177
8 0.1226593
9 0.1306207
10 0.1171781
11 0.1326746
12 0.1193065
13 0.1021766
14 0.09603537
15 0.09537092
16 0.0932705
17 0.0916439
18 0.09106756
19 0.08584962
20 0.07236134
21 0.0705647
22 0.07231209
23 0.07203045
24 0.0629142
25 0.061027
26 0.06071073
27 0.05929874
28 0.05999079
29 0.06004108
30 0.05614166
31 0.05567603
32 0.04760002
33 0.04918084
34 0.04828509
35 0.04456078
36 0.04337824
37 0.04240029
38 0.04339388
39 0.04331743
40 0.0415604
41 0.03997305
42 0.03915957
43 0.03717395
44 0.03682281
45 0.035641
46 0.03390101
47 0.0336439
48 0.03202548
49 0.03093809
50 0.02988861
51 0.03244007
52 0.0315305
53 0.03041975
54 0.0296673
55 0.02894504
56 0.02951627
57 0.02936413
58 0.02954912
59 0.02941718
60 0.02858119
61 0.02835918
62 0.02645953
63 0.02656899
64 0.02656398
65 0.02544819
66 0.02527689
67 0.02500054
68 0.0252879
69 0.02454519
70 0.02466019
71 0.02402546
72 0.02241241
73 0.02229492
74 0.02251894
75 0.0220242
76 0.02227532
77 0.02165573
78 0.02098222
79 0.02107343
80 0.02056456
81 0.01961997
82 0.01990699
83 0.01884644
84 0.01887395
85 0.01873936
86 0.01830727
87 0.01803548
88 0.01763909
89 0.01729781
90 0.01737328
91 0.01740244
92 0.01635399
93 0.01645105
94 0.01643932
95 0.01667558
96 0.01710184
97 0.01697416
98 0.01663057
99 0.0164269
};
\addplot [semithick, color2, dashed, mark=star, mark size=3, mark repeat=15, mark options={solid,fill opacity=0}, forget plot]
table {%
2 0.3211986
3 0.2560193
4 0.2517845
5 0.2353417
6 0.2469215
7 0.1559307
8 0.1241276
9 0.1297902
10 0.1187083
11 0.1356853
12 0.1220024
13 0.1070486
14 0.09954392
15 0.0964768
16 0.09132336
17 0.08884447
18 0.09097102
19 0.08706542
20 0.07266862
21 0.07137748
22 0.07271926
23 0.07149274
24 0.06246675
25 0.06074834
26 0.06041555
27 0.05771263
28 0.06136251
29 0.06131479
30 0.05527048
31 0.0554833
32 0.0457112
33 0.04908116
34 0.04832171
35 0.04504976
36 0.04323705
37 0.04099829
38 0.04266054
39 0.04435286
40 0.04148377
41 0.03943054
42 0.0386126
43 0.03837514
44 0.03849186
45 0.03593487
46 0.03415178
47 0.03400593
48 0.03278301
49 0.03129481
50 0.03137487
51 0.03268631
52 0.03180702
53 0.03065241
54 0.02997762
55 0.02827581
56 0.02949137
57 0.02996472
58 0.03030607
59 0.03016689
60 0.0282615
61 0.02877613
62 0.02687294
63 0.02686469
64 0.02731901
65 0.02522593
66 0.02540027
67 0.02599254
68 0.02499079
69 0.02368981
70 0.02430125
71 0.02271913
72 0.02058564
73 0.02008995
74 0.02107333
75 0.02037846
76 0.02094812
77 0.02037336
78 0.0205063
79 0.02014498
80 0.01975614
81 0.01958985
82 0.01988595
83 0.01859134
84 0.01849045
85 0.0191709
86 0.01859602
87 0.01816228
88 0.01729418
89 0.01665487
90 0.01725303
91 0.01691767
92 0.01611763
93 0.01677304
94 0.01684996
95 0.01645364
96 0.01761489
97 0.01736516
98 0.01704926
99 0.01708467
};
\end{axis}

\end{tikzpicture}

%% file: conclusions.tex
%!TEX root = main.tex
\revis{
\section{Conclusions and perspectives}
}
\label{sec:conclusions}

\revis{\subsection{Conclusions}}
Different optimal reductions of observations by projection in a Bayesian framework are investigated in this
work. The proposed methods are optimal in an information theoretic sense and aim at conserving the
information about the posterior distribution of interest for Gaussian linear models with correlated
additive noise.

Three optimization problems are proposed.  First, the Kullback-Leibler divergence between the
posterior distribution of the full and the reduced models is minimized.  This corresponds to an a
posteriori approach in the sense that a realization of the observations (a measurement) is required
to compute the optimal projection.  Second, we consider the minimization of the expected value of the
previous Kullback-Leibler divergence, where the expectation is taken with respect to the
observations.  As a consequence, no measurement is required to compute the optimal reduced space
and this strategy yields an a priori technique.  The last proposed approach is the maximization of
the mutual information between the projected observations and the parameters of interest. This last
approach is equivalent to the minimization of the entropy of the posterior distribution.

Solving the first two optimization problems requires specific numerical algorithms. We use in this
work the Riemannian trust-region algorithm on a manifold that take into account the invariance of
the problems. In contrast, the mutual information maximization only requires the solution to a
generalized eigenvalue problem. The computational cost and efficiency of the Riemannian algorithms 
will be addressed in a future work when large scale model will be considered.

Regarding the resulting posterior distributions, the three approaches are similar in terms of
(possibly expected) Kullback-Leibler divergence and mutual information, and perform much better, on the
considered examples, than the methods based on the principal component analysis of the
observations. We advocate therefore that the mutual information maximization is the most
appropriate approach for the determination of the optimal observation projection, given the balance
between accuracy and computational complexity.  For this particular approach, an a priori error
estimate on the mutual information loss is readily available as well as a bound on the number of
required projections. It is shown that no more projections than the rank of the linear model are
required, which is in particular lower than the number of parameters to be inferred.

\revis{
    Moreover, we addressed the linear Gaussian case in this work. However, the
    proposed approaches only require the first two moments of the distributions and
    have been successfully applied to nonlinear non-Gaussian examples, in which
optimality is no longer ensured.}

\revis{
    \subsection{Perspectives}

    In future works, the method will be applied to extreme hydrological flow problems
    (e.g.~\cite{Giraldi2017,SRAJ201482}). In particular, we plan to apply the
    approach to the framework of Ensemble Kalman filters
    (EnKF)~\cite{Evensen1994} for large datasets. The EnKF is a recursive
    Bayesian estimation technique for dynamical models of the form
    \begin{align*}
        X^{(k+1)} &= H X^{(k)} + L^{(k)},\\
        Y^{(k+1)} &= B X^{(k+1)} + E^{(k)},
    \end{align*}
    where $X^{(0)}$, $L^{(k)}$, and $E^{(k)}$ are independent Gaussian vectors. Note that
    the equation above is the same as \eqref{eq:linearmodel}.  To
    estimate the posterior distribution of $X^{(k+1)}$, the Kalman filter
    requires the inversion of the covariance matrix $C_{Y^{(k+1)}}$ at each
    iteration of the discrete dynamical system. However, in the EnKF,
    $C_{X^{(k+1)}}$ is estimated using a Monte-Carlo estimator with a sample
    size that can be much lower than the total number of observations $n$. As a
    consequence, the covariance of the forward state $X^{(k+1)}$ is low-rank
    and we showed in this paper that a low number of projections of
    the observations are enough to recover the mutual information between the
    estimated distribution of the state $X^{(k+1)}$ and the observations
    $Y^{(k+1)}$.

    Additional challenges arise when the datasets are high dimensional.
    Considering the mutual information based technique, the problem could be
    first tackled using high performance computing. Given that we a priori know
    an upper bound on the number of projections that is already low, we only
    need an efficient matrix product computation (e.g.~\cite{Drineas2006}) to
    implement the algorithm from~\cite{Absil2007}.
    %We emphasize that in the practical applications we
    %are considering, the use of a supercomputers for the simulation
    %of the physical model is mandatory, see e.g.~\cite{Theussl2016}.
    %Another approach could be based on the use of 
    %sketching techniques to decrease the size of the generalized
    %eigenvalue problem, see e.g.~\cite{Drineas2006} for a fast randomized
    %matrix product.
    Further developments are required to appropriately
    use these approaches in a
    streaming environment.
}

%% file: ack.tex
\section*{Acknowledgments}
This work is supported by King Abdullah University of Science and Technology
Awards CRG3-2156 and OSR-2016-RPP-3268.

%% file: proofs.tex
%!TEX root = main.tex
\appendix

\section{Proof of Proposition~\ref{prop:posteriordistrib}}
\label{proof:posteriordistrib}

According to Bayes' theorem, the posterior distribution is such that
\begin{equation*}
    f_X(x \mid Y=y) \propto f_Y(y \mid X=x) f_X(x),
\end{equation*}
or equivalently,
\begin{multline*}
    \log f_X(x\mid Y=y)
    =\log f_Y(y \mid X=x) + \log f_X(x) + k_0 \\
    = -\frac{1}{2} \left(\left(y - Bx - \me\right)^T\ce^{-1}\left(y - Bx -
            \me\right)\right. \\ \left. + \left(x
    - \mx\right)^T\cx^{-1}\left(x - \mx\right) \right) + k_1,
\end{multline*}
where $k_0$ and $k_1$ are constants. Because the log probability density function is quadratic with
respect to $x$, we conclude that the posterior distribution is also a multivariate normal
distribution, i.e. $P(X\mid Y=y)\sim \cN(\ms,\cs)$. This implies that, up to a constant $k_2$, the
following equality holds
\begin{equation*}
    \log f_X(x\mid Y=y) = -\frac{1}{2} \left(\left(x - \ms\right)^T \cs^{-1} \left(x -
    \ms\right)\right) + k_2.
\end{equation*}
Identifying the quadratic terms in $x$ and using the Woodbury matrix
identity~\cite[Equation~(29)]{householder1957} gives
\begin{align*}
    \cs^{-1} &= \cx^{-1} + B^T\ce^{-1}B = \cx^{-1}\left(\cx + \cax^T\ce^{-1}\cax\right)
    \cx^{-1}\\
    \text{and} \quad \cs &= \cx - \cx B^T \left(\ce + B\cx B^T\right)^{-1} B\cx = \cx - \cax^T
    \cy^{-1} \cax.
\end{align*}
Identifying the linear term w.r.t. $x$ yields
\begin{align*}
    \cs^{-1}\ms &= \cx^{-1}\mx + B^T\ce^{-1} \left(y - \me\right) \\
            & = B^T\ce^{-1}\left(y - \my\right) + B^T\ce^{-1}\ma + \cx^{-1}\mx.
\end{align*}
% With $\gs = \cs B^T \ce^{-1} = \cax^T (\rmI - \cy^{-1}\ca)\ce^{-1} = \cax^T\cy^{-1}$ and $\hs =
% \cs\cx^{-1}\mx +
% \gs \ma$, we finally have
% \begin{equation*}
%     \ms = \gs \left(y - \my\right) + \hs.
% \end{equation*}
We finally have
\begin{equation*}
    \ms = \gs \left(y - \my\right) + \hs,
\end{equation*}
with
\begin{align*}
    \gs &= \cs B^T \ce^{-1} = \cax^T (\rmI - \cy^{-1}\ca)\ce^{-1} =
    \cax^T\cy^{-1},\\
    \hs &=  \cs\cx^{-1}\mx + \gs \ma.
\end{align*}

For the posterior distribution of the reduced model, we substitute $(y-\my)$, $\cy$, $\me$,
$\ce$, $\ma$, $\cax$, and $\ca$ respectively by $V^T(y-\my)$, $V^T\cy V$, $V^T\me$, $V^T\ce V$,
$V^T\ma$, $V^T\cax$, and $V^T\ca V$ in the full model. The fact that $V$ is full-rank ensures
that $V^T\cy V$, $V^T\ce V$ and $$\cx + \cax^T V(V^T\ce V)^{-1}V^T\cax$$ are symmetric positive
definite matrices and hence are invertible.

\section{Proof of Proposition~\ref{prop:invariancemeancov}}
\label{proof:invariancemeancov}

For $M\in\glr$ we have
\begin{equation*}
    \cax VM(M^T V^T\ce VM)^{-1} M^T V^T\cax = \cax V (V^T\ce V)^{-1} V^T\cax,
\end{equation*}
so we deduce that $C_{VM} = \cv$. Moreover, given that
\begin{equation*}
    G_{VM} = \cax^T VM(M^T V^T\cy VM)^{-1} = \cax^T V(V^T\cy V)^{-1}M^{-T} = G_V
    M^{-T},
\end{equation*}
we conclude that $G_{VM} (VM)^T = G_V V^T$, $h_{VM}=h_V$ and finally $m_{VM} = \mv$.

\section{Proof of Proposition~\ref{prop:kld}}
\label{proof:kld}

Using the definition of the Kullback-Leibler divergence~\eqref{eq:defkld}, we have
\begin{multline*}
\kld{P(Z_0)}{P(Z_1)} =
\bE_Z\left(\log\left(\frac{\det(C_1)^{\frac{1}{2}}}{\det(C_0)^{\frac{1}{2}}}\right)\right.\\
    \left.-
    \frac{1}{2}\left(Z-m_0\right)^T
C_0^{-1}\left(Z-m_0\right) + \frac{1}{2}\left(Z-m_1\right)^T
C_1^{-1}\left(Z-m_1\right)\right). % \\
% &\quad=\frac{1}{2} \left(-\log\det\left(C_0C_1^{-1}\right) - \bE_Z\left(\left(Z-m_0\right)^T
% C_0^{-1}\left(Z-m_0\right) \right) + \bE_Z\left((Z-m_1)^T C_1^{-1}(Z-m_1)\right)\right).
\end{multline*}
Given that $Z\sim P(Z_0)$, we deduce that
\begin{align}
    \label{eq:tracetech}
    \bE_Z\left(\left(Z-m_0\right)^T C_0^{-1}\left(Z-m_0\right)\right) &=
    \bE_Z\left(\trace\left(\left(Z-m_0\right)^T C_0^{-1}\left(Z-m_0\right)\right)\right)\\
    &=\trace(\bE_Z((Z-m_0)(Z-m_0)^T) C_0^{-1}) \notag\\
    &=\trace(C_0 C_0^{-1}) = q. \notag
\end{align}
Moreover we have
\begin{align*}
    &\bE_Z\left(\left(Z-m_1\right)^T C_1^{-1}\left(Z-m_1\right)\right) \\
    &\quad= \bE_Z((Z-m_0 + m_0 - m_1)^T C_1^{-1} (Z - m_0 + m_0 - m_1)) \\
    &\quad= \bE_Z\left(\left(Z-m_0\right)^T C_1^{-1} \left(Z-m_0\right) + \left(m_0 - m_1\right)^T
    C_1^{-1}\left(m_0 - m_1\right) \right.\\
    &\qquad\qquad\qquad\left.+ 2 \left(Z-m_0\right)^T C_1^{-1}\left(m_0 - m_1\right)\right).
\end{align*}
Using the same trace technique as in Equation~\ref{eq:tracetech}, and using the fact that
$\bE_Z(Z)=m_0$, the term $\bE_Z((Z-m_1)^T C_1^{-1}(Z-m_1))$ is equal to
\begin{equation*}
    \bE_Z\left(\left(Z-m_1\right)^T C_1^{-1}\left(Z-m_1\right)\right) =
    \trace\left(C_0C_1^{-1}\right) + \left(m_0 - m_1\right)^T C_1^{-1}\left(m_0 - m_1\right),
\end{equation*}
yielding the final result
\begin{multline*}
    \kld{P(Z_0)}{P(Z_1)} =\\ \frac{1}{2}\left(\trace\left(C_0C_1^{-1}\right)
        -\log\det\left(C_0C_1^{-1}\right) - q + \left(m_0 - m_1\right)^T C_1^{-1}\left(m_0 -
    m_1\right)  \right).
\end{multline*}

\section{Proof of Theorem~\ref{th:existenceminkld}}
\label{proof:existenceminkld}

First, the map $\sJ_0$ is smooth ($\in\cC^{\infty}$) as the sum and composition of smooth
functions, noting that the determinant is always strictly positive.

Let $\pi:\bR^{n\times r}_* \to \Grn$ denotes the canonical projection defined by $\pi(V)=[V]$. Let
$\sK_0:\Grn\to\bR$ be the map defined by $\sJ_0(V) = \sK_0 \circ \pi(V)$.
$\sK_0$ is in fact the functional we are minimizing in Problem~\eqref{eq:min_kld}.

According to~\cite[Proposition~3.4.5]{Absil2009}, the smoothness of $\sJ_0$ implies that $\cK_0$ is
smooth and in particular continuous. According to~\cite[Lemma 5.1]{Milnor1974} $\Grn$ is compact,
the extreme value theorem concludes the proof.

\section{Proof of Proposition~\ref{prop:expkld}}
\label{proof:expkld}

Since only $\ms$ and $\mv$ depend on $Y$ in Equation~\eqref{eq:kld_as_fun_of_ldd_md}, the
expected Kullback-Leibler divergence admits the form
\begin{multline*}
    \bE_Y \left(\kld{P(X\mid Y)}{P(X\mid W=V^T Y)}\right) =\\ \frac{1}{2} \left(\ldd{\cs}{\cv} +
    \bE_Y \left(\md{\cv}{\ms}{\mv}\right)\right).
\end{multline*}
Note that $\ms = \gs (Y-\my) + \hs$ and $\mv = \gv V^T(Y-\my) + \hv$, hence
\begin{equation*}
    \ms - \mv = (\gs - \gv V^T) (Y-\my) + (\hs - \hv),
\end{equation*}
and
\begin{multline*}
    \bE_Y\left(\md{\cv}{\ms}{\mv}\right) =\\ \trace \left(\cv^{-1}\left(\gs -\gv
            V^T\right)\cy\left(\gs-\gv
    V^T\right)^T\right) + (\hs - \hv)\cv^{-1}(\hs - \hv),
\end{multline*}
which yields the final result.

\section{Proof of Theorem~\ref{th:maxmi}}
\label{proof:maxmi}

For a normally distributed $\bR^n$-valued random variable $Z\sim \cN(m_Z, C_Z)$, the entropy
$H(Z)$ is given by
\begin{equation*}
    H(Z) = \frac{1}{2}\log(\det(C_Z)) + \frac{n}{2}\log(2\pi e).
\end{equation*}
Given that $X$ and $W$ are normally distributed, we immediatly deduce
\begin{align*}
    H(X) &= \frac{1}{2} \log(\det(\cx)) + \frac{q}{2}\log(2\pi e), \\
    \text{and} \quad H(W) &= \frac{1}{2} \log(\det(V^T\cy V)) + \frac{r}{2}\log(2\pi e).
\end{align*}
In order to compute the joint-entropy $H(W,X)$, we need to characterize the covariance of $(W,X)$.
Note that we already know that $(W,X)$ is drawn according to a Gaussian distribution. In order to
obtain the covariance $C_{(W,X)}$, we identify the quadratic terms in the following
equality between the probability density functions:
\begin{equation*}
    \log f_{(W,X)}(W,X) = \log f_W(W\mid X) + \log f_X(X),
\end{equation*}
where the likelihood $f_W(W\mid X)$ is directly deduced from
Equation~\eqref{eq:reducedlinearmodel}. The conditional random distribution $P(W\mid X)$ follows
the Gaussian distribution $\cN(V^T (AX+\me), V^T\ce V)$.  Identifying the quadratic terms
yields
\begin{equation*}
    C_{(W,X)}^{-1} =
    \begin{pmatrix}
        \left(V^T\ce V\right)^{-1} & -\left(V^T\ce V\right)^{-1}V^T A \\
        - A^T V\left(V^T\ce V\right)^{-1} & \cx^{-1} + A^T V\left(V^T\ce V\right)^{-1} V^T A
    \end{pmatrix}.
\end{equation*}
According to~\cite[Section 9.1.2]{Petersen2012}, the determinant of the precision
matrix $C_{(W,X)}^{-1}$ is given by
\begin{equation*}
    \det\left(C_{(W,X)}^{-1}\right) = \det\left(\left(V^T\ce V\right)^{-1}\right) \det\left(\cx^{-1} \right).
\end{equation*}
We immediatly have
\begin{align*}
    H(W,X) &= \frac{1}{2}\log\det(C_{(W,X)}) + \frac{r + q}{2}\log(2\pi e)\\
           &= -\frac{1}{2}\log\det(C_{(W,X)}^{-1}) + \frac{r + q}{2}\log(2\pi e)\\
            &= -\frac{1}{2}\left(\log\det\left(\cx^{-1}\right) + \log\det\left(\left(V^T\ce
            V\right)^{-1}\right)\right) + \frac{r + q}{2}\log(2\pi e),
\end{align*}
and the mutual information reduces to
\begin{align*}
    \cI(W,X) &= \frac{1}{2}\left(\log\det\left(V^T\cy V\right) + \log\det\left(\left(V^T\ce
V\right)^{-1}\right)\right) \\
    &= \frac{1}{2}\log\det\left(\left(V^T\cy V\right)\left(V^T\ce V\right)^{-1}\right),
\end{align*}
which proves the first equality.

Regarding the entropy of the posterior distribution, we know that $P(X\mid W=V^T y) \sim
\cN(\mv,\cv)$, yielding
\begin{equation*}
    H(P(X\mid W=V^T y)) = \frac{1}{2} \log\det(\cv) + \frac{q}{2}\log(2\pi e),
\end{equation*}
and the entropy does not depend on the realization $y$ of $Y$. Using Equation~\eqref{eq:covreducedposterior}, we have
\begin{align*}
    &\log\det\left(\cv\right)\\
    &\quad= \log\det\left(\sqrtm{\cx}(\rmI -
    \sqrtm{\cx}B^T V\left(V^T\cy V\right)^{-1}V^T B\sqrtm{\cx})\sqrtm{\cx}\right)\\
    &\quad=\log\det\left(\rmI -
    \sqrtm{\cx}B^T V\left(V^T\cy V\right)^{-1}V^T B\sqrtm{\cx}\right) + \log\det \cx\\
    &\quad=\log\det\left(V^T\cy V - V^T B\cx B^T V\right) - \log\det \left(V^T\cy V\right) + \log\det
    \cx.
\end{align*}
The last equality is obtained using the identity $\det(\rmI + MN)=\det(\rmI + NM)$, with
$M=(V^T\cy V)^{-1/2}V^T B\cx^{1/2}$ and $N=M^T$ and factorizing the resulting expression by
$(V^T\cy V)^{-1/2}$ on the left and right in the determinant. We finally find that
\begin{align*}
    &\log\det\left(\cv\right) \\
    &\quad= \log\det\left(V^T\cy V - V^T \ca V\right) - \log\det \left(V^T\cy V\right)
    + \log\det \cx\\
    &\quad= \log\det\left(V^T\ce V\right) - \log\det \left(V^T\cy V\right) + \log\det \cx\\
    &\quad= \log\det \cx - \log\det\left(\left(V^T\cy V\right) \left(V^T\ce V\right)^{-1} \right),
\end{align*}
and the entropy is
\begin{multline*}
    H(P(X\mid W=V^T y)) =\\ -\frac{1}{2} \log\det\left(\left(V^T\cy V\right) \left(V^T\ce
    V\right)^{-1}\right) + \frac{1}{2}\log\det(\cx) + \frac{q}{2}\log(2\pi e),
\end{multline*}
that proves the second equality.

For the last part of the proof, we consider the maximization problem
\begin{equation}
    \max_{V\in\bR^{n\times r}_*} \log\det \left(\left(V^T\cy V\right)\left(V^T\ce V\right)^{-1}\right) =
    \max_{V\in\bR^{n\times r}_*} 2~\cI\left(W,X\right).
\end{equation}
First let us introduce a change of variable, setting $U = \ce^{\frac{1}{2}} V$. The
optimization problem becomes
\begin{equation*}
    \max_{U\in\bR^{n\times r}_*}\det \left(\left(U^T\ce^{-\frac{1}{2}}\cy\ce^{-\frac{1}{2}}
    U\right)\left(U^T U\right)^{-1}\right).
\end{equation*}
Then, the quantity $\cK(U) = \det ((U^T\cy U)(U^T U)^{-1})$ is invariant under
any invertible linear transformation on the right, meaning that $\cK(U) = \cK(UQ)$ for any
$Q\in\bR^{r\times r}$ invertible. With $\Stn$ denoting the Stiefel manifold defined by
\begin{equation*}
    \Stn = \{M \in \bR^{n\times r};\ M^T M = \rmI_r\},
\end{equation*}
there exists a matrix $\hat U \in \Stn$ such that $\cK(U) = \cK(\hat U)$. Such a matrix $\hat U$
can be computed using, for instance, a thin QR factorization. We can therefore consider the
following equivalent optimization problem
\begin{equation}
    \label{eq:maxdetu}
    \max_{U\in\Stn} \det
    \left(U^T\ce^{-\frac{1}{2}}\cy\ce^{-\frac{1}{2}} U\right).
\end{equation}
In order to conclude the proof, we need the following result.

\begin{lemma}
    \label{lem:maxdet}
    Let $K\in\bR^{n\times n}$ be a symmetric positive definite matrix with
    eigenvalues $(\lambda_i)_{i=1}^n$ in a decreasing order. Then we have
    \begin{equation}
        \label{eq:maxdet}
        \max_{U \in \Stn} \log\det\left(U^T K U\right) = \sum_{i=1}^r \log\lambda_i.
    \end{equation}
    Moreover, any solution to the optimization Problem~\eqref{eq:maxdet} is an invariant subspace
    of $K$ and a particular solution is given by the matrix $U$ whose columns are the eigenvectors
    of $K$ associated to the eigenvalues $(\lambda_i)_{i=1}^r$.
\end{lemma}
\begin{proof}
    First, a solution to Problem~\eqref{eq:maxdet} exists using the fact that $\sF:U\mapsto
    \log\det(U^T K U)$ is continuous and
    $\Stn$ is compact. It is closed as the inverse image of $\{0\}$ by the continuous function $U\mapsto
    U^T U-I$, and bounded because $\normt{U}_{\text{Fro}}^2 = r$ for all $U\in\Stn$. The extreme
    value theorem implies the existence of a maximizer.

    \newcommand\Us{U_\star}
    \newcommand\Ps{\Psi_\star}
    Let us introduce the map $\sH:\bR^{n\times r} \times \bR^{r\times r} \to \bR$ be defined by
    \begin{equation*}
        \sH(U, \Psi) = 2 \trace\left(\left(U^T U - I\right)\Psi\right),
    \end{equation*}
    and consider the Lagrangian function $\sL(U,\Psi) = \sF(U) + \sH(U,\Psi)$ associated to the
    constrained optimization Problem~\eqref{eq:maxdet}. An optimal solution $(\Us,\Ps)$ satisfies
    the equation
    \begin{equation*}
        \rmD_U\sL(\Us, \Ps)[\delta U] = 0, \quad \forall \delta U\in \bR^{n\times r},
    \end{equation*}
    where $\rmD_U\sL(\Us,\Ps)[\delta U]$ denotes the G\^ateaux derivative of the Lagrangian $\sL$ at $\Us$
    in the direction $\delta U$ with respect to the first parameter.
    Given the formula
    \begin{equation*}
        \frac{\rmd}{\rmd t} \log\det\left(U + t\delta U\right) = 2 \trace \left(\left(U +
        t\delta U\right)^{-1}\delta U\right),
    \end{equation*}
    we conclude that the G\^ateaux derivative $\rmD_U\sF(U)[\delta U]$ is
    \begin{equation*}
        \rmD_U\sF(U)[\delta U] = 2 \trace\left(\left(U^T K U\right)^{-1} U^T K \delta U\right),
    \end{equation*}
    and similarly we have
    \begin{equation*}
        \rmD_U\sH(U)[\delta U] = 2 \trace((\Ps + \Ps^T) U^T \delta U).
    \end{equation*}
    Hence, for all $\delta U\in \bR^{n\times r}$, a solution $\Us$ to Problem~\eqref{eq:maxdet} satisfies
    \begin{align*}
        \rmD_U\sL(\Us, \Ps)[\delta U] = 2 \trace \left(\left(\left(\Us^{T}K \Us\right)^{-1} \Us^{T}K
        + (\Ps + \Ps^T) \Us^T\right) \delta U\right) = 0.
    \end{align*}
    The result holding for all $\delta U$, we conclude that $\Us$
    satisfies
    \begin{align*}
        &\left(\Us^{T}K \Us\right)^{-1} \Us^{T}K
        + (\Ps + \Ps^T) \Us^T = 0\\
        \Leftrightarrow \qquad&K\Us = -\Us (\Ps + \Ps^T)\Us^T K\Us.
    \end{align*}
    Finally, multiplying this last equation on the left by $\Us^T$ and on the right by $(\Us^T K \Us)^{-1}$
    gives that $\Ps + \Ps^T = - \rmI_r$ and
    \begin{equation*}
        K\Us = \Us \Us^T K\Us,
    \end{equation*}
    meaning that $\Us$ spans an $r$-dimensional invariant subspace of $K$.

    \bigskip
    To conclude the proof, let $\cU$ be the $r$-dimensional subspace spanned
    by the columns of $\Us$, i.e.\ $\cU = \ran\Us$, and consider $K$ as a linear map on $\bR^n$.

    $K$ being diagonalizable, the restriction $K_{|\cU}$ of $K$ to its invariant subspace $\cU$ is
    also diagonalizable. Hence there exists an orthonormal basis of $\cU$ formed of eigenvectors of
    $K_{|\cU}$ and therefore of eigenvectors of $K$. Given the invariance $\sF(UQ) = \sF(U)$ for
    every orthogonal matrix $Q\in\bR^{r\times r}$, we can arbitrary set the columns of $\Us$ to be
    eigenvectors of $K$. As a consequence, the determinant is
    \begin{equation*}
        \log\det(\Us^T K \Us) = \sum_{i \in \cI} \log\lambda_i,
    \end{equation*}
    where $\cI$ is a subset of $\{1,\hdots,n\}$ such that $\#\cI=r$. The sum is maximized by
    picking the $r$ largest eigenvalues $(\lambda_i)_{i=1}^r$, and therefore a solution $\Us$ is
    given by a matrix whose columns corresponds to $r$ eigenvectors associated to the dominant
    eigenvalues.
\end{proof}

Since $\ce^{-\frac{1}{2}}\cy\ce^{-\frac{1}{2}}$ is symmetric positive definite,
Lemma~\ref{lem:maxdet} gives first that a solution to Problem~\eqref{eq:maxdetu} is given by the
matrix $U$ whose columns are the dominant eigenvectors of $\ce^{-\frac{1}{2}}\cy\ce^{-\frac{1}{2}}$.
Using the equality $V = \ce^{-\frac{1}{2}} U$, we finally find that a solution to
Problem~\eqref{eq:maxdet} is given by the matrix $V$ whose columns are $r$ dominant eigenvectors
associated to the generalized eigenvalue problem
\begin{align*}
    \cy v = \lambda \ce v, \quad \lambda \in\bR, \ v \in \bR^n.
\end{align*}